\documentclass[a4paper,reqno,11pt]{amsart}
\usepackage[T1]{fontenc}
\usepackage[utf8]{inputenc}
\usepackage{amsmath,amssymb,amsfonts,amsthm,amscd}
\usepackage{mathrsfs}
\usepackage{calligra}
\usepackage{bm}
\usepackage{bbm}

\usepackage{booktabs}   
\usepackage{tabularx}   

\usepackage{tikz}
\usetikzlibrary{arrows.meta, positioning, calc}
\usepackage{pgf}
\usepackage{tikz-cd}
\usepackage{tkz-graph}
\usepackage[all,2cell]{xy}

\usepackage{graphicx}
\usepackage{multicol}
\usepackage{enumerate}
\usepackage{comment}
\usepackage{rotating}
\usepackage[normalem]{ulem}
\usepackage{xcolor}

\usetikzlibrary{decorations.markings,arrows}

\tikzset{
  midarrow/.style={
    postaction=decorate,
    decoration={
      markings,
      mark=at position .5 with {\arrow{stealth}}
    }
  }
}

\allowdisplaybreaks

\UseAllTwocells
\SilentMatrices

\numberwithin{equation}{section}

\usepackage[
  pagebackref,
  colorlinks=true,
  linkcolor=red,
  citecolor=teal,
  urlcolor=blue,
  hypertexnames=false
]{hyperref}

\newtheorem{thm}{Theorem}[section]
\newtheorem{cor}[thm]{Corollary}
\newtheorem{lem}[thm]{Lemma}
\newtheorem{prop}[thm]{Proposition}

\newtheorem{ques}[thm]{Question}

\theoremstyle{definition}
\newtheorem{dfn}[thm]{Definition}

\newtheorem{rmk}[thm]{Remark}
\newtheorem{rmks}[thm]{Remarks}
\newtheorem{example}[thm]{Example}
\newtheorem{examples}[thm]{Examples}

\theoremstyle{remark}

\def\b{\beta}

\newcommand{\image}{\operatorname{Im}}

\newcommand{\Obj}{\operatorname{Obj}}

\newcommand{\MCE}{\operatorname{MCE}}

\newcommand{\lspan}{\operatorname{span}}

\newcommand{\gr}{\operatorname{gr}}

\def\-{\text{-}}

\newcommand{\Hom}{\operatorname{Hom}}

\newcommand{\Typ}{\operatorname{Typ}}

\newcommand{\Ker}{\operatorname{Ker}}

\newcommand{\KP}{\operatorname{KP}}

\newcommand{\N}{\mathbf{n}}
\newcommand{\M}{\mathbf{m}}
\newcommand{\bP}{\mathbf{p}}
\newcommand{\Q}{\mathbf{q}}

\newcommand{\id}{\operatorname{id}}

\begin{document}

\pagestyle{headings}

\title[Remarks on the Farthing--Muhly--Yeend inverse semigroup]{Remarks on the Farthing--Muhly--Yeend inverse semigroup of a finitely aligned higher-rank graph}

\author{Promit Mukherjee}
\address{Promit Mukherjee: Department of Mathematics, Jadavpur University, Kolkata-700032, India} \email{\url{promitmukherjeejumath@gmail.com}}
\date{November 3, 2026}



\begin{abstract}
Given a finitely aligned higher-rank $k$-graph $\Lambda$, Farthing, Muhly and Yeend (Semigroup Forum 2005) constructed an inverse semigroup $S_\Lambda$ to provide a groupoid model for the $C^*$-algebra of $\Lambda$. In this paper, we extend some of the existing results related to the semigroup structure of a graph inverse semigroup to the more general setting of the inverse semigroup $S_\Lambda$. We obtain geometric descriptions of the Green's relations on $S_\Lambda$, which effectively generalize the existing characterization of Green's relations on a graph inverse semigroup. We provide a simple and useful description of the Cuntz--Krieger ideal of the Booleanization of $S_\Lambda$ and, based on this, identify a certain $X$-to-join Booleanization of $S_\Lambda$, where the abstract and concrete Cuntz--Krieger relations agree. Besides the natural partial order, $S_\Lambda$ is also ordered by set-inclusion. This helps us to realize $S_\Lambda$ as a semigroup \emph{graded-ordered} by $\mathbb{Z}^k$ and, subsequently, gives rise to a full inverse subsemigroup $S_\Lambda^h$ of $S_\Lambda$, which is graded by $\mathbb{Z}^k$ in the usual sense. We show that the Kumjian--Pask algebra of $\Lambda$ is graded isomorphic to the quotient of the contracted semigroup algebra of $S_\Lambda^h$ by the tight ideal. In the final part of the paper, we show that a $k$-graph morphism $p:\Lambda_1\longrightarrow \Lambda_2$ that is injective and satisfies a certain path lifting property induces, in a natural way, an injective $0$-preserving semigroup homomorphism from $S_{\Lambda_1}$ to $S_{\Lambda_2}$. In view of this, we build a suitable categorical setup to conclude that the isomorphism established between the Kumjian--Pask algebra of $\Lambda$ and the tight ideal quotient of the semigroup algebra of $S_\Lambda^h$ is natural. 
\end{abstract}

\maketitle


\section{Introduction}\label{sec Intro}
Over the years, the class of inverse semigroups has become a perfect fulcrum in the development of algebras with combinatorial origin. Realizing combinatorial algebras as groupoid convolution algebras is effective in understanding the structure of these algebras. Inverse semigroups arising from the combinatorial data often help to identify the right groupoid models---in most cases, the groupoid is either the universal groupoid of the inverse semigroup or, a suitable restriction thereof, e.g., the tight groupoid \cite{Exel}. In the literature, this inverse semigroup-theoretic pathway has been used extensively to study combinatorial Cuntz--Krieger $C^*$-algebras such as graph $C^*$-algebras \cite{Paterson}, $C^*$-algebras of higher-rank graphs \cite{Farthing}, $C^*$-algebras of categories of paths \cite{Spielberg, Milan} or, more generally, $C^*$-algebras of left cancellative small categories \cite{LCSC, Ortega}, and $C^*$-algebras of self-similar graphs and $k$-graphs \cite{Pardo, Larki}.   

In \cite{Kumjian-Pask}, Kumjian and Pask defined higher-rank graphs (or $k$-graphs) as a higher-dimensional generalization of the path categories of directed graphs. They restricted attention to row-finite higher-rank graphs with no sources and studied their universal $C^*$-algebras. Raeburn, Sims and Yeend successively weakened this restriction, first by allowing sources and imposing a local convexity condition \cite{Raeburn}, and then by considering the concept of \emph{finite alignment} \cite{Yeend}, which says that any two paths can be extended in a finite number of ways to a common larger path of minimum possible degree. What they defined are known as finitely aligned higher-rank graphs. They also introduced and studied in \cite{Yeend} the $C^*$-algebra of these combinatorial objects. In \cite{Pino}, Aranda Pino et al. introduced \emph{Kumjian--Pask algebras}, which extend Leavitt path algebras of directed graphs \cite{Abrams} to higher dimensions and can be regarded as purely algebraic analogues of the $C^*$-algebras of row-finite higher-rank graphs without sources. These algebras are generalized in the setting of finitely aligned $k$-graphs by Clark and Pangalela \cite{Clark}.   

Inspired by the work of Paterson \cite{Paterson}, in order to provide a groupoid model for the $C^*$-algebra of a finitely aligned $k$-graph, Farthing, Muhly and Yeend \cite{Farthing} defined an inverse semigroup $S_\Lambda$. By considering the action of $S_\Lambda$ given by prefix exchange on the path space of $\Lambda$, they constructed the corresponding groupoid of germs. This groupoid, when restricted to the closed invariant subspace of the unit space consisting of the boundary paths, models the $C^*$-algebra of $\Lambda$. Unlike ordinary directed graphs, two paths in a general finitely aligned $k$-graph ($k\ge 2$) can have more than one minimal common extension. To accommodate this in the multiplication of the semigroup, the elements of $S_\Lambda$ were naturally chosen to be certain finite subsets of pairs of paths instead of single pairs of paths. This makes $S_\Lambda$ a rather large semigroup compared to the graph inverse semigroup of a directed graph \cite{Ash-Hall, Lawson, Meakin-Wang, Mesyan, Paterson}. Given the comparatively large size of this inverse semigroup, which reflects the complexity of the underlying higher-rank graph, it is natural to explore the extent to which the structural properties of graph inverse semigroups extend to this more general class of inverse semigroups. This is the prime motivation of this paper. 

In this paper, we first study the Green's relations on the inverse semigroup $S_\Lambda$ and obtain complete descriptions of the relations $\mathscr{L}$, $\mathscr{R}$, $\mathscr{H}$ and $\mathscr{D}$ (Propositions \ref{pro the relations L,R,H} and \ref{pro necessary and sufficient condition for D}). For the $\mathscr{J}$ relation, we restrict ourselves to certain elements of $S_\Lambda$, which we call \emph{divergent}, and obtain necessary and sufficient conditions to determine when two divergent elements are $\mathscr{J}$-related (Proposition \ref{pro partial orders on equivalence classes J} and Corollary \ref{cor J relation on divergent sets}). The graph inverse semigroup $S(\Gamma)$ of a directed graph $\Gamma$ embeds inside the Farthing--Muhly--Yeend inverse semigroup of the corresponding path category $\Gamma^*$. In view of this, we note that our characterizations of the Green's relations on $S_\Lambda$ perfectly extend the existing characterizations of Green's relations on a graph inverse semigroup obtained by Mesyan and Mitchell \cite{Mesyan}.

We then explore some fundamental concepts of $S_\Lambda$ related to idempotents. In particular, we show that $S_\Lambda$ is a fundamental inverse semigroup, and provide geometric descriptions of compatibility, orthogonality and the natural partial order using the combinatorics of paths in $\Lambda$ (Propositions \ref{pro S_Lambda is fundamental}, \ref{pro describing compatibility} and \ref{pro natural partial order on S_Lambda}). The heart of the Cuntz--Krieger algebras of combinatorial objects (originated in the works of Cuntz and Krieger \cite{CK}) comprises certain relations determined by the combinatorial data. These are formally known as the \emph{Cuntz--Krieger relations}. In \cite{LawsonCK}, Lawson and Vdovina, using the notion of covers of elements, abstractified these relations in the setting of inverse semigroups. They explicitly verified that abstract Cuntz--Krieger relations agree with the concrete ones in the case of tight Booleanization of polycyclic inverse monoids or, more generally, of graph inverse semigroups of directed graphs. Motivated to understand the relationship between abstract and concrete Cuntz--Krieger relations in the context of the Farthing--Muhly--Yeend inverse semigroup $S_\Lambda$, we obtain a simpler description of the Cuntz--Krieger ideal of the universal Booleanization $\mathsf{B}(S_\Lambda)$ of $S_\Lambda$ in terms of generators and relations (Theorem \ref{th description of Cuntz-Krieger ideal}). More precisely, we show that any arbitrary cover-to-join relation can be deduced if we know the relations originating from $(i)$ the internal covers of elements of $S_\Lambda$ given by partitions and $(ii)$ the covers of vertices of $\Lambda$ given by paths in exhaustive subsets. This helps us to conclude that if we perform an $X$-to-join Booleanization \cite{GannaX} of $S_\Lambda$, with the base set $X$ determined by the internal cover relations of $(i)$ (Remark \ref{rem abstract CK relations and X-to-join}), then the remaining abstract Cuntz--Krieger relations are precisely the concrete Cuntz--Krieger relations that define the higher-rank Cuntz--Krieger algebras.

Being a subset of the power set of $\Lambda$, the Farthing--Muhly--Yeend inverse semigroup $S_\Lambda$ is naturally equipped with $\subseteq$, the partial order of set-inclusion. The ordered structure of $S_\Lambda$ with respect to $\subseteq$ appears to be an important feature of this semigroup in many respects. For example, the order $\subseteq$ is included in the natural partial order $\le$ on $S_\Lambda$, and the existence of joins with respect to $\subseteq$ ensures the existence of joins with respect to $\le$; in fact, both joins coincide (Proposition \ref{pro natural partial order and set-inclusion}). The join-closed ideals of the ordered semigroup $(S_\Lambda,\subseteq)$ are in inclusion-preserving one-to-one correspondence with the hereditary subsets of $\Lambda^0$ (Theorem \ref{th structure thm for U-closed ideals}). Most importantly, in view of this order of set-inclusion, the elements of $S_\Lambda$ can be decomposed based on the natural $\mathbb{Z}^k$-grading induced by the degree functor $d:\Lambda\longrightarrow \mathbb{N}^k$. This motivates us to formulate the notion of a \emph{graded-ordered semigroup}. Roughly speaking, it is an ordered semigroup with a graded ring-like structure (sum of elements replaced by existential join). We provide some interesting examples of graded-ordered semigroups (Examples \ref{ex graded-ordered semigroup}), which include the multiplicative reduct of an additively cancellative zero-sum free group graded semiring as well as any inverse semigroup graded in the sense of Hazrat and Mesyan \cite{HazMes}. It follows that the inverse semigroup $S_\Lambda$ is naturally graded-ordered by $\mathbb{Z}^k$ and, consequently, the union of all the homogeneous components $S_\Lambda^h$ is a $\mathbb{Z}^k$-graded full inverse subsemigroup of $S_\Lambda$. Using the graded inverse semigroup $S_\Lambda^h$, we are able to prove that the Kumjian--Pask algebra of $\Lambda$ is graded isomorphic to the tight ideal quotient of a contracted inverse semigroup algebra, namely $\KP_\mathsf{k}(\Lambda)\cong^{\gr} \mathsf{k}_0S_\Lambda^h/{\mathcal{T}_\mathsf{k}(S_\Lambda^h)}$, where $\mathcal{T}_\mathsf{k}(S_\Lambda^h)$ is the tight ideal of $S_\Lambda^h$ (Theorem \ref{th the relationship between S_Lambda and KPA}). This is possibly new in the literature of Kumjian--Pask algebras and extends to higher dimensions, the natural realization of a Leavitt path algebra as the contracted semigroup algebra of the tight Booleanization of the corresponding graph inverse semigroup.

The final theme of this paper is to set up a categorical framework in which one can interpret the isomorphism between $\KP_\mathsf{k}(\Lambda)$ and $\mathsf{k}_0S_\Lambda^h/{\mathcal{T}_\mathsf{k}(S_\Lambda^h)}$. Given two finitely aligned $k$-graphs $\Lambda_1$ and $\Lambda_2$, we ask the question: when can a $k$-graph morphism $p:\Lambda_1\longrightarrow \Lambda_2$ be extended naturally to a semigroup homomorphism between the corresponding Farthing--Muhly--Yeend inverse semigroups? In contrast to the graph inverse semigroup case (see \cite[Theorem 20]{Mesyan}), we here observe that imposition of the injectivity condition on $p$ is, in general, not enough to induce a semigroup homomorphism (Example \ref{ex orthogonal map does not produce semigroup hom}). Together with injectivity, if $p$ has the \emph{$r$-path lifting property}, then it naturally induces an injective $0$-preserving semigroup homomorphism from $S_{\Lambda_1}$ to $S_{\Lambda_2}$ (Theorem \ref{th answer to the query}). We consider the category of finitely aligned $k$-graphs together with morphisms having these properties and exhibit two functors to the category of $\mathbb{Z}^k$-graded algebras. One is the Kumjian--Pask algebra functor (Theorem \ref{th KP functor}) and the other is a composition of three functors defined between several intermediate categories, which, given a $k$-graph $\Lambda$, first produces the graded-ordered inverse semigroup $S_\Lambda$, then the graded inverse semigroup $S_\Lambda^h$, and finally the graded algebra $\mathsf{k}_0S_\Lambda^h/{\mathcal{T}_\mathsf{k}(S_\Lambda^h)}$. Our final remark is that these two functors are naturally isomorphic (Remark \ref{rem the natural isomorphism}). Thus, the obtained relationship between discrete higher-rank graph algebras and the inverse semigroups of higher-rank graphs is not merely an objectwise correspondence but is also stable under the dynamics created by appropriate $k$-graph morphisms. 

In a recent discussion concerning the character space of the semilattice of idempotents of $S_\Lambda$, Larki drew our attention to the work of Ortega and Pardo \cite{Ortega}. In \cite{Ortega}, the authors defined some inverse semigroups to model $C^*$-algebras of finitely aligned left cancellative small categories as tight $C^*$-algebras. Since finitely aligned higher-rank graphs form a special class of such categories, some of the structural features of the Farthing--Muhly--Yeend inverse semigroup, investigated in this paper, have natural counterparts in the framework of \cite{Ortega}. However, by exploiting the combinatorial and geometric structure of higher-rank graphs, we here provide an explicit analysis of the Farthing--Muhly--Yeend inverse semigroup. To the best of our knowledge, such an analysis, including the graded-order-theoretic approach to this inverse semigroup and the functorial study of the relationship with Kumjian--Pask algebras, has not appeared previously in the literature.

We now briefly describe the organization of the paper. In Section \ref{sec prelim}, we recall the necessary background on inverse semigroups and higher-rank graphs and provide an overview of the Farthing--Muhly--Yeend inverse semigroup of a $k$-graph. In Section \ref{sec Green's relations}, we obtain the characterizations of the Green's relations on $S_\Lambda$. Section \ref{sec partial order, covers, joins} examines further structural facets of $S_\Lambda$, including fundamentality, compatibility of elements, the natural partial order, covers and the Cuntz–-Krieger ideal of its Booleanization. In Section \ref{sec Relation to KP-algebra}, we study the ordered structure of $S_\Lambda$ with respect to the partial order of set-inclusion, establish this as a graded-ordered inverse semigroup and, based on the results of the previous section, prove a graded isomorphism between the Kumjian--Pask algebra of $\Lambda$ and the tight ideal quotient of the contracted semigroup algebra of the graded core of $S_\Lambda$. Section \ref{sec extending morphisms} is devoted to the formalization of several categories and functors in view of the semigroup homomorphisms induced by suitable $k$-graph morphisms, and to the application of these to establish the naturality of the graded isomorphism obtained in the previous section.   

\section{Fundamental concepts and the main object of study}\label{sec prelim}
In this section, we collect the prerequisite definitions and notions related to inverse semigroups and higher-rank graphs, as well as recall our main object of study: the inverse semigroup associated with a higher-rank graph.  

\subsection{Semigroups}\label{ssec semigroups} 
A nonempty set $S$ together with an associative binary operation is called a \emph{semigroup}. A semigroup is called a \emph{monoid} if it contains an identity with respect to the binary operation. A \emph{zero} in a semigroup $S$ is an element $0$ such that $0s=s0=0$ for all $s\in S$. If a zero exists, then it is obviously unique. A nonempty subset $T$ of a semigroup $S$ is called a \emph{subsemigroup} if $T$ is closed under the operation of $S$. A subsemigroup is termed as \emph{full} if it contains the set of idempotents $\mathcal{E}(S)$ of $S$. A subsemigroup $J$ of $S$ is called a \emph{left (right, or two-sided) ideal} if it absorbs elements of $S$ from the left (the right or both sides). For any $x\in S$, $S^1 x:=Sx\cup \{x\}$ is the principal left ideal of $S$ generated by $x$. Principal right and two-sided ideals are analogously denoted by $xS^1$ and $S^1 x S^1$, respectively.

One can consider the following five equivalence relations on a semigroup $S$, which are collectively known as the \emph{Green's relations}.

$(i)$ $\mathscr{L}$: $x~\mathscr{L}~y\Leftrightarrow S^1 x=S^1 y$;

$(ii)$ $\mathscr{R}$: $x~\mathscr{R}~y\Leftrightarrow xS^1=yS^1$;

$(iii)$ $\mathscr{H}:=\mathscr{L}\cap \mathscr{R}$;

$(iv)$ $\mathscr{D}:=\mathscr{L}\vee \mathscr{R}$;

$(v)$ $\mathscr{J}$: $x~\mathscr{J}~y\Leftrightarrow S^1 x S^1=S^1 y S^1$.

The relation $\le_{\mathscr{L}}$ defined by $[x]_{\mathscr{L}}\le_{\mathscr{L}}[y]_{\mathscr{L}}$ if and only if $S^1x\subseteq S^1y$ is a partial order on the set of equivalence classes of $\mathscr{L}$. Similarly, one has the partial orders $\le_{\mathscr{R}}$ and $\le_{\mathscr{J}}$.

For two semigroups $S$ and $T$, a map $\phi:S\longrightarrow T$ is said to be a \emph{homomorphism} if $\phi(s_1s_2)=\phi(s_1)\phi(s_2)$ for all $s_1,s_2\in S$. A homomorphism $\phi:S\longrightarrow T$ between semigroups with zeros is called \emph{$0$-preserving} if $\phi(0)=0$. 

If $S$ is a semigroup and $\leqslant$ is a partial order on $S$ such that for any $s,x,y\in S$ with $x\leqslant y$ we have $sx\leqslant sy$ and $xs\leqslant ys$, then the pair $(S,\leqslant)$ is called an \emph{ordered semigroup}. An \emph{order ideal} of an ordered semigroup $(S,\leqslant)$ is an ideal $J$ that is also a down-set with respect to the partial order $\leqslant$. An order ideal $J$ is called \emph{join-closed} if $J$ is closed under existing joins, i.e., for any two elements $x,y\in J$ admitting a join in $S$ with respect to $\leqslant$, $x\vee y\in J$. Two elements $s,t$ in an ordered semigroup $(S,\leqslant)$ are called \emph{comparable} if either $s\leqslant t$ or $t\leqslant s$. If neither of these happens, $s$ and $t$ are called \emph{incomparable}, in which case we write $s~||~t$. For $s\in S$, $s^{\downarrow}$ (resp., $s^{\uparrow}$) denotes the principal down-set (resp., principal up-set) of $s$ with respect to the partial order $\leqslant$.

For any semigroup $S$ and a commutative unital ring $\mathsf{k}$, the \emph{semigroup algebra} of $S$ is denoted by $\mathsf{k}S$ and is defined as the free associative $\mathsf{k}$-algebra generated by the underlying set of $S$ subject to the relations given by the binary operation on $S$. If a semigroup $S$ has a $0$, then the quotient of $\mathsf{k}S$ by the ideal $J_0:=\{a0~|~a\in \mathsf{k}\}$ is denoted by $\mathsf{k}_0S$ and is called the \emph{contracted semigroup algebra} of $S$. For any $s\in S$, we denote the image of $s$ in $\mathsf{k}_0S$ by $\overline{s}$ with the understanding that $\overline{0}=0$.

If $S$ is a semigroup with $0$ and $G$ is an abelian group, then $S$ is called $G$-graded (see \cite[Definition 2.1]{HazMes}) if there is a family $\{S_x~|~x\in G\}$ of subsets of $S$ such that 

$(i)$ $S=\displaystyle{\bigcup_{x\in G}}~S_x$;

$(ii)$ $S_x S_y\subseteq S_{xy}$ for all $x,y\in G$;

$(iii)$ $S_x\cap S_y=\{0\}$ for any two distinct $x,y\in G$.

For any $0\neq s\in S$, we write $\deg(s)=x$ if $s\in S_x$. From $(ii)$ it follows that $\deg(e)=\epsilon$ for all $e\in \mathcal{E}(S)$, where $\epsilon$ is the identity of $G$. A homomorphism $\phi:S\longrightarrow T$ between $G$-graded semigroups $S$ and $T$ is called a \emph{graded homomorphism} if $\phi(S_x)\subseteq T_x$ for all $x\in G$. 

When $S$ is a $G$-graded semigroup with $0$, it is easy to see that the grading of $S$ induces a natural $G$-grading on the contracted semigroup algebra $\mathsf{k}_0S$, where the homogeneous components are given by \[(\mathsf{k}_0S)_x:=\lspan_\mathsf{k}\left(\{\overline{s}~|~s\in S_x\}\right)\] for all $x\in G$. 

\subsection{Inverse semigroups}\label{ssec invsem}
For a detailed and comprehensive exposition of the theory of inverse semigroups, we refer the reader to the wonderful book of Lawson \cite{Lawson}. 

A semigroup $S$ is called an \emph{inverse semigroup} if for each $s\in S$, there is a unique element denoted by $s^{-1}\in S$ such that $ss^{-1}s=s$ and $s^{-1}ss^{-1}=s^{-1}$. The unary operation $s\longmapsto s^{-1}$ is often called the \emph{inversion}. An \emph{inverse subsemigroup} of an inverse semigroup is a subsemigroup which is closed under the inversion. 

The Green's relations $\mathscr{L}$, $\mathscr{R}$ and $\mathscr{D}$ on an inverse semigroup can be simplified substantially as described in the following lemma.

\begin{lem}\label{lem L and R on inverse semigroups}
Let $S$ be an inverse semigroup and $x,y\in S$. Then 

$(i)$ $[x]_{\mathscr{L}} \le_{\mathscr{L}} [y]_{\mathscr{L}}$ if and only if $x^{-1}x\le y^{-1}y$,

$(ii)$ $[x]_{\mathscr{R}} \le_{\mathscr{R}} [y]_{\mathscr{R}}$ if and only if $xx^{-1}\le yy^{-1}$,

$(iii)$ $x~\mathscr{L}~y$ if and only if $x^{-1}x=y^{-1}y$,

$(iv)$ $x~\mathscr{R}~y$ if and only if $xx^{-1}=yy^{-1}$.

$(v)$ $x~\mathscr{D}~y$ if and only if there exists $z\in S$ such that $x^{-1}x=z^{-1}z$ and $zz^{-1}=yy^{-1}$.
\end{lem}
Using the above lemma, it is easy to observe that if $T$ is an inverse subsemigroup of an inverse semigroup $S$ and $x,y\in T$, then $x~\mathscr{L}~y$ (resp., $x~\mathscr{R}~y$) in $T$ if and only if $x~\mathscr{L}~y$ (resp., $x~\mathscr{R}~y$) in $S$. 

In an inverse semigroup $S$, $\mathcal{E}(S)$ forms a commutative inverse subsemigroup of $S$, which is for many reasons, the most important subsemigroup related to the structure of $S$. The idempotents give rise to a partial order $\le$ on $S$ which is defined as follows. For any two $s,t\in S$, $s\le t$ if and only if any of the following equivalent conditions hold:

$(i)$ $s=te$ for some $e\in \mathcal{E}(S)$;

$(ii)$ $s=ft$ for some $f\in \mathcal{E}(S)$;

$(iii)$ $s=ss^{-1}t$;

$(iv)$ $s=ts^{-1}s$.

The relation $\le$ is referred to as the \emph{natural partial order} on $S$. The natural partial order respects the multiplication and inversion of $S$, consequently, $(S,\le)$ is an ordered inverse semigroup. The restriction of $\le$ to $\mathcal{E}(S)$ has a simpler description: for $e,f\in \mathcal{E}(S)$, $e\le f$ if and only if $ef=fe=e$. In view of this, $\mathcal{E}(S)$ becomes a \emph{lower-semilattice}, where the meet of two idempotents $e,f$ is simply their product $ef$. Two elements $s,t$ of $S$ are called 

$(i)$ \emph{compatible} if both $s^{-1}t$ and $st^{-1}$ are idempotents;

$(ii)$ \emph{orthogonal} if $s^{-1}t=st^{-1}=0$. 

If $s,t$ are compatible (resp., orthogonal), then we write $s\sim t$ (resp., $s\perp t$). A \emph{compatible} subset of $S$ is a subset in which any two elements are compatible. The inverse semigroup $S$ is called \emph{Boolean} \cite{LawBool} if $(\mathcal{E}(S),\le)$ is a generalized Boolean algebra and any two compatible elements in $S$ have a join with respect to $\le$.

If $S$ is an inverse semigroup, then a \emph{character} of $\mathcal{E}(S)$ is a nonzero $0$-preserving semilattice homomorphism $\chi:\mathcal{E}(S)\longrightarrow \textbf{2}$, where \textbf{2} is the two-element Boolean algebra. The space of all characters of $\mathcal{E}(S)$ denoted by $\widehat{\mathcal{E}(S)}$ forms a locally compact Hausdorff topological space under the subspace topology inherited from the product topology on $\textbf{2}^{\mathcal{E}(S)}$, or equivalently, the topology of pointwise convergence. There is a natural action (see \cite{Exel}) (called the \emph{spectral action}) $\varphi$ of $S$ on the space $\widehat{\mathcal{E}(S)}$ defined by $\varphi_s(\chi)(e):=\chi(s^{-1}e s)$ for all $s\in S$, $e\in \mathcal{E}(S)$ and $\chi\in \widehat{\mathcal{E}(S)}$ with $\chi(s^{-1}s)=1$. The groupoid of germs \cite[Proposition 4.17]{Exel} of the spectral action is denoted by $\mathcal{G}(S)$ and is called the \emph{universal groupoid} of $S$. This is an ample groupoid. The dual Boolean inverse semigroup under the non-commutative Stone duality \cite{LawNC} is denoted by $\mathsf{B}(S)$ and is called the \emph{Booleanization} of $S$. This can also be obtained via a universal construction starting from $S$ (see \cite{LawBool} for details). 

We now recall the notion of covers in an inverse semigroup $S$ with $0$. A \emph{cover} of $s\in S$ is a subset $\mathcal{Z}\subseteq s^{\downarrow}$ such that for any $0\neq t\le s$, there is a $z\in Z$ such that $t\wedge z\neq 0$ (since $t,z\le s$, we have $t\sim z$ and hence their meet exists by \cite[Lemma 1.4.11]{Lawson}). For two inverse semigroups $S_1,S_2$ with zeros, Donsig and Milan \cite{Milan} defined a homomorphism $\phi:S_1\longrightarrow S_2$ to be \emph{cover-to-join} if for any cover $\mathcal{Z}$ of $s\in S_1$, the elements $\phi(z)$, $z\in \mathcal{Z}$ admit a join in $S_2$ and $\phi(s)=\displaystyle{\bigvee_{z\in \mathcal{Z}}}\phi(z)$. 

Suppose $S$ is an inverse semigroup with $0$ and $\mathsf{B}(S)$ is the Booleanization of $S$. For any element $a\in S$, a finite cover $\mathcal{Z}$ of $a$ is a compatible set in $S$ and hence remains so when viewed as a subset of $\mathsf{B}(S)$ under the canonical map from $S$ to $\mathsf{B}(S)$. By the property of Boolean inverse semigroups, $\displaystyle{\bigvee_{z\in \mathcal{Z}}}~z$ exists in $\mathsf{B}(S)$. Moreover, since $\displaystyle{\bigvee_{z\in \mathcal{Z}}}~z\le a$, we can consider the relative complement $a\setminus \left(\displaystyle{\bigvee_{z\in \mathcal{Z}}}~z\right)$. In \cite{LawsonCK}, Lawson and Vdovina named the relations $a\setminus \left(\displaystyle{\bigvee_{z\in \mathcal{Z}}}~z\right)=0$ or equivalently, 
\begin{equation}\label{eq ACK}
a=\displaystyle{\bigvee_{z\in \mathcal{Z}}}~z
\end{equation}
the \emph{abstract Cuntz--Krieger relations} on $\mathsf{B}(S)$. The additive ideal generated by the relations (\ref{eq ACK}) is called the \emph{Cuntz--Krieger ideal} of $\mathsf{B}(S)$ in \cite{LawsonCK}. The quotient of $\mathsf{B}(S)$ by the Cuntz--Krieger ideal is known as the \emph{tight Booleanization} of $S$ and is denoted by $\mathsf{B}_{tight}(S)$ (see \cite{GannaX,Ganna}). Equivalently, $\mathsf{B}_{tight}(S)$ is the dual Boolean inverse semigroup of the \emph{tight groupoid} $\mathcal{G}_{tight}(S)$ of $S$, which is defined as the restriction of the universal groupoid $\mathcal{G}(S)$ to the closed $\varphi$-invariant subspace $\widehat{\mathcal{E}(S)}_{tight}$ of $\widehat{\mathcal{E}(S)}$ consisting of the \emph{tight} characters \cite{Exel}. By \cite[Theorem 2.2]{Milan}, a tight character of $\mathcal{E}(S)$ is a character of $\mathcal{E}(S)$ that is also cover-to-join.  

\subsection{Higher-rank graphs}\label{ssec higher-rank graphs}
Throughout this section, $\mathbb{N}$ stands for the set of nonnegative integers. Given a positive integer $k$, a \emph{higher-rank graph of rank $k$} or a $k$-graph is defined as a countably small category $\Lambda$ together with a functor $d:\Lambda\longrightarrow \mathbb{N}^k$ ($\mathbb{N}^k$ is viewed as a small category with one object) such that for any $\lambda\in \Lambda$ and $\N,\M\in \mathbb{N}^k$ with $d(\lambda)=\N+\M$, there are unique $\alpha,\beta\in \Lambda$ satisfying $d(\alpha)=\N$, $d(\beta)=\M$ and $\lambda=\alpha\beta$. This specific property of the functor $d$ is called the \emph{unique factorization property}.

Suppose $\Lambda$ is a $k$-graph. The functor $d$ is called the \emph{degree} functor. The domain and the codomain maps are respectively denoted by $s$ and $r$. The elements of $\Lambda$ are usually called \emph{paths}. For any path $\lambda$, $r(\lambda)$, $s(\lambda)$ and $d(\lambda)$ are respectively called the \emph{range}, the \emph{source} and the \emph{degree} of $\lambda$. By considering the \emph{arrows only} approach to categories, we can view $\Obj(\Lambda)$ as a subset of $\Lambda$ and then $\Obj(\Lambda)=\{\lambda\in \Lambda~|~r(\lambda)=\lambda=s(\lambda)\}$. For any $\N\in \mathbb{N}^k$, we denote $\Lambda^\N:=d^{-1}(\N)$. From the unique factorization property, it follows that $\Lambda^0=\Obj(\Lambda)$. 

We fix some elementary notations which will be used frequently in the paper. For $v\in \Lambda^0$, \[v\Lambda:=\{\lambda\in \Lambda~|~r(\lambda)=v\},~~\Lambda v:=\{\lambda\in \Lambda~|~s(\lambda)=v\}.\] More generally, for any two subsets $X,Y$ of $\Lambda$, \[XY:=\{\lambda\mu~|~\lambda\in X,\mu\in Y~\text{and}~s(\lambda)=r(\mu)\}.\] If $X=\{\lambda\}$ (resp., $Y=\{\mu\}$), then we simply write $XY$ as $\lambda Y$ (resp., $X\mu$). The above rule defines a semigroup structure on $\mathcal{P}(\Lambda)$. The $k$-graph $\Lambda$ is called \emph{row-finite} if the set $v\Lambda^\N$ is finite for all $v\in \Lambda^0$ and $\N\in \mathbb{N}^k$.

If $\Lambda_1,\Lambda_2$ are $k$-graphs with respective degree functors $d_1$ and $d_2$, then a \emph{$k$-graph morphism} from $\Lambda_1$ to $\Lambda_2$ is a functor $p:\Lambda_1\longrightarrow \Lambda_2$ such that $d_2(p(\lambda))=d_1(\lambda)$ for all $\lambda\in \Lambda_1$. 

Suppose $\Gamma=(\Gamma^0,\Gamma^1,r_\Gamma,s_\Gamma)$ is a directed graph, where $\Gamma^0,\Gamma^1$ are respectively the sets of vertices and edges, and $r_\Gamma$, $s_\Gamma$ are respectively the range and source maps. The set of all finite paths $\Gamma^*$ forms a small category (called the \emph{path category}), where $r(\alpha):=s_\Gamma(\alpha)$, $s(\alpha):=r_\Gamma(\alpha)$ and the composition is the usual concatenation of paths. The definition of a $k$-graph is modelled on the path categories and it turns out that these are exactly the higher-rank graphs of rank $1$. For $k\ge 2$, the geometric structure of a $k$-graph is far more complex. To know more about the concepts of general $k$-graphs (including the notions of $1$-skeleton, infinite paths, geometric properties like cofinality, aperiodicity and the $C^*$-algebra of higher-rank graphs), the reader may go through \cite{kNotes}. 

We now recall a certain class of $k$-graphs, which generalize both row-finite $k$-graphs and $1$-graphs. To do this, we need to be accustomed with some notations. In what follows, we view $\mathbb{N}^k$ as an ordered monoid with respect to the coordinatewise ordering. 
\begin{itemize}
    \item[$\bullet$] Let $\lambda\in \Lambda$ and $0\le \M \le \N\le d(\lambda)$. By $\lambda(\M,\N)$, we denote the unique path of $\Lambda$ (exists by the unique factorization property) such that $\lambda=\alpha \lambda(\M,\N)\beta$ for some $\alpha\in \Lambda^{\M}$ and $\beta\in \Lambda^{d(\lambda)-\N}$.

    \item[$\bullet$] For any $\lambda,\mu\in \Lambda$, \[\MCE(\lambda,\mu):=\{\tau\in \Lambda~|~\tau(0,d(\lambda))=\lambda,\tau(0,d(\mu))=\mu~\text{and}~d(\tau)=d(\lambda)\vee d(\mu)\}.\] 

    \item[$\bullet$] For any $\lambda,\mu\in \Lambda$, \[\Lambda^{\min}(\lambda,\mu):=\{(\alpha,\beta)\in \Lambda\times \Lambda~|~\lambda\alpha=\mu\beta~\text{and}~\lambda\alpha\in \MCE(\lambda,\mu)\}.\]

    \item[$\bullet$] For $X,Y\subseteq \Lambda$, \[\Lambda^{\min}(X,Y):=\displaystyle{\bigcup_{\lambda\in X,\mu\in Y}}\Lambda^{\min}(\lambda,\mu).\]
\end{itemize}
Define a relation $\leqslant$ on $\Lambda$ as follows: $\lambda\leqslant \tau$ if and only if $\tau(0,d(\lambda))=\lambda$ ($cf.$ \cite[Definition 1.2]{Ortega}). Then $\leqslant$ is easily seen to be a partial order on $\Lambda$. A \emph{minimal common extension} of $\lambda$ and $\mu$ is a minimal element of $\lambda^{\uparrow}\cap \mu^{\uparrow}$ with respect to $\leqslant$. It is not hard to see that a minimal common extension of $\lambda$ and $\mu$ has degree $d(\lambda)\vee d(\mu)$. In view of this $\MCE(\lambda,\mu)$ is the set of all minimal common extensions of $\lambda$ and $\mu$. The $k$-graph $\Lambda$ is called \emph{finitely aligned} if $\MCE(\lambda,\mu)$ is finite (or equivalently, $|\Lambda^{\min}(\lambda,\mu)|<\infty$) for all $\lambda,\mu\in \Lambda$.

If $\Gamma$ is a directed graph, then $\lambda,\mu\in \Gamma^*$ have a minimal common extension if and only if they are comparable under the partial order $\leqslant$ defined in the previous paragraph. In this case, the larger path is the minimal common extension. It follows that $|\MCE(\lambda,\mu)|\le 1$ for all $\lambda,\mu\in \Gamma^*$, and hence the path category $\Gamma^*$ is finitely aligned as a $1$-graph.  

Given a finitely aligned $k$-graph $\Lambda$, Clark and Pangalela \cite{Clark} defined the Kumjian--Pask algebra of $\Lambda$, which is the discrete analogue of the $C^*$-algebra of $\Lambda$ \cite{Yeend}. We end this subsection by recalling its definition.

For $v\in \Lambda^0$, a subset $X$ of $v\Lambda$ is called \emph{exhaustive} if for any $\lambda\in v\Lambda$, there is a $\tau\in X$ such that $\Lambda^{\min}(\lambda,\tau)\neq \emptyset$. For any $\lambda\in \Lambda$, fix a formal symbol $\lambda^*$ (which is called the \emph{ghost path} of $\lambda$) and set $v^*=v$ for all $v\in \Lambda^0$. Extend $r$ and $s$ to ghost path as follows: $r(\lambda^*):=s(\lambda)$ and $s(\lambda^*):=r(\lambda)$. For $\lambda,\mu\in \Lambda$ with $s(\lambda)=r(\mu)$, define $\mu^*\lambda^*:=(\lambda\mu)^*$. 

Suppose $\mathsf{k}$ is a commutative unital ring and $A$ is a $\mathsf{k}$-algebra. Then a \emph{Kumjian--Pask $\Lambda$ family} (see \cite[Definition 3.1]{Clark}) in $A$ is a collection of elements $\{s_\lambda,s_{\mu^*}\in A~|~\lambda,\mu\in \Lambda\}$, which satisfies the following conditions:
\begin{itemize}
\item[(KP1)] $s_u s_v=\delta_{u,v}s_u$ for all $u,v\in \Lambda^0$;

\item[(KP2)] $s_\lambda s_\mu=s_{\lambda\mu}$ and $s_{\mu^*} s_{\lambda^*}=s_{(\lambda\mu)^*}$ for all $\lambda,\mu\in \Lambda$ with $s(\lambda)=r(\mu)$;

\item[(KP3)] $s_{\lambda^*}s_\mu=\displaystyle{\sum_{(\alpha,\beta)\in \Lambda^{\min}(\lambda,\mu)}}s_\alpha s_{\beta^*}$ for all $\lambda,\mu\in \Lambda$;

\item[(KP4)] $\displaystyle{\prod_{\tau\in X}}\left(s_v-s_\tau s_{\tau^*}\right)=0$ for all $v\in \Lambda^0$ and finite exhaustive subset $X$ of $v\Lambda$. 
\end{itemize}
The universal $\mathsf{k}$-algebra generated by a Kumjian--Pask $\Lambda$ family $\{s_\lambda,s_{\mu^*}~|~\lambda,\mu\in \Lambda\}$ (see the proof of \cite[Theorem 3.7]{Clark} for the construction) is called the \emph{Kumjian--Pask algebra} of $\Lambda$ and is denoted by $\KP_\mathsf{k}(\Lambda)$. This algebra is naturally $\mathbb{Z}^k$-graded, where the homogeneous components are given by \[\KP_\mathsf{k}(\Lambda)_\N:=\lspan_\mathsf{k}\left(\{s_\lambda s_{\mu^*}~|~\lambda,\mu\in \Lambda~\text{and}~d(\lambda)-d(\mu)=\N\}\right)\] for all $\N\in \mathbb{Z}^k$.

\subsection{Farthing--Muhly--Yeend inverse semigroup of a $k$-graph}\label{ssec The FMY inverse semigroup}
We now describe the main object of study of this paper. Fix a finitely aligned $k$-graph $\Lambda$. We write \[\Lambda *_s \Lambda:=\{(\lambda,\mu)\in \Lambda\times \Lambda~|~s(\lambda)=s(\mu)\}.\] Two pairs $(\lambda,\mu)$ and $(\xi,\eta)$ in $\Lambda *_s \Lambda$ are said to be \emph{componentwise orthogonal} if $\Lambda^{\min}(\lambda,\xi)=\emptyset=\Lambda^{\min}(\mu,\eta)$. In this case, we write $(\lambda,\mu) \perp (\xi,\eta)$. 

We now recall the definition of the inverse semigroup associated with $\Lambda$. 

\begin{dfn}\label{def the inverse semigroup S_Lambda}
(\cite[Definition 4.1]{Farthing}) Let $\Lambda$ be a finitely aligned $k$-graph and \[S_\Lambda:=\{F\subseteq \Lambda *_s \Lambda~|~|F|< \infty~\text{and}~(\lambda,\mu)\perp (\xi,\eta)~\text{for distinct pairs}~ (\lambda,\mu),(\xi,\eta)\in F\}.\] For $F,G\in S_\Lambda$, define \[FG:=\displaystyle{\bigcup_{\substack{(\lambda,\mu)\in F\\(\xi,\eta)\in G}}}\left\{(\lambda \alpha,\eta\beta)~|~(\alpha,\beta)\in \Lambda^{\min}(\mu,\xi)\right\}\] and \[F^{-1}:=\left\{(\mu,\lambda)\in \Lambda *_s \Lambda~|~(\lambda,\mu)\in F\right\}.\] Then $S_\Lambda$ is an inverse semigroup with respect to the above operations (see \cite[Proposition 4.3]{Farthing}). We call $S_\Lambda$ the \emph{Farthing--Muhly--Yeend inverse semigroup} of $\Lambda$. 
\end{dfn}

Note that $\emptyset$ is vacuously in $S_\Lambda$ and it plays the role of an (absorbing) zero. If $(\lambda,\mu),(\xi,\eta)\in \Lambda *_s \Lambda$, then for the ease of writing, we set \[(\lambda,\mu)\cdot (\xi,\eta):=\{(\lambda \alpha,\eta\beta)~|~(\alpha,\beta)\in \Lambda^{\min}(\mu,\xi)\}.\] In view of this notation, 
\begin{equation}\label{eq the multiplication}
FG=\displaystyle{\bigcup_{\substack{(\lambda,\mu)\in F\\(\xi,\eta)\in G}}}(\lambda,\mu)\cdot (\xi,\eta).
\end{equation}

Suppose $(\lambda,\mu),(\xi,\eta)\in \Lambda *_s \Lambda$. It is interesting to note that $\{(\lambda,\mu)\}\perp \{(\xi,\eta)\}$ as elements of $S_\Lambda$ if and only if $(\lambda,\mu)\perp (\xi,\eta)$. 

Viewing any subset $F\subseteq \Lambda *_s \Lambda$ as a relation on $\Lambda$, we can define the \emph{domain} and \emph{range} of $F$ in the usual way: \[\mathsf{D}(F):=\{\lambda\in \Lambda~|~(\lambda,\mu)\in F~\text{for some}~\mu\in \Lambda\},~~\mathsf{R}(F):=\{\mu\in \Lambda~|~(\lambda,\mu)\in F~\text{for some}~\lambda\in \Lambda\}.\] Observe that the map $b_F$, sending each element $\lambda\in \mathsf{D}(F)$ to the unique element $\mu\in \mathsf{R}(F)$ such that $(\lambda,\mu)\in F$, is a bijection from $\mathsf{D}(F)$ to $\mathsf{R}(F)$. For any $F\in S_\Lambda$, $\mathsf{D}(F^{-1})=\mathsf{R}(F)$ and $\mathsf{R}(F^{-1})=\mathsf{D}(F)$. Hence, we have $b_{F^{-1}}=b_F^{-1}$. 

We call a subset $X$ of $\Lambda$ \emph{pairwise orthogonal} if for any two distinct $\lambda,\mu\in X$, we have $\Lambda^{\min}(\lambda,\mu)=\emptyset$. It is clear that for any $F\in S_\Lambda$, both $\mathsf{D}(F)$ and $\mathsf{R}(F)$ are pairwise orthogonal subsets of $\Lambda$.   

It is easy to identify the idempotents in $S_\Lambda$. 
\begin{lem}\label{lem idempotents of the semigroup}
\textup{($cf.$ \cite[Corollary 3.5]{Larki})} Let $\Lambda$ be a finitely aligned $k$-graph. Then $\mathcal{E}(S_\Lambda)=S_\Lambda \cap \mathcal{P}(\Delta_\Lambda)$.  
\end{lem}
\begin{proof}
If $F\in S_\Lambda \cap \mathcal{P}(\Delta_\Lambda)$, then $F$ is of the form \[F=\{(\lambda_i,\lambda_i)~|~i=1,2,\ldots,n~\text{and}~\Lambda^{\min}(\lambda_i,\lambda_j)=\emptyset~\text{for}~i\neq j\}.\] It is then clear that $F^2=F$. Now assume that $F\in \mathcal{E}(S_\Lambda)\setminus \{\emptyset\}$. Let $(\lambda,\mu)\in F$. Since $F=F^2$, there exist $(\gamma,\delta),(\xi,\eta)\in F$ such that $(\lambda,\mu)\in (\gamma,\delta)\cdot (\xi,\eta)$. Note that $(\eta,\xi)\in F^{-1}=F$. If $(\gamma,\delta)\neq (\eta,\xi)$, then we would have $(\gamma,\delta)\perp (\eta,\xi)$, whence $\Lambda^{\min}(\delta,\xi)=\emptyset$. But this implies $(\gamma,\delta)\cdot (\xi,\eta)=\emptyset$, which is not possible. Therefore, we must have $(\gamma,\delta)=(\eta,\xi)$, which implies $\gamma=\eta$ and $\delta=\xi$. Then $(\lambda,\mu)\in (\gamma,\delta)\cdot (\delta,\gamma)=\{(\gamma,\gamma)\}$, and we have $\lambda=\mu=\gamma$. Hence, $F\subseteq \Delta_\Lambda$ and the lemma follows.   
\end{proof}

Let $\Gamma$ be a directed graph. Recall the \emph{graph inverse semigroup} \cite{Ash-Hall, Jones, Mesyan} $S(\Gamma)$ of $\Gamma$ is the semigroup with $0$ generated by $\{v,e,f^{-1}~|~v\in \Gamma^0,e,f\in\Gamma^1\}$ (where $e^{-1}$ is a formal symbol for each $e\in \Gamma^1$) subject to the relations:

$(i)$ $uv=\delta{u,v}u$ for all $u,v\in \Gamma^0$;

$(ii)$ $s_\Gamma(e)e=e=er_\Gamma(e)$ and $r_\Gamma(e)e^{-1}=e^{-1}=e^{-1}s_\Gamma(e)$ for all $e\in \Gamma^1$;

$(iii)$ $e^{-1}f=\delta_{e,f}r_\Gamma(e)$ for all $e,f\in \Gamma^1$. 

For any $\lambda=e_1 e_2\cdots e_n\in \Gamma^*$, let $\lambda^{-1}:=e_n^{-1} e_{n-1}^{-1} \cdots e_1^{-1}$. Also for any $v\in \Gamma^0$, let $v^{-1}:=v$. Then any nonzero element of $S(\Gamma)$ can be uniquely written as $\alpha\beta^{-1}$ for paths $\alpha,\beta$ with $r_\Gamma(\alpha)=r_\Gamma(\beta)$. It is interesting to note the following connection.
\begin{lem}\label{lem the relation between gis and S_Lambda}
For any directed graph $\gamma$, the map $\psi_\Gamma:S(\Gamma)\longrightarrow S_{\Gamma^*}$ defined by $\psi_\Gamma(\alpha\beta^{-1}):=\{(\alpha,\beta)\}$ and $\psi_\Gamma(0):=\emptyset$ is an injective semigroup homomorphism. Consequently, $S(\Gamma)$ can be considered as an inverse subsemigroup of $S_{\Gamma^*}$.
\end{lem}
One may realize that $\image(\psi_\Gamma)$ is in fact the Paterson's model \cite{Paterson} for the graph inverse semigroup. The map $\psi_\Gamma$ is surjective if and only if $\Gamma$ is either the trivial graph with a single vertex and no edges or the rose with one petal. In the first case, $S_{\Gamma^*}(\cong S(\Gamma))$ is the $2$-element Boolean algebra $\textbf{2}$ and in the second case, it is the bicyclic monoid. 

We conclude this section by describing the Farthing--Muhly--Yeend inverse semigroup of a simple looking $2$-graph. Let $\Delta:=\mathbb{N}^2$, viewed as a category with $1$ object. Then $\Delta$ together with $d:=\text{id}_{\mathbb{N}^2}$ is a finitely aligned $2$-graph. Note that for any two $\N,\M\in \Delta$, $\Delta^{\min}(\N,\M)=\{\left((\N\vee \M)-\N,(\N\vee \M)-\M\right)\}$. It follows that any element of $S_\Delta$ is either a singleton of the form $\{(\N,\M)\}$ or $\emptyset$. So we can identify $S_\Delta$ with $\mathbb{N}^2\times \mathbb{N}^2$. Moreover, \[\{(\N,\M)\}\{(\bP,\Q)\}=\{\left(\N+(\M\vee \bP)-\M,\Q+(\M\vee \bP)-\bP\right)\},\] which is analogous to the multiplication of the bicyclic monoid. Hence, $S_\Delta$ can be viewed as the \emph{two-dimensional version} of the bicyclic monoid. 

\section{Green's relations on \texorpdfstring{$S_\Lambda$}{SΛ}}\label{sec Green's relations}
In this section, we aim to characterize the Green's relations on the inverse semigroup $S_\Lambda$. We start with the relations $\mathscr{L}$, $\mathscr{R}$ and $\mathscr{H}$.

\begin{prop}\label{pro the relations L,R,H}
Let $F,G\in S_\Lambda$. Then the following hold.

$(i)$ $F~\mathscr{L}~G$ if and only if $\mathsf{R}(F)=\mathsf{R}(G)$.

$(ii)$ $F~\mathscr{R}~G$ if and only if $\mathsf{D}(F)=\mathsf{D}(G)$.

$(iii)$ $F~\mathscr{H}~G$ if and only if $\mathsf{D}(F)=\mathsf{D}(G)$ and $\mathsf{R}(F)=\mathsf{R}(G)$.
\end{prop}
\begin{proof}
$(i)$ We observe that $F^{-1}F=\Delta_{\mathsf{R}(F)}$. Let $\mu\in \mathsf{R}(F)$. Then there exists $\lambda\in \Lambda$ such that $(\lambda,\mu)\in F$. It is clear that $(\mu,\lambda)\cdot (\lambda,\mu)=\{(\mu,\mu)\}$. Thus $(\mu,\mu)\in F^{-1}F$, which proves that $\Delta_{\mathsf{R}(F)}\subseteq F^{-1}F$. Now let $(\lambda,\mu)\in F$ and $(\eta,\delta)\in F^{-1}$ be such that $(\eta,\delta)\cdot (\lambda,\mu)\neq \emptyset$. Then $\Lambda^{\min}(\delta,\lambda)\neq \emptyset$. Since $(\delta,\eta)\in F$, we must have $(\lambda,\mu)=(\delta,\eta)$. Then $(\eta,\delta)\cdot (\lambda,\mu)=(\mu,\lambda)\cdot (\lambda,\mu)=\{(\mu,\mu)\}$. Since $\mu\in \mathsf{R}(F)$, it follows that $F^{-1}F=\displaystyle{\bigcup_{\substack{(\eta,\delta)\in F^{-1}\\(\lambda,\mu)\in F}}}(\eta,\delta)\cdot (\lambda,\mu)\subseteq \Delta_{\mathsf{R}(F)}$. Hence the desired equality follows. It is now evident that $F^{-1}F=G^{-1}G$ if and only if $\mathsf{R}(F)=\mathsf{R}(G)$. The result now follows from Lemma \ref{lem L and R on inverse semigroups}$(iii)$. 

$(ii)$ Similar to $(i)$.

$(iii)$ Follows from $(i)$ and $(ii)$.
\end{proof}
\begin{rmks}\label{rem recovering L,R of graph inverse semigroups}
$(i)$ Let $\Gamma$ be a directed graph and $\alpha\beta^{-1},\gamma\delta^{-1}\in S(\Gamma)$. Then $F:=\{(\alpha,\beta)\}$, $G:=\{(\gamma,\delta)\}$ are in $S_{\Gamma^*}$. Note that $\mathsf{D}(F)=\{\alpha\}$, $\mathsf{R}(F)=\{\beta\}$, $\mathsf{D}(G)=\{\gamma\}$ and $\mathsf{R}(G)=\{\delta\}$. Using Proposition \ref{pro the relations L,R,H}, we have $F~\mathscr{L}~G$ if and only if $\beta=\delta$, and $F~\mathscr{R}~G$ if and only if $\alpha=\gamma$. Thus one recovers (in view of Lemma \ref{lem the relation between gis and S_Lambda}) the characterizations of the Green's relations $\mathscr{L}$, $\mathscr{R}$ and $\mathscr{H}$ on the graph inverse semigroup $S(\Gamma)$ obtained in \cite[Corollary 2]{Mesyan}. 

$(ii)$ Unlike the graph inverse semigroup case \cite[Corollary 2(4)]{Mesyan}, the Green's relation $\mathscr{H}$ on $S_\Lambda$ is, in general, not the diagonal relation on $S_\Lambda$. For a simple example, consider the finitely aligned $2$-graph $\Lambda$ with the following $1$-skeleton:
\[
\begin{tikzpicture}[>=stealth, line width=1.2pt]

    \fill (0,0) circle (2pt);
    \fill (2,0) circle (2pt);
    \fill (0,-2) circle (2pt);
    \fill (2,-2) circle (2pt);

    \draw[->, blue] (2,0) -- (0,0)
        node[midway, above] {\textcolor{black}{$\lambda$}};

    \draw[->, blue] (2,-2) -- (0,-2);

    \draw[->, dashed, red, line width=1.2pt]
        (0,0) -- (0,-2);
        
    \draw[->, dashed, red, line width=1.2pt]
        (2,0) -- (2,-2)
        node[midway, right] {\textcolor{black}{$\mu$}};

    \node[left] at (0,0) {$w$};
    \node[right] at (2,0) {$z$};
    \node[left] at (0,-2) {$u$};
    \node[right] at (2,-2) {$v$};

\end{tikzpicture}
\]
where solid blue edges have degree $(1,0)$ and the dashed red edges have degree $(0,1)$. Let $F=\{(\lambda,\lambda),(\mu,\mu)\}$ and $G=\{(\lambda,\mu),(\mu,\lambda)\}$. Then $F,G\in S_\Lambda$ and $\mathsf{D}(F)=\mathsf{D}(G)=\{\lambda,\mu\}=\mathsf{R}(F)=\mathsf{R}(G)$. Hence $F~\mathscr{H}~G$ by Proposition \ref{pro the relations L,R,H}$(iii)$; however $F\neq G$.\qed 
\end{rmks}

We now proceed to characterize Green's $\mathscr{D}$ relation on $S_\Lambda$. For any $F\in S_\Lambda$, we write \[\mathfrak{s}(F):=s(\mathsf{D}(F))=s(\mathsf{R}(F)).\] The following lemma provides a necessary condition for two elements of $S_\Lambda$ to be $\mathscr{D}$-related.
\begin{lem}\label{lem Necessary condition for D}
Let $\Lambda$ be a finitely aligned $k$-graph and $F,G\in S_\Lambda$. If $F~\mathscr{D}~G$, then $\mathfrak{s}(F)=\mathfrak{s}(G)$.   
\end{lem}
\begin{proof}
Let $F~\mathscr{D}~G$. Then there is an $H\in S_\Lambda$ such that $F~\mathscr{L}~H$ and $H~\mathscr{R}~G$. Using Proposition \ref{pro the relations L,R,H}, we then have $\mathsf{R}(F)=\mathsf{R}(H)$ and $\mathsf{D}(H)=\mathsf{D}(G)$. These imply \[\mathfrak{s}(F)=s(\mathsf{R}(F))=s(\mathsf{R}(H))=s(\mathsf{D}(H))=s(\mathsf{D}(G))=\mathfrak{s}(G).\qedhere \] 
\end{proof}
We now provide an example to show that the necessary condition stated in Lemma \ref{lem Necessary condition for D} is not always sufficient.
\begin{example}\label{ex Necessary but no sufficient}
Let $\Gamma$ be the directed graph shown below and $\Lambda=\Gamma^*$ the path category of $\Gamma$ viewed as a $1$-graph. 
\[
\begin{tikzpicture}[
    >=stealth,
    every node/.style={font=\large},
    vertex/.style={circle, fill=black, inner sep=1.5pt}
]

\node[vertex] (w1) at (1,1.5) {};
\node[vertex] (w2) at (1,0.6) {};
\node[vertex] (u)  at (3,1.1) {};

\node[vertex] (v)  at (5,1.1) {};
\node[vertex] (z1) at (7,0.6) {};
\node[vertex] (z2) at (7,1.5) {};

\node[left=3pt of w1] {\normalsize{$w_1$}};
\node[left=3pt of w2] {\normalsize{$w_2$}};
\node[below=5pt of u] {\normalsize{$u$}};

\node[below=5pt of v] {\normalsize{$v$}};
\node[right=3pt of z1] {\normalsize{$z_1$}};
\node[right=3pt of z2] {\normalsize{$z_2$}};

\draw[->, thick] (w1) -- node[above] {\normalsize{$e_1$}} (u);
\draw[->, thick] (w2) -- node[below] {\normalsize{$e_2$}} (u);

\draw[->, thick]
    (v) to[out=155,in=25]
    node[above] {\normalsize{$e_3$}} (u);

\draw[->, thick]
    (u) to[out=-25,in=205]
    node[below] {\normalsize{$f_3$}} (v);

\draw[->, thick] (z1) -- node[below] {\normalsize{$f_1$}} (v);
\draw[->, thick] (z2) -- node[above] {\normalsize{$f_2$}} (v);

\end{tikzpicture}
\]
Consider the idempotents $F=\{(e_1,e_1),(f_2,f_2),(f_3,f_3)\}$ and $G=\{(f_1,f_1),(e_2,e_2),(e_3,e_3)\}$ in $S_\Lambda$. Then $\mathfrak{s}(F)=\mathfrak{s}(G)=\{u,v\}$. Now if $F~\mathscr{D}~G$, then there is an $H\in S_\Lambda$ such that $F~\mathscr{L}~H$ and $H~\mathscr{R}~G$, which, in view of Proposition \ref{pro the relations L,R,H}, implies that \[\mathsf{R}(H)=\mathsf{R}(F)=\{e_1,f_2,f_3\},~~\mathsf{D}(H)=\mathsf{D}(G)=\{f_1,e_2,e_3\}.\] Since none of the pairs $(e_2,f_2)$, $(e_2,f_3)$, $(e_3,f_2)$ and $(e_3,f_3)$ is in $\Lambda *_s\Lambda$, we must have $(e_2,e_1),(e_3,e_1)\in H$. But this is not possible in view of the definition of $S_\Lambda$. So $F$ and $G$ are not $\mathscr{D}$-related. 
\end{example}
In order to obtain a necessary and sufficient condition characterizing the $\mathscr{D}$-relation, we define the following. Let $F\in S_\Lambda$ and $m_F:\mathfrak{s}(F)\longrightarrow \mathbb{Z}^+$ be the map such that \[m_F(u):=|\{(\lambda,\mu)\in F~|~s(\lambda)=s(\mu)=u\}|\] for all $u\in \mathfrak{s}(F)$.

Observe that for any $u\in \mathfrak{s}(F)$, the map $\theta:\mathsf{D}(F)u\longrightarrow \mathsf{R}(F)u$ sending each $\lambda\in \mathsf{D}(F)u$ to the unique element $\mu$ such that $(\lambda,\mu)\in F$ is a bijection and $|\mathsf{D}(F)u|=|\mathsf{R}(F)u|=m_F(u)$. 
\begin{prop}\label{pro necessary and sufficient condition for D}
Let $\Lambda$ be a finitely aligned $k$-graph and $F,G\in S_\Lambda$. Then $F~\mathscr{D}~G$ if and only if $m_F=m_G$ as functions.     
\end{prop}
\begin{proof}
Suppose $F~\mathscr{D}~G$. That the domains of $m_F$ and $m_G$ are equal is exactly Lemma \ref{lem Necessary condition for D}. Let $u\in \mathfrak{s}(F)$. Choose any $H\in S_\Lambda$ such that $F~\mathscr{L}~H$ and $H~\mathscr{R}~G$. Then using Proposition \ref{pro the relations L,R,H} and the observation above Proposition \ref{pro necessary and sufficient condition for D}, we have \[m_F(u)=|\mathsf{R}(F)u|=|\mathsf{R}(H)u|=|\mathsf{D}(H)u|=|\mathsf{D}(G)u|=m_G(u).\] 

Now assume that the condition holds. For each $u\in \mathfrak{s}(F)=\mathfrak{s}(G)$, since $|\mathsf{D}(G)u|=m_G(u)=m_F(u)=|\mathsf{R}(F)u|$, we can fix a bijection $f_u:\mathsf{D}(G)u\longrightarrow \mathsf{R}(F)u$. Now define \[H:=\{(\lambda,f_u(\lambda))~|~u\in \mathfrak{s}(F),\lambda\in \mathsf{D}(G)u\}.\] We show that $H\in S_\Lambda$. Clearly $H$ is finite. Suppose $(\lambda,f_u(\lambda))$ and $(\xi,f_v(\xi))$ are distinct pairs in $H$. If $u\neq v$, then $\lambda\neq \xi$ and $f_u(\lambda)\neq f_v(\xi)$. On the other hand, if $u=v$, then also $\lambda\neq \xi$ and $f_u(\lambda)\neq f_u(\xi)$ since $f_u$ is a bijection. Hence in any case, since $\lambda,\xi\in \mathsf{D}(G)$ and $f_u(\lambda),f_v(\xi)\in \mathsf{R}(F)$, it follows that $\Lambda^{\min}(\lambda,\xi)=\Lambda^{\min}(f_u(\lambda),f_v(\xi))=\emptyset$. Therefore, $(\lambda,f_u(\lambda))\perp (\xi,f_v(\xi))$. This shows that $H\in S_\Lambda$. Now \[\mathsf{R}(H)=\{f_u(\lambda)~|~u\in \mathfrak{s}(F),\lambda\in \mathsf{D}(G)u\}=\displaystyle{\bigsqcup_{u\in \mathfrak{s}(F)}}\mathsf{R}(F)u=\mathsf{R}(F),\] and \[\mathsf{D}(H)=\displaystyle{\bigsqcup_{u\in \mathfrak{s}(G)}}\mathsf{D}(G)u=\mathsf{D}(G).\] By Proposition \ref{pro the relations L,R,H}, $F~\mathscr{L}~H$ and $H~\mathscr{R}~G$ and consequently, $F~\mathscr{D}~G$. Hence the sufficiency follows. 
\end{proof}
\begin{rmk}\label{rem recovering D of graph inverse semigroups}
If $\Lambda$ is the path category of a directed graph $\Gamma$ and $F=\{(\alpha,\beta)\}$, $G=\{(\gamma,\delta)\}$ in $\image(\psi)\subseteq S_\Lambda$ (see Lemma \ref{lem the relation between gis and S_Lambda}), then $\mathfrak{s}(F)=\mathfrak{s}(G)$ amounts to say that $r_\Gamma(\alpha)=r_\Gamma(\gamma)$. The reader may find that this is exactly the necessary and sufficient condition for $\alpha\beta^{-1}$ and $\gamma\delta^{-1}$ to be $\mathscr{D}$-related in $S(\Gamma)$ (see \cite[Corollary 2(5)]{Mesyan}). Hence, Proposition \ref{pro necessary and sufficient condition for D} perfectly extends the characterization of Green's $\mathscr{D}$ relation on graph inverse semigroups to the inverse semigroups of finitely aligned higher-rank graphs.\qed
\end{rmk}

We are left with the characterization of Green's $\mathscr{J}$ relation on $S_\Lambda$. In general, it is not easy to tell whether two arbitrary elements of $S_\Lambda$ are $\mathscr{J}$-related. Thus we focus on specific elements of $S_\Lambda$ defined below.
\begin{dfn}\label{def divergent set}
Let $\Lambda$ be a finitely aligned $k$-graph and $S_\Lambda$ the Farthing--Muhly--Yeend inverse semigroup of $\Lambda$. An element $F\neq \emptyset$ of $S_\Lambda$ is called \emph{divergent} if the map $s|_{\mathsf{D}(F)}:\mathsf{D}(F)\longrightarrow \mathfrak{s}(F)$ is bijective.     
\end{dfn}

Before describing the $\mathscr{J}$ relation on divergent elements, we notice a simple but useful way of expressing these elements. First observe that if $F\in S_\Lambda$ is divergent, then \[F_\mathsf{D}:=\{(\lambda,s(\lambda))~|~\lambda\in \mathsf{D}(F)\},\] and \[F_\mathsf{R}:=\{(\mu,s(\mu))~|~\mu\in \mathsf{R}(F)\}\] are elements of $S_\Lambda$. Now, a moment's thought reflects that $F$ can be written as 
\begin{equation}\label{eq div1}
F=F_{\mathsf{D}}\Delta_{\mathfrak{s}(F)}F_{\mathsf{R}}^{-1}.    
\end{equation} 
Since $\mathsf{R}(F_\mathsf{D})=\mathfrak{s}(F)=\mathsf{R}(F_\mathsf{R})$ and $\Delta_{\mathfrak{s}(F)}\in \mathcal{E}(S_\Lambda)$, multiplying (\ref{eq div1}) by $F_\mathsf{D}^{-1}$ from the left and by $F_\mathsf{R}$ from the right, we obtain 
\begin{equation}\label{eq div2}
\Delta_{\mathfrak{s}(F)}=F_{\mathsf{D}}^{-1}F F_\mathsf{R}.
\end{equation}
From (\ref{eq div1}) and (\ref{eq div2}), it follows that $F~\mathscr{J}~\Delta_{\mathfrak{s}(F)}$. 
 
We now present the following theorem which describes when two divergent elements are $\mathscr{J}$-related.
\begin{prop}\label{pro partial orders on equivalence classes J}
Let $\Lambda$ be any finitely aligned $k$-graph and $F,G$ divergent elements in $S_\Lambda$. Then $[F]_{\mathscr{J}}\le_{\mathscr{J}} [G]_{\mathscr{J}}$ if and only if there exists a pairwise orthogonal subset $\mathcal{C}$ of $\Lambda$ such that for any $u\in \mathfrak{s}(F)$, $\mathfrak{s}(G)\mathcal{C}u\neq \emptyset$.
\end{prop}
\begin{proof}
Suppose $[F]_{\mathscr{J}}\le_{\mathscr{J}} [G]_{\mathscr{J}}$. By the above observation, we have $\Delta_{\mathfrak{s}(F)}\in S_\Lambda \Delta_{\mathfrak{s}(G)} S_\Lambda$. This implies that $\Delta_{\mathfrak{s}(F)}=H \Delta_{\mathfrak{s}(G)}K$ for some $H,K\in S_\Lambda$. Choose any $u\in \mathfrak{s}(F)$. Then $(u,u)\in \Delta_{\mathfrak{s}(F)}=H \Delta_{\mathfrak{s}(G)}K$. Hence there are $(\lambda_u,\mu_u)\in H$, $v\in \mathfrak{s}(G)$ and $(\xi_u,\eta_u)\in K$ such that \[(u,u)\in (\lambda_u,\mu_u)\cdot (v,v)\cdot (\xi_u,\eta_u).\] We must have $r(\mu_u)=v=r(\xi_u)$. Then $(\lambda_u,\mu_u)\cdot (v,v)=\{(\lambda_u,\mu_u)\}$ and consequently, $u=\lambda_u \alpha=\eta_u\beta$ for some $(\alpha,\beta)\in \Lambda^{\min}(\mu_u,\xi_u)$. This implies that $\lambda_u=\eta_u=u$ and so $s(\mu_u)=s(\lambda_u)=u$. Now let $\mathcal{A}:=\{\mu_u~|~u\in \mathfrak{s}(F)\}$. By the construction, it is clear that $\mathfrak{s}(G)\mathcal{A}u\neq \emptyset$ for all $u\in \mathfrak{s}(F)$. Also since $\mathcal{A}\subseteq \mathsf{R}(H)$, $\mathcal{A}$ is pairwise orthogonal. 

Conversely assume that there exists $\mathcal{C}\subseteq \Lambda$ such that $\mathcal{C}$ is pairwise orthogonal and $\mathfrak{s}(G)\mathcal{C}u\neq \emptyset$ for all $u\in \mathfrak{s}(F)$. Again by the observation made prior to this proposition, it is enough to show that $\Delta_{\mathfrak{s}(F)}\in S_\Lambda \Delta_{\mathfrak{s}(G)} S_\Lambda$. For each $u\in \mathfrak{s}(F)$, we fix a $\lambda_u\in \mathfrak{s}(G)\mathcal{C}u$ and define \[F:=\{(u,\lambda_u)~|~u\in \mathfrak{s}(F)\}.\] Since $\mathcal{C}$ is pairwise orthogonal, $F\in S_\Lambda$. We claim that $\Delta_{\mathfrak{s}(F)}=F \Delta_{\mathfrak{s}(G)} F^{-1}$. For any $u\in \mathfrak{s}(F)$, $(r(\lambda_u),r(\lambda_u))\in \Delta_{\mathfrak{s}(G)}$ and \[(u,\lambda_u)\cdot (r(\lambda_u),r(\lambda_u))\cdot (\lambda_u,u)=\{(u,u)\}.\] This implies $(u,u)\in F\Delta_{\mathfrak{s}(G)}F^{-1}$ and hence $\Delta_{\mathfrak{s}(F)}\subseteq F\Delta_{\mathfrak{s}(G)}F^{-1}$. Again since 
\[(u,\lambda_u)\cdot (v,v)\cdot (\lambda_w,w):=
	\left\{
	\begin{array}{ll}
		 (u,u)  & \mbox{if } u=w~\text{and}~r(\lambda_u)=v, \\
	   \emptyset  & \mbox{otherwise},
	\end{array}
	\right.
\] it follows that $F\Delta_{\mathfrak{s}(G)}F^{-1}\subseteq \Delta_{\mathfrak{s}(F)}$. Therefore, the claim is established and we have $\Delta_{\mathfrak{s}(F)}\in S_\Lambda \Delta_{\mathfrak{s}(G)} S_\Lambda$. 
\end{proof}

Suppose $X,Y\subseteq \Lambda^0$. We say that $X$ and $Y$ are \emph{orthogonally connected} if there are pairwise orthogonal subsets $\mathcal{A},\mathcal{B}\subseteq \Lambda$ such that for any $u\in X$ and $v\in Y$, both $Y\mathcal{A}u$ and $X\mathcal{B}v$ are nonempty. In view of this definition, Proposition \ref{pro partial orders on equivalence classes J} immediately yields the following corollary.
\begin{cor}\label{cor J relation on divergent sets}
Let $\Lambda$ be any finitely aligned $k$-graph and $F,G$ divergent elements of $S_\Lambda$. Then $F~\mathscr{J}~G$ if and only if $\mathfrak{s}(F)$ and $\mathfrak{s}(G)$ are orthogonally connected. 
\end{cor}

We illustrate the use of Proposition \ref{pro partial orders on equivalence classes J} and Corollary \ref{cor J relation on divergent sets} via some examples.
\begin{examples}\label{ex illustrating J characterization}
$(i)$ Consider any directed graph $\Gamma$ and let $\Lambda=\Gamma^*$. Suppose $\alpha\beta^{-1},\gamma\delta^{-1}\in S(\Gamma)$. Then $F=\{(\alpha,\beta)\}$ and $G=\{(\gamma,\delta)\}$ are elements in $S_\Lambda$. These are trivially divergent elements. Hence we can apply Proposition \ref{pro partial orders on equivalence classes J} to conclude that $[F]_{\mathscr{J}}\le_{\mathscr{J}} [G]_{\mathscr{J}}$ in $S_\Lambda$ if and only if there is a pairwise orthogonal subset $\mathcal{C}$ of $\Lambda$ such that $r_\Gamma(\gamma)\mathcal{C}r_\Gamma(\alpha)\neq \emptyset$. This is equivalent to the existence of a path $\tau$ in $\Gamma$ such that $s_\Gamma(\tau)=r_\Gamma(\gamma)$ and $r_\Gamma(\tau)=r_\Gamma(\alpha)$. Similarly Corollary \ref{cor J relation on divergent sets} implies that $F~\mathscr{J}~G$ if and only if there is a path from $r_\Gamma(\gamma)$ to $r_\Gamma(\alpha)$ and vice-versa, i.e., $r_\Gamma(\alpha)$ and $r_\Gamma(\gamma)$ lie in the same strongly connected component of $\Gamma$. This is exactly the characterization of Green's $\mathscr{J}$ relation on $S(\Gamma)$ (see \cite[Corollary 2(3)]{Mesyan}). 

$(ii)$ Consider the $1$-skeleton of the $2$-graph $\Omega_{2,(2,1)}$ shown below: 
\[
\begin{tikzpicture}[
    >=stealth,
    thick,
    every node/.style={font=\large}
]

\coordinate (A) at (0,2);
\coordinate (B) at (2,2);
\coordinate (C) at (4,2);

\coordinate (D) at (0,0);
\coordinate (E) at (2,0);
\coordinate (F) at (4,0);

\fill (A) circle (2pt);
\fill (B) circle (2pt);
\fill (C) circle (2pt);

\fill (D) circle (2pt);
\fill (E) circle (2pt);
\fill (F) circle (2pt);

\draw[blue, thick, ->] (C) -- (B);
\draw[blue, thick, ->] (B) -- (A);

\draw[blue, thick, ->] (F) -- (E);
\draw[blue, thick, ->] (E) -- (D);

\draw[red, thick, dashed, ->] (A) -- (D);
\draw[red, thick, dashed, ->] (B) -- (E);
\draw[red, thick, dashed, ->] (C) -- (F);

\node[above left=1pt] at (A) {\normalsize{$(0,1)$}};
\node[above=1pt]      at (B) {\normalsize{$(1,1)$}};
\node[above right=1pt] at (C) {\normalsize{$(2,1)$}};

\node[below left=1pt]  at (D) {\normalsize{$(0,0)$}};
\node[below=1pt]       at (E) {\normalsize{$(1,0)$}};
\node[below right=1pt] at (F) {\normalsize{$(2,0)$}};

\node[left=1pt]  at ($(A)!0.5!(D)$) {\normalsize{$\alpha$}};
\node[right=1pt] at ($(B)!0.5!(E)$) {\normalsize{$\beta$}};
\node[right=1pt] at ($(C)!0.5!(F)$) {\normalsize{$\gamma$}};

\node[below=0.5pt] at ($(A)!0.5!(B)$) {\normalsize{$\lambda$}};
\node[below=0.5pt] at ($(B)!0.5!(C)$) {\normalsize{$\mu$}};

\end{tikzpicture}
\]
Let $F=\{(\mu,\gamma),(\lambda,\beta)\}$ and $G=\{(\alpha,\alpha)\}$. Then $\mathfrak{s}(F)=\{(2,1),(1,1)\}$ and $\mathfrak{s}(G)=\{(0,1)\}$. Clearly $F,G$ are divergent elements in $S_{\Omega_{2,(2,1)}}$. Note that $\lambda$ is the only path from $(1,1)$ to $(0,1)$, and $\lambda\mu$ is the only path from $(2,1)$ to $(0,1)$. Therefore, if there is a subset $\mathcal{C}$ of $\Omega_{2,(2,1)}$ such that $\mathfrak{s}(G)\mathcal{C}u\neq \emptyset$ for all $u\in \mathfrak{s}(F)$, then $\lambda,\lambda\mu$ must be in $\mathcal{C}$. But then $\mathcal{C}$ cannot be pairwise orthogonal. So by Proposition \ref{pro partial orders on equivalence classes J}, $[F]_{\mathscr{J}}\nleq_{\mathscr{J}} [G]_{\mathscr{J}}$.

In this same example, if we replace $G$ by $G'=\{(\alpha,\alpha),(\lambda,\beta)\}$ and take $\mathcal{C}=\{\lambda,\mu\}$, then $\mathcal{C}$ is pairwise orthogonal, $\mathfrak{s}(G')\mathcal{C}(2,1)=\{\mu\}$ and $\mathfrak{s}(G')\mathcal{C}(1,1)=\{\lambda\}$. Thus the condition of Proposition \ref{pro partial orders on equivalence classes J} is satisfied and we have $[F]_{\mathscr{J}}\le_{\mathscr{J}} [G']_{\mathscr{J}}$. 
\end{examples}

\section{Natural partial order, covers, joins and abstract Cuntz--Krieger relations}\label{sec partial order, covers, joins}
In this section, we first explore some idempotent-related concepts in the setting of the Farthing--Muhly--Yeend inverse semigroups of a higher-rank graph. 

Recall that an inverse semigroup $S$ is called \emph{fundamental} if the maximal idempotent-separating congruence on $S$ equals the diagonal relation on $S$. This is equivalent to saying that $Z(\mathcal{E}(S))=\mathcal{E}(S)$, where $Z(\mathcal{E}(S))$ is the centralizer of $\mathcal{E}(S)$ (see \cite[Proposition 2.6]{Ganna}). Mesyan and Mitchell in \cite{Mesyan} proved that the graph inverse semigroup $S(\Gamma)$ of any directed graph $\Gamma$ is fundamental. The proof uses the fact that $S(\Gamma)$ is a combinatorial inverse semigroup \cite{Jones}, i.e., Green's $\mathscr{H}$ relation on $S(\Gamma)$ coincides with the diagonal relation. However, the inverse semigroup of a finitely aligned $k$-graph is in general not combinatorial as we saw in Remark \ref{rem recovering L,R of graph inverse semigroups}$(ii)$. Despite this issue, we have the following result.

\begin{prop}\label{pro S_Lambda is fundamental}
Let $\Lambda$ be any finitely aligned $k$-graph. Then the inverse semigroup $S_\Lambda$ is fundamental. 
\end{prop}
\begin{proof}
We show that $Z(\mathcal{E}(S_\Lambda))=\mathcal{E}(S_\Lambda)$. Since $S_\Lambda$ is an inverse semigroup, we already have $\mathcal{E}(S_\Lambda)\subseteq Z(\mathcal{E}(S_\Lambda))$. Now let $F\in Z(\mathcal{E}(S_\Lambda))$ and $(\lambda,\mu)\in F$. Consider the element $E=\{(\lambda,\lambda)\}\in \mathcal{E}(S_\Lambda)$. Then $EF=FE$. Since for any $\eta\in \mathsf{D}(F)\setminus \{\lambda\}$, $\Lambda^{\min}(\lambda,\eta)=\emptyset$ and $\mu$ is the only path such that $(\lambda,\mu)\in F$, we have \[EF=\displaystyle{\bigcup_{(\eta,\delta)\in F}}(\lambda,\lambda)\cdot (\eta,\delta)=(\lambda,\lambda)\cdot (\lambda,\mu)=\{(\lambda,\mu)\}.\] Hence $(\lambda,\mu)\in FE$. Then there exist $(\eta,\delta)\in F$ and $(\alpha,\beta)\in \Lambda^{\min}(\delta,\lambda)$ such that $(\lambda,\mu)=(\eta\alpha,\lambda\beta)$. This gives 
\begin{equation}\label{eq presup1}
\mu=\lambda\beta.    
\end{equation}
Again letting $E'=\{(\mu,\mu)\}\in \mathcal{E}(S_\Lambda)$ and proceeding as above, we obtain 
\begin{equation}\label{eq presup2}
\lambda=\mu\gamma    
\end{equation}
for some $\gamma\in \Lambda$. Combining (\ref{eq presup1}) and (\ref{eq presup2}), we have $\lambda=\lambda\beta\gamma$, which implies $\beta=\gamma=s(\lambda)$. Therefore, $\mu=\lambda s(\lambda)=\lambda$. This shows that $F\subseteq \Delta_{\Lambda}$ and hence $F\in \mathcal{E}(S_\Lambda)$ by Lemma \ref{lem idempotents of the semigroup}. Since $F$ was chosen arbitrarily from $Z(\mathcal{E}(S_\Lambda))$, this implies that $Z(\mathcal{E}(S_\Lambda))\subseteq \mathcal{E}(S_\Lambda)$ and thus completes the proof. 
\end{proof}

Compatibility of elements in an inverse semigroup (see Subsection \ref{ssec invsem}) is another important concept related to idempotents as this gives a necessary condition for the existence of join and is also used to define Boolean inverse semigroups. We now describe when two elements of the inverse semigroup $S_\Lambda$ are compatible. 
\begin{prop}\label{pro describing compatibility}
Let $\Lambda$ be a finitely aligned $k$-graph and $F,G\in S_\Lambda$. Then

$(i)$ $F\perp G$ if and only if $F\cap G=\emptyset$ and $F\cup G\in S_\Lambda$. 

$(ii)$ $F\sim G$ if and only if for any $\lambda\in \mathsf{D}(F)$ and $\xi\in \mathsf{D}(G)$, $\Lambda^{\min}(\lambda,\xi)=\Lambda^{\min}(b_F(\lambda),b_G(\xi))$.   
\end{prop}
\begin{proof}
$(i)$ First assume that $F\cap G=\emptyset$ and $F\cup G\in S_\Lambda$. Since $F,G\subseteq F\cup G$, these imply that $\Lambda^{\min}(\mathsf{D}(F),\mathsf{D}(G))=\emptyset=\Lambda^{\min}(\mathsf{R}(F),\mathsf{R}(G))$. It follows that $F^{-1}G=FG^{-1}=\emptyset$ and so $F\perp G$. Conversely, if $F\perp G$, then we must have $F\cap G=\emptyset$, otherwise if $(\lambda,\mu)\in F\cap G$, then $\{(\mu,\mu)\}=(\mu,\lambda)\cdot (\lambda,\mu)\subseteq F^{-1}G$ and we would have $F^{-1}G\neq \emptyset$. Now let $(\lambda,\mu)\in F$ and $(\xi,\eta)\in G$. Since $F^{-1}G=FG^{-1}=\emptyset$, $\Lambda^{\min}(\lambda,\xi)=\Lambda^{\min}(\mu,\eta)=\emptyset$. Since $F,G\in S_\Lambda$, this is sufficient to conclude that $F\cup G\in S_\Lambda$. 

$(ii)$ First we show that the condition is sufficient. Let $(\lambda,\mu)\in F$ and $(\xi,\eta)\in G$. Let $(x,y)\in (\mu,\lambda)\cdot (\xi,\eta)$. Then $x=\mu\alpha$ and $y=\eta\beta$ for some $(\alpha,\beta)\in \Lambda^{\min}(\lambda,\xi)$. By the condition, $(\alpha,\beta)\in \Lambda^{\min}(\mu,\eta)$ and hence $x=y$. By Lemma \ref{lem idempotents of the semigroup}, this implies that $(\mu,\lambda)\cdot (\xi,\eta)\in \mathcal{E}(S_\Lambda)$ for all $(\lambda,\mu)\in F$ and $(\xi,\eta)\in G$. Therefore, $F^{-1}G\in \mathcal{E}(S_\Lambda)$ by (\ref{eq the multiplication}). Similarly, it can be shown that $FG^{-1}\in \mathcal{E}(S_\Lambda)$. Hence $F\sim G$.

Now assume that $F\sim G$. Fix any $\lambda\in \mathsf{D}(F)$ and $\xi\in \mathsf{D}(G)$. For simplicity, we write $\mu:=b_F(\lambda)$ and $\eta:=b_G(\xi)$. Let $(\alpha,\beta)\in \Lambda^{\min}(\lambda,\xi)$. Then $(\mu\alpha,\eta\beta)\in F^{-1}G$. Since $F^{-1}G$ is an idempotent, we have $\mu\alpha=\eta\beta$. This implies $\Lambda^{\min}(\mu,\eta)\neq \emptyset$ and 
\begin{equation}\label{comp1}
d(\alpha)=d(\mu\alpha)-d(\mu)\ge \left(d(\mu)\vee d(\eta)\right)-d(\mu).
\end{equation}
Now choose any $(\gamma,\delta)\in \Lambda^{\min}(\mu,\eta)$. Then $(\lambda\gamma,\xi\delta)\in FG^{-1}$. Since $FG^{-1}$ is also an idempotent, $\lambda\gamma=\xi\delta$. But then 
\begin{equation}\label{comp2}
\left(d(\mu)\vee d(\eta)\right)-d(\mu)=d(\gamma)=d(\lambda\gamma)-d(\lambda)\ge \left(d(\lambda)\vee d(\xi)\right)-d(\lambda)=d(\alpha)
\end{equation}
Using the inequalities (\ref{comp1}) and (\ref{comp2}), we have $d(\alpha)=\left(d(\mu)\vee d(\eta)\right)-d(\mu)$, whence $(\alpha,\beta)\in \Lambda^{\min}(\mu,\eta)$. Therefore, we have shown that whenever $F\sim G$, then $\Lambda^{\min}(\lambda,\xi)\subseteq \Lambda^{\min}(b_F(\lambda),b_G(\xi))$ for any $\lambda\in \mathsf{D}(F)$ and $\xi\in \mathsf{D}(G)$. Now $F\sim G$ implies that $F^{-1}\sim G^{-1}$. Using the previous containment for $b_F(\lambda)\in \mathsf{D}(F^{-1})$ and $\b_G(\xi)\in \mathsf{D}(G^{-1})$, we obtain \[\Lambda^{\min}(b_F(\lambda),b_G(\xi))\subseteq \Lambda^{\min}\left(b_{F^{-1}}(b_F(\lambda)),b_{G^{-1}}(b_G(\xi))\right)=\Lambda^{\min}(\lambda,\xi).\] This completes the proof.
\end{proof}
  
\begin{rmk}\label{rem compatibility consequences}
As a consequence of Proposition \ref{pro describing compatibility}, one can describe compatibility of two elements in a graph inverse semigroup. Let $\Gamma$ be a directed graph and $\alpha\beta^{-1},\gamma\delta^{-1}\in S(\Gamma)$. Then $\{(\alpha,\beta)\},\{(\gamma,\delta)\}\in \image(\psi_\Gamma)$. Note that $\alpha\beta^{-1}\sim \gamma\delta^{-1}$ in $S(\Gamma)$ if and only if $\{(\alpha,\beta)\}\sim \{(\gamma,\delta)\}$ in $S_{\Gamma^*}$. By Proposition \ref{pro describing compatibility}, the latter holds if and only if $(\Gamma^*)^{\min}(\alpha,\gamma)=(\Gamma^*)^{\min}(\beta,\delta)$. If these sets are empty, then $(\alpha,\beta)\perp (\gamma,\delta)$. If these are nonempty, then either $\alpha=\gamma\nu$, $\beta=\delta\nu$ for some $\nu\in \Gamma^*$ or $\gamma=\alpha\tau$, $\delta=\beta\tau$ for some $\tau\in \Gamma^*$. Incorporating this with \cite[Lemma 15(2)]{Mesyan}, it follows that two elements of a graph inverse semigroup are compatible if and only if they are either orthogonal or comparable.\qed     
\end{rmk}

We now delve into the ordered structure of the inverse semigroup $S_\Lambda$ with respect to the natural partial order and prove some fundamental results, which will lay the foundation for the next section. In the following proposition, we provide a geometric description of the natural partial order on $S_\Lambda$. 

\begin{prop}\label{pro natural partial order on S_Lambda}
Let $\Lambda$ be a finitely aligned $k$-graph and $F,G\in S_\Lambda$. Then $F\le G$ if and only if for each $(\lambda,\mu)\in F$, there exist unique $\beta\in \Lambda$ and $(\xi,\eta)\in G$ such that $\lambda=\xi\beta$ and $\mu=\eta\beta$. 
\end{prop}
\begin{proof}
Suppose $F\le G$. Then $F=EG$ for some $E\in \mathcal{E}(S_\Lambda)$. Let $(\lambda,\mu)\in F$. In view of Lemma \ref{lem idempotents of the semigroup} and (\ref{eq the multiplication}), we can say that there are $(\gamma,\gamma)\in E$ and $(\xi,\eta)\in G$ such that $(\lambda,\mu)\in (\gamma,\gamma)\cdot (\xi,\eta)$. It follows that $(\lambda,\mu)=(\gamma\alpha,\eta\beta)$ for some $(\alpha,\beta)\in \Lambda^{\min}(\gamma,\xi)$. Hence, $\lambda=\gamma\alpha=\xi\beta$ and $\mu=\eta\beta$. If there is another pair $(\xi',\eta')\in G$ and a path $\beta'\in \Lambda$ satisfying the same equalities, then we have \[\lambda=\xi\beta=\xi'\beta',~~\mu=\eta\beta=\eta'\beta'.\] These imply that $\lambda(0,d(\xi)\vee d(\xi'))\in \MCE(\xi,\xi')$ and $\mu(0,d(\eta)\vee d(\eta'))\in \MCE(\eta,\eta')$, and consequently, $(\xi,\eta)$ and $(\xi',\eta')$ are not componentwise orthogonal. Since they both are in $G$, we must have $(\xi,\eta)=(\xi',\eta')$. Then we have $\xi\beta=\xi\beta'$, which implies that $\beta=\beta'$. Hence, the uniqueness follows. 

Now suppose the condition holds. Let $E:=\Delta_{\mathsf{D}(F)}$. Since $F\in S_\Lambda$, $\mathsf{D}(F)$ is a finite pairwise orthogonal subset of $\Lambda$. This, together with Lemma \ref{lem idempotents of the semigroup}, implies that $E\in \mathcal{E}(S_\Lambda)$. We claim that $F=EG$. Let $(\lambda,\mu)\in F$. By the condition, there are $(\xi,\eta)\in G$ and $\beta\in \Lambda$ such that $\lambda=\xi\beta$ and $\mu=\eta\beta$. Then $(s(\lambda),\beta)\in \Lambda^{\min}(\lambda,\xi)$ and hence, $(\lambda,\mu)\in (\lambda,\lambda)\cdot (\xi,\eta)\subseteq EG$. Therefore, $F\subseteq EG$. Now for the reverse inclusion, let $(\gamma,\gamma)\in E$ and $(\tau,\delta)\in G$ be such that $(\gamma,\gamma)\cdot (\tau,\delta)\neq \emptyset$. Then $\Lambda^{\min}(\gamma,\tau)\neq \emptyset$. Since $\gamma\in \mathsf{D}(F)$, $(\gamma,\nu)\in F$ for some $\nu\in \Lambda$. Then there are $(\xi,\eta)\in G$ and $\beta\in \Lambda$ such that $\gamma=\xi\beta$ and $\nu=\eta\beta$. From these observations, it follows that $\Lambda^{\min}(\xi,\tau)\neq \emptyset$. But since $(\xi,\eta),(\tau,\delta)\in G$, we should have $\xi=\tau$ and $\eta=\delta$. Then $\Lambda^{\min}(\gamma,\tau)=\{(s(\gamma),\beta)\}$, and we have $(\gamma,\gamma)\cdot (\tau,\delta)=\{(\gamma,\nu)\}\subseteq F$. Then by (\ref{eq the multiplication}), $EG\subseteq F$. Therefore, $F=EG$ and so $F\le G$.   
\end{proof}
The above proposition immediately yields the following corollary, which describes the natural partial order when restricted to the semilattice of idempotents. The reader may compare the first part of this corollary with \cite[Proposition 2.10]{Ortega} and \cite[Lemma 4.4]{Larki}.

\begin{cor}\label{cor partial order on semilattice}
Let $\Lambda$ be a finitely aligned $k$-graph and $F,G\in \mathcal{E}(S_\Lambda)$. Then $F\le G$ if and only if for each $\lambda\in \mathsf{D}(F)$, there exist unique $\xi\in \mathsf{D}(G)$ and $\beta\in \Lambda$ such that $\lambda=\xi\beta$. The semilattice of idempotents $\mathcal{E}(S_\Lambda)$ admits a top element with respect to $\le$ if and only if $|\Lambda^0|<\infty$. In this case, $\Delta_{\Lambda^0}$ is the top element.
\end{cor}
\begin{proof}
The first part is immediate from Lemma \ref{lem idempotents of the semigroup} and Proposition \ref{pro natural partial order on S_Lambda}. 

Let $|\Lambda^0|<\infty$ and $F\in \mathcal{E}(S_\Lambda)$. Since for any $\lambda\in \mathsf{D}(F)$, $r(\lambda)\in \mathsf{D}(\Delta_{\Lambda^0})$ and $\lambda=r(\lambda)\lambda$, it follows from the first part that $F\le \Delta_{\Lambda^0}$. This implies that $\Delta_{\Lambda^0}$ is the top element of $\mathcal{E}(S_\Lambda)$ with respect to the natural partial order. 

Now assume that $\mathcal{E}(S_\Lambda)$ has a top element $F$. Then for any $u\in \Lambda^0$, $\{(u,u)\}\le F$. By the first part, this implies that $u=\xi\beta$ for some $\xi\in \mathsf{D}(F)$ and $\beta\in \Lambda$. But this implies $\xi=\beta=u$. Therefore, $(u,u)=(\xi,\xi)\in F$ for all $u\in \Lambda^0$. Since $|F|<\infty$, it follows that $|\Lambda^0|<\infty$. 
\end{proof}
With the help of the above, we can easily describe the partial orders $\le_{\mathscr{L}}$ and $\le_{\mathscr{R}}$ defined respectively on the set of equivalence classes of the Green's relations $\mathscr{L}$ and $\mathscr{R}$. 
\begin{prop}\label{pro partial orders on equivalence classes L,R}
Let $F,G\in S_\Lambda$. Then

$(i)$ $[F]_{\mathscr{L}}\le_{\mathscr{L}} [G]_{\mathscr{L}}$ if and only if for each $\mu\in \mathsf{R}(F)$, there exist $\eta\in \mathsf{R}(G)$ and $\beta\in \Lambda$ such that $\mu=\eta\beta$, and

$(ii)$ $[F]_{\mathscr{R}}\le_{\mathscr{R}} [G]_{\mathscr{R}}$ if and only if for each $\lambda\in \mathsf{D}(F)$, there exist $\xi\in \mathsf{D}(G)$ and $\beta\in \Lambda$ such that $\lambda=\xi\beta$.
\end{prop}
\begin{proof}
By Lemma \ref{lem L and R on inverse semigroups}, $[F]_{\mathscr{L}}\le_{\mathscr{L}} [G]_{\mathscr{L}}$ if and only if $F^{-1}F\le G^{-1}G$, where $\le$ is the natural partial order. Following the proof of Proposition \ref{pro the relations L,R,H} $(i)$, we have $F^{-1}F=\Delta_{\mathsf{R}(F)}$ and $G^{-1}G=\Delta_{\mathsf{R}(G)}$. Hence $\mathsf{D}(F^{-1}F)=\mathsf{R}(F)$ and $\mathsf{D}(G^{-1}G)=\mathsf{R}(G)$. Now $(i)$ follows from the first part of Corollary \ref{cor partial order on semilattice}. Since $[F]_{\mathscr{R}}\le_{\mathscr{R}} [G]_{\mathscr{R}}$ if and only if $FF^{-1}\le GG^{-1}$, the proof of $(ii)$ is equally straightforward.
\end{proof}

With the complete description of the natural partial order on $S_\Lambda$ in hand, we now turn our attention to covers of elements of $S_\Lambda$. First we fix some notations. Suppose $F\in S_\Lambda$ and $\mathcal{Z}$ is a cover of $F$. For each $(\lambda,\mu)\in F$ and $G\in \mathcal{Z}$, we write \[G_{(\lambda,\mu)}:=\{(x,y)\in G~|~x=\lambda \alpha~\text{and}~y=\mu\alpha~\text{for some}~\alpha\in \Lambda\}.\] Also, for any $(\lambda,\mu)\in F$, we set \[\mathcal{Z}_{(\lambda,\mu)}:=\{G_{(\lambda,\mu)}~|~G\in \mathcal{Z}\}.\]

By considering $S_\Lambda$ as the inverse semigroup of the self-similar $k$-graph $(\{e\},\Lambda)$, parts $(ii)$ and $(iii)$ of the following proposition can be compared with the statement mentioned in \cite[Remark 4.5]{Larki} about covers of idempotents. However, we feel that there should exist a detailed proof of this, which was missing in \cite{Larki}, and we provide it for covers of arbitrary elements of $S_\Lambda$.

\begin{prop}\label{pro remarks on covers of an element}
Let $\Lambda$ be a finitely aligned $k$-graph and $\emptyset\neq F\in S_\Lambda$. Then the following hold.

$(i)$ If $F=\displaystyle{\bigcup_{i=1}^{n}}~F_i$, where $F_i\in S_\Lambda$, $i=1,2,\ldots,n$, then $\mathcal{C}:=\{F_i~|~i=1,2,\ldots,n\}$ is a cover of $F$.

$(ii)$ For any finite cover $\mathcal{Z}$ of $F$ and $(\lambda,\mu)\in F$, $\mathcal{Z}_{(\lambda,\mu)}$ is a finite cover of $\{(\lambda,\mu)\}$.

$(iii)$ For any finite cover $\mathcal{Z}$ of $F$ and $G\in \mathcal{Z}$, $G=\displaystyle{\bigsqcup_{(\lambda,\mu)\in F}}G_{(\lambda,\mu)}$. 
\end{prop}
\begin{proof}
$(i)$ Suppose $\emptyset \neq H\le F$. Choose any $(\lambda,\mu)\in H$. By Proposition \ref{pro natural partial order on S_Lambda}, there exist $(\xi,\eta)\in F$ and $\beta\in \Lambda$ such that $\lambda=\xi\beta$ and $\mu=\eta\beta$. Since $F=\displaystyle{\bigcup_{i=1}^{n}}~F_i$, $(\xi,\eta)\in F_j$ for some $j\in \{1,2,\ldots,n\}$. It follows that $\{(\lambda,\mu)\}\le H,F_j$. Hence $\{(\lambda,\mu)\}\le H\wedge F_j$, which shows that $H\wedge F_j\neq \emptyset$. Therefore, $\mathcal{C}$ is a finite cover of $F$. 

$(ii)$ For any $G\in \mathcal{Z}$ and $(\lambda,\mu)\in F$, since $G_{(\lambda,\mu)}\subseteq G$, $G_{(\lambda,\mu)}\in S_\Lambda$. By definition, $G_{(\lambda,\mu)}\le \{(\lambda,\mu)\}$ and so $\mathcal{Z}_{(\lambda,\mu)}\subseteq \{(\lambda,\mu)\}^{\downarrow}$. Let $\emptyset \neq H\le \{(\lambda,\mu)\}$. Then $H\le F$ and since $\mathcal{Z}$ is a cover of $F$, there exists $G\in \mathcal{Z}$ such that $H\wedge G\neq \emptyset$. Choose any $(x,y)\in H\wedge G$. Then $\{(x,y)\}\le H\le \{(\lambda,\mu)\}$, which implies that $x=\lambda\alpha$ and $y=\mu\alpha$ for some $\alpha\in \Lambda$. Again since $\{(x,y)\}\le G$, there are $(\xi,\eta)\in G$ and $\beta\in \Lambda$ such that $x=\xi\beta$ and $y=\eta\beta$. Now $G\le F$ and so we have $(\gamma,\delta)\in F$ and $\nu\in \Lambda$ such that $\xi=\gamma\nu$ and $\eta=\delta\nu$. Combining all these, we get $\lambda\alpha=\gamma\nu\beta$ and $\mu\alpha=\delta\nu\beta$. But since distinct pairs in $F$ are componentwise orthogonal, we must have $(\lambda,\mu)=(\gamma,\delta)$. It follows that $(\xi,\eta)\in G_{(\lambda,\mu)}$, whence $(\xi,\eta)\le G_{(\lambda,\mu)}$. Now we have $\{(x,y)\}\le H,G_{(\lambda,\mu)}$. This implies that $H\wedge G_{(\lambda,\mu)}\neq \emptyset$ and thus proves $(ii)$. 

$(iii)$ Let $(x,y)\in G$. Since $G\le F$, there exist $(\lambda,\mu)\in F$ and $\beta\in \Lambda$ such that $x=\lambda\beta$ and $y=\mu\beta$. Hence $(x,y)\in G_{(\lambda,\mu)}$. This proves the inclusion $G\subseteq \displaystyle{\bigcup_{(\lambda,\mu)\in F}}G_{(\lambda,\mu)}$. The other inclusion is obvious. Hence we have $G=\displaystyle{\bigcup_{(\lambda,\mu)\in F}}G_{(\lambda,\mu)}$. To show that the union is a disjoint union, we work by contradiction. Suppose $(\lambda,\mu)$ and $(\gamma,\delta)$ are distinct pairs in $F$ with $G_{(\lambda,\mu)}\cap G_{(\gamma,\delta)}\neq \emptyset$ and choose an element $(x,y)$ from this intersection. Then there exist $\alpha,\beta\in \Lambda$ such that $x=\lambda\alpha=\gamma\beta$, $y=\mu\alpha=\delta\beta$. Then $\Lambda^{\min}(\lambda,\gamma),\Lambda^{\min}(\mu,\delta)\neq \emptyset$, which contradicts that $(\lambda,\mu)$ and $(\gamma,\delta)$ are componentwise orthogonal. This completes the proof.  
\end{proof}

\begin{rmk}\label{rem an alternate proof of inner cover}
The reader may note that Part $(i)$ of Proposition \ref{pro remarks on covers of an element} can also be proved by using Proposition \ref{pro natural partial order and set-inclusion} of the next section.\qed   
\end{rmk} 

By using Proposition \ref{pro remarks on covers of an element}, we now provide a concrete and simpler description of the Cuntz--Krieger ideal of $\mathsf{B}(S_\Lambda)$. This will help us to unveil the relationship between \emph{abstract} and \emph{concrete} Cuntz--Krieger relations in the context of inverse semigroups of higher-rank graphs and their discrete algebras (see Remark \ref{rem abstract CK relations and X-to-join}). 

\begin{thm}\label{th description of Cuntz-Krieger ideal}
Let $\Lambda$ be any finitely aligned $k$-graph and $\mathscr{I}_{\textup{CK}}$ the Cuntz--Krieger ideal of $\mathsf{B}(S_\Lambda)$. Then $\mathscr{I}_{\textup{CK}}$ is the additive ideal of $\mathsf{B}(S_\Lambda)$ generated by the relations:

\textup{(I)} $F\cup G=F\vee G$, for all $F,G\in \mathcal{E}(S_\Lambda)$ with $F\cup G\in S_\Lambda$ and $F\cap G=\emptyset$;

\textup{(II)} $\{(v,v)\}=\displaystyle{\bigvee_{\tau\in X}} \{(\tau,\tau)\}$, for all $v\in \Lambda^0$ and finite exhaustive subset $X$ of $v\Lambda$.
\end{thm}
\begin{proof}
Let $\mathcal{J}$ denote the additive ideal of $\mathsf{B}(S_\Lambda)$ generated by (I) and (II). Let $F,G\in \mathcal{E}(S_\Lambda)$ be such that $F\cup G\in S_\Lambda$ and $F\cap G=\emptyset$. Then by Proposition \ref{pro remarks on covers of an element}$(i)$, $\{F,G\}$ is a cover of $F\cup G$ in $S_\Lambda$ and so 
\begin{equation}\label{eq I in CK}
(F\cup G)\setminus (F\vee G)\in \mathscr{I}_{\textup{CK}}.    
\end{equation}
Now let $v\in \Lambda^0$ and $X$ a finite exhaustive subset of $v\Lambda$. Let $\mathcal{C}=\left\{\{(\tau,\tau)\}~|~\tau\in X\right\}$. By Corollary \ref{cor partial order on semilattice}, $\mathcal{C}\subseteq \{(v,v)\}^{\downarrow}$. If $\emptyset\neq H\le \{(v,v)\}$, then choosing any $(x,x)\in H$, we have $x\in v\Lambda$. Since $X$ is exhaustive, there is a $\lambda\in X$ such that $\MCE(x,\lambda)\neq \emptyset$. Now choosing any $\gamma\in \MCE(x,\lambda)$, it follows that $\{(\gamma,\gamma)\}\le \{(x,x)\}\le H$ and $\{(\gamma,\gamma)\}\le \{(\lambda,\lambda)\}$. So $H\wedge \{(\lambda,\lambda)\}\neq \emptyset$. This shows that $\mathcal{C}$ is a cover of $\{(v,v)\}$. Hence 
\begin{equation}\label{eq II in CK}
\{(v,v)\}\setminus \left(\displaystyle{\bigvee_{\tau\in X}}\{(\tau,\tau)\}\right)\in \mathscr{I}_{\textup{CK}}.   
\end{equation}
From (\ref{eq I in CK}) and (\ref{eq II in CK}), it follows that $\mathcal{J}\subseteq \mathscr{I}_{\textup{CK}}$.
 
For the other inclusion, we need to show that any relation on $\mathsf{B}(S_\Lambda)$ of the form (\ref{eq ACK}) follows from the relations (I) and (II). In view of \cite[Lemma 11.1]{LawsonCK} (see \cite[Lemma 3.1(1)]{LawNC} for a proof), it suffices to prove the relations (\ref{eq ACK}) only for idempotents of $S_\Lambda$. So assume $F\in \mathcal{E}(S_\Lambda)$ and $\mathcal{Z}$ is any finite cover of $F$. Let $\lambda\in \mathsf{D}(F)$. Then by Proposition \ref{pro remarks on covers of an element}$(ii)$, $\mathcal{Z}_{(\lambda,\lambda)}$ is a finite cover of $\{(\lambda,\lambda)\}$. Applying \cite[Lemma 4.6(2)]{Larki}, it follows that \[\left\{\{(\tau,\tau)\}~|~(\lambda\tau,\lambda\tau)\in G_{(\lambda,\lambda)}~\text{for some}~G\in \mathcal{Z}\right\}\] is a finite cover of $\{(s(\lambda),s(\lambda))\}$ and consequently, by \cite[Lemma 4.6(3)]{Larki}, the set \[X_\lambda:=\{\tau~|~(\lambda\tau,\lambda\tau)\in G_{(\lambda,\lambda)}~\text{for some}~G\in \mathcal{Z}\}\] is a finite exhaustive subset of $s(\lambda)\Lambda$. Then by (II), we have 
\begin{equation}\label{eq CKS1}
\{(s(\lambda),s(\lambda))\}=\displaystyle{\bigvee_{\tau\in X_\lambda}}\{(\tau,\tau)\}.   
\end{equation}
Multiplying both sides of (\ref{eq CKS1}) by $\{(\lambda,s(\lambda))\}$ from the left and by $\{(s(\lambda),\lambda)\}$ from the right and using the fact that in a Boolean inverse semigroup, multiplication distributes over finite compatible joins, we get 
\begin{align}\label{eq CKS2}
\{(\lambda,\lambda)\}&=\{(\lambda,s(\lambda))\}\left(\displaystyle{\bigvee_{\tau\in X_\lambda}}\{(\tau,\tau)\}\right)\{(s(\lambda),\lambda)\} \nonumber \\ 
&=\displaystyle{\bigvee_{\tau\in X_\lambda}}\left((\lambda,s(\lambda))\cdot (\tau,\tau)\cdot (s(\lambda),\lambda)\right) \nonumber\\
&=\displaystyle{\bigvee_{\tau\in X_\lambda}}\{(\lambda\tau,\lambda\tau)\}=\displaystyle{\bigvee_{G\in \mathcal{Z}}}\left(\displaystyle{\bigvee_{x\in \mathsf{D}(G_{(\lambda,\lambda)})}}\{(x,x)\}\right)=\displaystyle{\bigvee_{G\in \mathcal{Z}}}G_{(\lambda,\lambda)}.
\end{align}
The last equality above follows from relations (I). Finally by using relation (I), (\ref{eq CKS2}) and Proposition \ref{pro remarks on covers of an element}$(iii)$, we obtain 
\begin{align*}
F=\displaystyle{\bigvee_{\lambda\in \mathsf{D}(F)}}\{(\lambda,\lambda)\}=\displaystyle{\bigvee_{\lambda\in \mathsf{D}(F)}}\left(\displaystyle{\bigvee_{G\in \mathcal{Z}}} G_{(\lambda,\lambda)}\right)=\displaystyle{\bigvee_{G\in \mathcal{Z}}}\left(\displaystyle{\bigvee_{\lambda\in \mathsf{D}(F)}}G_{(\lambda,\lambda)}\right)=\displaystyle{\bigvee_{G\in \mathcal{Z}}}~G, 
\end{align*}
which is the desired relation. Therefore, $\mathscr{I}_{\textup{CK}}\subseteq \mathcal{J}$ and we are done. 
\end{proof}
The following remark is worth mentioning.
\begin{rmk}\label{rem abstract CK relations and X-to-join}
If $F,G\in S_\Lambda$ are such that $F\cup G\in S_\Lambda$ and $F\cap G=\emptyset$, then $F,G\le F\cup G$ by Proposition \ref{pro natural partial order on S_Lambda}. Clearly, $F$ and $G$ are compatible in $\mathsf{B}(S_\Lambda)$ and we have $F\vee G\le F\cup G$. We have no reason to believe that $F\cup G\le F\vee G$ holds in $\mathsf{B}(S_\Lambda)$. 


Looking at the Cuntz--Krieger relations \[\displaystyle{\prod_{\tau\in X}}(s_v-s_{\lambda}s_{\lambda}^*)=0\] that define a Cuntz--Krieger $\Lambda$-family (see \cite{Farthing,Yeend}) and the extra relations (I) of Theorem \ref{th description of Cuntz-Krieger ideal}, we can say that the abstract Cuntz--Krieger relations on $\mathsf{B}(S_\Lambda)$ do not coincide immediately with all and only the concrete Cuntz--Krieger relations.

Let \[X_{\cup}:=\left\{\left(F,\{F_1,F_2,\ldots,F_n\}\right)\in \mathcal{E}(S_\Lambda)\times \mathcal{P}_{fin}(\mathcal{E}(S_\Lambda))~|~F=\displaystyle{\bigsqcup_{i=1}^{n}}~F_i\right\}\] and consider the \emph{$X_{\cup}$-to-join Booleanization} \cite{GannaX} $\mathsf{B}_{X_{\cup}}(S_\Lambda)$ of $S_\Lambda$. By \cite[Corollary 5.4$(2)$]{GannaX}, it follows that $\mathsf{B}_{X_{\cup}}(S_\Lambda)\cong \mathsf{B}(S_\Lambda)/{\mathcal{J}_{\textup{I}}}$, where $\mathcal{J}_{\textup{I}}$ is the additive ideal generated by relations (I). Then the abstract Cuntz--Krieger relations on $\mathsf{B}(S_\Lambda)$ essentially reduce to the relations (II) on $\mathsf{B}_{X_{\cup}}(S_\Lambda)$, and hence, this particular $X$-to-join Booleanization of $S_\Lambda$ gives the perfect setting for the agreement of abstract and concrete Cuntz--Krieger relations. 

One can see that the set $X_{\cup}$ is an $S_\Lambda$-invariant subset \cite[Definition 3.7]{GannaX} of $\mathcal{E}(S_\Lambda)\times \mathcal{P}_{fin}(\mathcal{E}(S_\Lambda))$. Indeed, if $F=\displaystyle{\bigsqcup_{i=1}^{n}}~F_i$, then a technical computation yields that $G\in S_\Lambda$, $G^{-1}FG=\displaystyle{\bigsqcup_{i=1}^{n}}~G^{-1}F_i G$ for any $G\in S_\Lambda$. Thus, $\widehat{\mathcal{E}(S_\Lambda)}_{X_\cup}$, the set of all $X_\cup$-to-join characters on $\mathcal{E}(S_\Lambda)$, is a closed invariant subset of the whole character space $\widehat{\mathcal{E}(S_\Lambda)}$. Consequently, $\mathsf{B}_{X_\cup}(S_\Lambda)$ is the dual Boolean inverse semigroup of the ample groupoid $\mathcal{G}(S_\Lambda)|_{\widehat{\mathcal{E}(S_\Lambda)}_{X_\cup}}$.\qed
\end{rmk}

\section{Graded core and the relationship with Kumjian--Pask algebras}\label{sec Relation to KP-algebra}
As the elements of the inverse semigroup $S_\Lambda$ are finite sets, apart from the natural partial order, one can also equip this semigroup with the partial order of set-inclusion. If $F,G\in S_\Lambda$ are such that $F\subseteq G$, then from Proposition \ref{pro natural partial order on S_Lambda} it follows that $F\le G$, where $\le$ is the natural partial order on $S_\Lambda$. Obviously, the converse is not true as one may find many examples; for an instance, if $\lambda\in \Lambda\setminus \Lambda^0$, then $\{(\lambda,\lambda)\}\le \{(r(\lambda),r(\lambda)\}$ but $\{(\lambda,\lambda)\}\nsubseteq \{(r(\lambda),r(\lambda)\}$.

\begin{rmk}\label{rem set-inclusion for graph inverse semigroup}
If $\Gamma$ is a directed graph and we view the graph inverse semigroup $S(\Gamma)$ as a subsemigroup of $S_{\Gamma^*}$ via the canonical embedding $\psi_\Gamma$ of Lemma \ref{lem the relation between gis and S_Lambda}, then the partial order $\subseteq$ of $S_{\Gamma^*}$ reduces to the equality relation on $S(\Gamma)$. As a consequence, two distinct elements $\alpha\beta^{-1}$ and $\gamma\delta^{-1}$ of $S(\Gamma)$ cannot have a join with respect to $\subseteq$ in $S(\Gamma)$. However, they may have a join in the bigger semigroup $S_{\Gamma^*}$ (see Proposition \ref{pro natural partial order and set-inclusion} below). \qed      
\end{rmk}

We now record a useful proposition, which tells us precisely when two elements of $S_\Lambda$ admit a join with respect to $\subseteq$.   
\begin{prop}\label{pro natural partial order and set-inclusion}
Let $\Lambda$ be a finitely aligned $k$-graph and $F,G\in S_\Lambda$. Then $F,G$ have a join with respect to $\subseteq$ if and only if $F\cup G\in S_\Lambda$. Moreover, in this case, both the join and meet of $F$ and $G$ with respect to the natural partial order exist and are equal to $F\cup G$ and $F\cap G$ respectively.
\end{prop}
\begin{proof}
Let $H$ be the join of $F,G$ with respect to $\subseteq$. Then $F\subseteq H$ and $G\subseteq H$. This implies $F\cup G\subseteq H$. Since $H\in S_\Lambda$, it follows that for any two distinct pairs $(\lambda,\mu)$ and $(\xi,\eta)$ in $F\cup G$, $(\lambda,\mu)\perp (\xi,\eta)$. Therefore, $F\cup G\in S_\Lambda$. The converse is obvious. 

For the final statement let $F\cup G\in S_\Lambda$. By Proposition \ref{pro natural partial order on S_Lambda}, $F,G\le F\cup G$. Let $K\in S_\Lambda$ be such that $F,G\le K$. Let $(\lambda,\mu)\in F\cup G$. If $(\lambda,\mu)\in F$, then since $F\le K$, there exist $(\gamma,\delta)\in K$ and $\beta\in \Lambda$ such that $\lambda=\gamma\beta$ and $\mu=\delta\beta$. Similar thing happens if $(\lambda,\mu)\in G$ since $G\le K$. Therefore, $F\cup G\le K$ by Proposition \ref{pro natural partial order on S_Lambda}. This shows that $F\cup G$ is the join of $F$ and $G$ with respect to $\le$. For the meet, observe that $F\cap G$ is obviously an element of $S_\Lambda$ and $F\cap G\le F,G$. Now let $\emptyset\neq H\in S_\Lambda$ be such that $H\le F,G$. Let $(x,y)\in H$. Then there exist $(\xi,\eta)\in F$, $(\tau,\nu)\in G$ and $\alpha,\beta\in \Lambda$ such that $x=\xi\alpha=\tau\beta$ and $y=\eta\alpha=\nu\beta$. It follows that $\Lambda^{\min}(\xi,\tau),\Lambda^{\min}(\eta,\nu)\neq \emptyset$. Since $F\cup G\in S_\Lambda$ and $(\xi,\eta),(\tau,\nu)\in F\cup G$, we must have $(\xi,\eta)=(\tau,\nu)$, whence $(\xi,\eta)\in F\cap G$. Applying Proposition \ref{pro natural partial order on S_Lambda} again, it follows that $H\le F\cap G$ and hence $F\cap G$ is the meet of $F$ and $G$ with respect to $\le$.   
\end{proof}

The usefulness of the above lemma lies in the second part of the statement as it provides a sufficient condition for two arbitrary elements of $S_\Lambda$ to have a join with respect to the natural partial order $\le$. The following example makes this even more apparent.
\begin{example}\label{ex two elements with no join}
Suppose $\Lambda$ is a $2$-graph whose $1$-skeleton contains the following zigzag of paths as a subgraph. 
\[
\begin{tikzpicture}[
    vertex/.style={circle, fill=black, inner sep=2pt},
    bluearrow/.style={-{Stealth[length=5pt]}, thick, blue},
    redarrow/.style={-{Stealth[length=5pt]}, thick, red, dashed},
    every node/.style={font=\large}
]

\node[vertex, label=above:\normalsize{$a$}] (a) at (1.25,1.25) {};
\node[vertex, label=below:\normalsize{$u$}] (u) at (0,0) {};
\node[vertex, label=below:\normalsize{$v$}] (v) at (2.5,0) {};

\node[vertex, label=above:\normalsize{$b$}] (b) at (3.75,1.25) {};
\node[vertex, label=below:\normalsize{$w$}] (w) at (5,0) {};

\node[vertex] (x) at (-1.25,1.25) {};
\node[vertex] (y) at (6.25,1.25) {};

\draw[bluearrow] (a) -- node[left=2pt] {\normalsize{\textcolor{black}{$\alpha$}}} (u);
\draw[redarrow]  (a) -- node[right=2pt] {\normalsize{\textcolor{black}{$\beta$}}} (v);

\draw[bluearrow] (b) -- node[left=2pt] {\normalsize{\textcolor{black}{$\delta$}}} (v);
\draw[redarrow]  (b) -- node[right=2pt] {\normalsize{\textcolor{black}{$\gamma$}}} (w);

\draw[redarrow]  (x) -- node[right=2pt] {} (u);
\draw[bluearrow]  (y) -- node[right=2pt] {} (w);

\end{tikzpicture}
\]
Let $F:=\{(\alpha,\beta)\}$ and $G:=\{(\gamma,\delta)\}$. If $\Lambda^{\min}(\beta,\delta)=\emptyset$, then $F\cup G\in S_\Lambda$ and by Proposition \ref{pro natural partial order and set-inclusion}, it is the join of $F$ and $G$ with respect to both the partial orders $\subseteq$ and $\le$. On the other hand, if $\Lambda^{\min}(\beta,\delta)\neq \emptyset$, then $F\cup G\notin S_\Lambda$. We observe that $F$ and $G$ have no common upper bound (and hence no join) with respect to $\le$. For if $H\in S_\Lambda$ is such that $F,G\le H$, then there exist $(\lambda,\mu), (\xi,\eta)\in H$ and $x,y\in \Lambda$ such that $\alpha=\lambda x$, $\beta=\mu x$, $\gamma=\xi y$ and $\delta=\eta y$. But then choosing any $(\nu,\tau)\in \Lambda^{\min}(\beta,\delta)$, we have \[\mu x \nu=\beta \nu=\delta \tau=\eta y \tau,\] which implies $\Lambda^{\min}(\mu,\eta)\neq \emptyset$ and consequently, $(\lambda,\mu)=(\xi,\eta)$. Then $u=r(\alpha)=r(\lambda)=r(\xi)=r(\gamma)=w$, which is clearly not the case. 
\end{example}

The partial order $\subseteq$ is clearly compatible with the multiplication and the inversion of $S_\Lambda$, and hence we can consider $(S_\Lambda,\subseteq)$ as an ordered inverse semigroup. The next proposition identifies the join-closed ideals of this ordered semigroup. 

\begin{prop}\label{pro u-closed ideals are join closed}
Let $\Lambda$ be a finitely aligned $k$-graph and $\mathcal{J}$ an ideal of $S_\Lambda$. Then $\mathcal{J}$ is a join-closed ideal of $(S_\Lambda,\subseteq)$ if and only if for $F,G\in \mathcal{J}$, $F\cup G\in S_\Lambda$ implies $F\cup G\in \mathcal{J}$.    
\end{prop}
\begin{proof}
If $\mathcal{J}$ is a join-closed ideal, then obviously it has the said property. For the other direction, we just need to observe that under the condition described in the statement, $\mathcal{J}$ is an order ideal with respect to $\subseteq$. Let $F\in \mathcal{J}$ and $G\in S_\Lambda$ be such that $G\subseteq F$. Since $\emptyset\in \mathcal{J}$, we may assume $G$ to be nonempty. Let $(\lambda,\mu)\in G$. Then $(\lambda,\mu)\in F$ and \[\{(\lambda,\mu)\}=(\lambda,\lambda)\cdot (\lambda,\mu)\cdot (\mu,\mu)\subseteq \{(\lambda,\lambda)\}F\{(\mu,\mu)\}\in \mathcal{J}\] since $\mathcal{J}$ is an ideal. Finally since $\displaystyle{\bigcup_{(\lambda,\mu)\in G}} \{(\lambda,\mu)\}=G\in S_\Lambda$ and $\mathcal{J}$ is $\cup$-closed, it follows that $G\in \mathcal{J}$. Hence $\mathcal{J}$ is an order ideal and we are done.  
\end{proof}

In view of the above proposition, we refer to the join closed ideals of $(S_\Lambda,\subseteq)$ as \emph{$\cup$-closed}. It is important to note that not every ideal of $S_\Lambda$ is $\cup$-closed. For example, consider the $2$-graph $\Lambda$ with the following $1$-skeleton: 
\[
\begin{tikzpicture}[>=stealth, line width=1.2pt]

    \fill (0,0) circle (2pt);
    \fill (0,-2) circle (2pt);
    \fill (2,-2) circle (2pt);

    \draw[->, blue] (2,-2) -- node[below] {\textcolor{black}{$f$}} (0,-2);

    \draw[->, dashed, red, line width=1.2pt]
        (0,0) -- node[left] {\textcolor{black}{$e$}} (0,-2);

    \node[] (a) at (-0.3,0) {$v$};
    \node[] (b) at (-0.2,-2.2) {$u$};
    \node[] (c) at (2.3,-2) {$w$};

\end{tikzpicture}
\]
Then it can be checked that \[\mathcal{J}:=\left\{\emptyset,\{(v,v)\},\{(w,w)\},\{(v,v),(w,w)\},\{(v,e)\},\{(e,v)\},\{(w,f)\}, \{(f,w)\},\{(e,e)\},\{(f,f)\}\right\}\] is an ideal of $S_\Lambda$ which is not $\cup$-closed since $\{(e,e)\}, \{(f,f)\}\in \mathcal{J}$ but $\{(e,e),(f,f)\}\notin \mathcal{J}$.

Instead of giving particular examples of $\cup$-closed ideals in the Farthing--Muhly--Yeend inverse semigroup of a $k$-graph, we shall present an ample source of producing such ideals. Given a $k$-graph $\Lambda$, recall that a subset $\mathsf{H}$ of $\Lambda^0$ is called \emph{hereditary} if for $\lambda\in \Lambda$, $r(\lambda)\in \mathsf{H}$ implies $s(\lambda)\in \mathsf{H}$. We now prove two lemmas which will lead us to a structure theorem for $\cup$-closed ideals of $S_\Lambda$.

\begin{lem}\label{lem ideal to hereditary subset}
Let $\Lambda$ be a finitely aligned $k$-graph and $\mathcal{J}$ an ideal of $S_\Lambda$. Let \[\mathsf{H}_\mathcal{J}:=\{v\in \Lambda^0~|~\{(v,v)\}\in \mathcal{J}\}.\] Then $\mathsf{H}_\mathcal{J}$ is a hereditary subset of $\Lambda^0$.    
\end{lem}
\begin{proof}
Let $\lambda\in \Lambda$ be such that $r(\lambda)\in \mathsf{H}_\mathcal{J}$. Then $\{(r(\lambda),r(\lambda))\}\in \mathcal{J}$. Note that \[(s(\lambda),\lambda)\cdot (r(\lambda),r(\lambda))=\{(s(\lambda)s(\lambda),r(\lambda)\lambda)\}=\{(s(\lambda),\lambda)\}\] and \[(s(\lambda),\lambda)\cdot (\lambda,s(\lambda))=\{(s(\lambda),s(\lambda))\}.\] From these it follows that \[\{(s(\lambda),s(\lambda))\}=\{(s(\lambda),\lambda)\}\{(r(\lambda),r(\lambda))\}\{(\lambda,s(\lambda))\}\in \mathcal{J}\] since $\mathcal{J}$ is an ideal. Therefore, $s(\lambda)\in \mathsf{H}_\mathcal{J}$ and so $\mathsf{H}_\mathcal{J}$ is hereditary.
\end{proof}
\begin{lem}\label{lem Hereditary subset to join closed ideal}
Let $\Lambda$ be a finitely aligned $k$-graph and $\mathsf{H}$ any hereditary subset of $\Lambda^0$. Let \[\mathcal{J}_\mathsf{H}:=\{F\in S_\Lambda~|~\mathfrak{s}(F)\subseteq \mathsf{H}\}.\] Then $\mathcal{J}_\mathsf{H}$ is a $\cup$-closed ideal of $S_\Lambda$.
\end{lem}
\begin{proof}
Suppose $F,G\in \mathcal{J}_\mathsf{H}$ are such that $F\cup G\in S_\Lambda$. Since \[\mathfrak{s}(F\cup G)=s(\mathsf{D}(F\cup G))=s(\mathsf{D}(F))\cup s(\mathsf{D}(G))=\mathfrak{s}(F)\cup \mathfrak{s}(G)\subseteq \mathsf{H},\] we have $F\cup G\in \mathcal{J}_\mathsf{H}$. Thus $\mathcal{J}_\mathsf{H}$ is $\cup$-closed. It remains to show that $\mathcal{J}_\mathsf{H}$ is indeed an ideal of $S_\Lambda$. For this, let $F\in \mathcal{J}_\mathsf{H}$ and $K\in S_\Lambda$. Suppose $v\in \mathfrak{s}(KF)$. Then $v=s(\lambda)=s(\mu)$ for some $(\lambda,\mu)\in KF$. Now there exist $(\gamma,\eta)\in K$ and $(\xi,\delta)\in F$ such that $(\lambda,\mu)\in (\gamma,\eta)\cdot (\xi,\delta)$. This implies that $\lambda=\gamma\alpha$ and $\mu=\delta\beta$ for some $(\alpha,\beta)\in \Lambda^{\min}(\eta,\xi)$. Since $F\in \mathcal{J}_\mathsf{H}$, $r(\beta)=s(\xi)\in \mathfrak{s}(F)\subseteq \mathsf{H}$. Now since $\mathsf{H}$ is hereditary, we have $v=s(\mu)=s(\beta)\in \mathsf{H}$. Therefore, $\mathfrak{s}(KF)\subseteq \mathsf{H}$, whence $KF\in \mathcal{J}_\mathsf{H}$. A similar line of argument can be employed to show that $FK\in \mathcal{J}_\mathsf{H}$. Hence the lemma follows. 
\end{proof}

Let us denote the set of all hereditary subsets of $\Lambda^0$ by $\mathscr{H}er(\Lambda)$ and the set of all $\cup$-closed ideals of $S_\Lambda$ by $\mathscr{U}id(S_\Lambda)$. Both these sets are partially ordered by set-inclusion. We now present the following result, which connects these two posets in a nice way.

\begin{thm}\label{th structure thm for U-closed ideals}
Let $\Lambda$ be a finitely aligned $k$-graph. Then the map 
\begin{align*}
\phi: \mathscr{H}er(\Lambda)&\longrightarrow \mathscr{U}id(S_\Lambda)\\
\mathsf{H} &\longmapsto \mathcal{J}_\mathsf{H}
\end{align*} 
is an order isomorphism with inverse given by
\begin{align*}
\psi: \mathscr{U}id(S_\Lambda) &\longrightarrow \mathscr{H}er(\Lambda)\\
\mathcal{J} &\longmapsto \mathsf{H}_\mathcal{J}.
\end{align*}
\end{thm}
\begin{proof}
By Lemmas \ref{lem ideal to hereditary subset} and \ref{lem Hereditary subset to join closed ideal}, the maps $\phi$ and $\psi$ are well defined. Again, it is easy to see that both $\phi$ and $\psi$ preserve the order. Therefore, we only need to show that $\phi$ and $\psi$ are inverses of each other. 

Let $\mathsf{H}$ be a hereditary subset of $\Lambda^0$. Let $v\in \mathsf{H}_{\mathcal{J}_\mathsf{H}}$. Then $F:=\{(v,v)\}\in \mathcal{J}_\mathsf{H}$. But this implies that $v\in \mathfrak{s}(F)\subseteq \mathsf{H}$. Thus $\mathsf{H}_{\mathcal{J}_\mathsf{H}}\subseteq \mathsf{H}$. Conversely, if $u\in \mathsf{H}$, then $\{(u,u)\}\in \mathcal{J}_\mathsf{H}$ and consequently, $u\in \mathsf{H}_{\mathcal{J}_\mathsf{H}}$. Hence $\mathsf{H}\subseteq \mathsf{H}_{\mathcal{J}_\mathsf{H}}$ and we have $\psi\circ \phi=\id_{\mathscr{H}er(\Lambda)}$.

Now suppose $\mathcal{J}$ is any $\cup$-closed ideal of $S_\Lambda$. Let $\emptyset \neq F\in \mathcal{J}$ and $u\in \mathfrak{s}(F)$. Choose any $(\lambda,\mu)\in F$ such that $u=s(\lambda)=s(\mu)$. Now 
\begin{align*}
\{(s(\lambda),\lambda)\}F\{(\mu,s(\mu))\}&=\displaystyle{\bigcup_{(\alpha,\beta)\in F}}(s(\lambda),\lambda)\cdot (\alpha,\beta)\cdot (\mu,s(\mu)\\
&=(s(\lambda),\lambda)\cdot (\lambda,\mu)\cdot (\mu,s(\mu))\\
&=\{(s(\lambda),s(\mu))\}=\{(u,u)\},
\end{align*}
where the second equality follows since $\Lambda^{\min}(\lambda,\alpha)=\emptyset$ for all $\alpha\in \mathsf{D}(F)\setminus \{\lambda\}$ and $\Lambda^{\min}(\beta,\mu)=\emptyset$ for all $\beta\in \mathsf{R}(F)\setminus \{\mu\}$. Since $\mathcal{J}$ is an ideal, it follows that $\{(u,u)\}\in \mathcal{J}$, whence $u\in \mathsf{H}_\mathcal{J}$. This implies that $\mathfrak{s}(F)\subseteq \mathsf{H}_\mathcal{J}$ and so $F\in \mathcal{J}_{\mathsf{H}_\mathcal{J}}$. Therefore, we have the inclusion $\mathcal{J}\subseteq \mathcal{J}_{\mathsf{H}_\mathcal{J}}$. For the other containment, let $\emptyset\neq G\in \mathcal{J}_{\mathsf{H}_\mathcal{J}}$. Then for any $(\xi,\eta)\in G$, $s(\xi)=s(\eta)\in \mathfrak{s}(G)\subseteq \mathsf{H}_\mathcal{J}$. This implies that $\{(s(\lambda),s(\mu))\}\in \mathcal{J}$. Since $\mathcal{J}$ is an ideal, we have \[\{(\lambda,\mu)\}=\{(\lambda,s(\lambda))\}\{(s(\lambda),s(\mu))\}\{(s(\mu),\mu)\}\in \mathcal{J}\] and since $\mathcal{J}$ is $\cup$-closed, $G=\displaystyle{\bigcup_{(\lambda,\mu)\in G}} \{(\lambda,\mu)\}\in \mathcal{J}$. Hence $\mathcal{J}_{\mathsf{H}_\mathcal{J}}\subseteq \mathcal{J}$ and $\phi\circ \psi=\id_{\mathscr{U}id(S_\Lambda)}$. 
\end{proof}

The graph inverse semigroup $S(\Gamma)$ of any directed graph $\Gamma$ is equipped with a natural $\mathbb{Z}$-grading (see \cite[\S 8]{HazMes}), which is induced by lengths of paths and corresponds to the natural $\mathbb{Z}$-grading of the Leavitt path algebra $L_\mathsf{k}(\Gamma)$, when the latter is viewed as the quotient of the contracted semigroup algebra of $S(\Gamma)$ modulo the (CK2) relations. 

Due to the larger size of the Farthing--Muhly--Yeend inverse semigroup of a $k$-graph $\Lambda$ compared to a graph inverse semigroup, the geometric grading of graph inverse semigroups cannot be extended directly to $S_\Lambda$ in order to make it a graded inverse semigroup. Interestingly, the ordered structure of $S_\Lambda$ with respect to $\subseteq$ helps us to realize $S_\Lambda$ as a different kind of \emph{graded} semigroups, which we define below.

\begin{dfn}\label{def graded-ordered semigroup}
Let $G$ be an abelian group and $(S,\leqslant)$ an ordered semigroup with $0$ such that $0$ is the least element with respect to $\leqslant$. We say that $S$ is \emph{graded-ordered by $G$} if there exist join-closed subsets $S_x$ of $S$, $x\in G$ such that 

$(i)$ $0\in S_x$ for all $x\in G$;

$(ii)$ $S_x S_y\subseteq S_{xy}$ for all $x,y\in G$;

$(iii)$ every element $0\neq s\in S$ can be written uniquely as $s=\displaystyle{\bigvee_{x\in G}}~s_x$, where for each $x\in G$, $s_x\in S_x$ and $s_x=0$ for all but finitely many $x\in G$.

An ordered inverse semigroup $(S,\leqslant)$ with $0$ is said to be \emph{graded-ordered by $G$} if in addition to the above three conditions, we also have $S_x^{-1}\subseteq S_{x^{-1}}$ for all $x\in G$.   
\end{dfn}

Before giving some examples of graded-ordered semigroups, we record some basic properties about their components. 

\begin{prop}\label{lem basics of graded core of S_Lambda}
Let $(S,\leqslant)$ be an ordered semigroup graded-ordered by $G$. Then the following hold.

$(i)$ For distinct $x,y\in G$, $0\neq s\in S_x$ and $0\neq t\in S_y$, $s~||~t$. 

$(ii)$ For distinct $x,y\in G$, $S_x\cap S_y=\{0\}$.

$(iii)$ Each $S_x$ is a down-set with respect to $\leqslant$.
\end{prop}
\begin{proof}
$(i)$ Let $x,y$ be distinct elements of $G$, $0\neq s\in S_x$ and $0\neq t\in S_y$. If $s\le t$, then we can write $t$ in two different ways, namely $t=0\vee t$ (viewing $0\in S_x$ and $t\in S_y$) and $t=s\vee t$ (viewing $s\in S_x$ and $t\in S_y$), which contradicts Condition $(iii)$ of Definition \ref{def graded-ordered semigroup}. So $s\nleq t$. By a similar argument, $t\nleq s$ and hence $s$ and $t$ are incomparable.

$(ii)$ This directly follows from $(i)$.

$(iii)$ Let $s\in S_x$ and $t\in S$ be such that $t\le s$. If $t=0$, then $t\in S_x$ by Condition $(i)$ of Definition \ref{def graded-ordered semigroup}. If $t\neq 0$, then by Condition $(iii)$, we can write $t=\displaystyle{\bigvee_{y\in G}}t_y$, where $t_y\in S_y$ for all $y\in G$ and only finitely many $t_y\neq 0$. Now $s=s\vee t$, together with the uniqueness mentioned in Condition $(iii)$, implies that $t_y=0$ for all $y\neq x$. It follows that $t=t_x\in S_x$. 
\end{proof}
\begin{examples}\label{ex graded-ordered semigroup}
$(i)$ Let $S=\mathbb{N}[t]$ be the monoid of all polynomials with nonnegative integral coefficients under usual multiplication of polynomials. For $f(t),g(t)\in S$, define a relation $\leqslant$ as follows: $f(t)\leqslant g(t)$ if $g(t)=f(t)+h(t)$ for some $h(t)\in S$. Then $\leqslant$ is a partial order on $S$ and the zero polynomial is the least element with respect to $\leqslant$. Any two elements of $S$ have a join with respect to $\leqslant$, namely \[\left(\displaystyle{\sum_{i\in \mathbb{N}}}a_i t^i\right)\vee \left(\displaystyle{\sum_{i\in \mathbb{N}}}b_i t^i\right)=\displaystyle{\sum_{i\in \mathbb{N}}}(a_i\vee b_i)t^i.\] Now for each $n\in \mathbb{Z}$, set $S_n:=\{at^n~|~a\in \mathbb{N}\}$ if $n\in \mathbb{N}$ and $S_n:=\{0\}$ if $n<0$. One can check that Condition $(i)$--$(iii)$ of Definition \ref{def graded-ordered semigroup} are satisfied. Hence $(S,\leqslant)$ is graded-ordered by $\mathbb{Z}$.

$(ii)$ The example in $(i)$ can be generalized as follows. Suppose $(S,+,\cdot)$ is an additively cancellative zero-sum free semiring (see \cite{Golan}). Then the algebraic preorder $\le^+$ on the commutative monoid $(S,+)$ (defined by $s\le^+ t$ if and only if $t=s+q$ for some $q\in S$) is a partial order and the multiplicative reduct $(S,\cdot)$ is an ordered semigroup with respect to $\le^+$. Now if $G$ is an abelian group and $S$ is $G$-graded as a semiring (see \cite[Definition 2.4]{MNS}), then we observe that the homogeneous components $S_x$, $x\in G$ are join-closed with respect to $\le^+$. To see this, let $s,t\in S_x$ such that $s\vee t$ exists. Since $S_x$ is an additive submonoid of $(S,+)$, $s+t\in S_x$. Moreover, $s,t\le^+ s+t$, which implies $(s\vee t)\le^+(s+t)$ and so $s+t=(s\vee t)+q$ for some $q\in S$. But then for any $y\in G\setminus \{x\}$, $(s\vee t)_y+q_y=0$, which gives $(s\vee t)_y=0$ since $S$ is zero-sum free. Thus $s\vee t\in S_x$. Condition $(i)$ and $(ii)$ of Definition \ref{def graded-ordered semigroup} follow from the definition of $G$-grading. Again for Condition $(iii)$, it suffices to observe that the finite join of homogeneous elements of different degrees is nothing but their sum. For this, let $s_i\in S_{x_i}$, $i=1,2,\ldots,n$ with $x_i\neq x_j$ for $i\neq j$. Clearly $s_i\le^+ \displaystyle{\sum_{j=1}^{n}}~s_j$. If $t\in S$ is such that $s_i\le^+ t$ for all $i$, then there are $q^{(i)}\in S$ such that $t=s_i+q^{(i)}$. Then $t_{x_i}=s_i+(q^{(i)})_{x_i}$ for each $i=1,2,\ldots,n$ and \[\displaystyle{\sum_{i=1}^{n}}~s_i\le^+ \displaystyle{\sum_{i=1}^{n}}~\left(s_i+(q^{(i)})_{x_i}\right)\le^+ \displaystyle{\sum_{x\in G}}~t_x=t.\] 

$(iii)$ Let $T$ be any semigroup, $G$ an abelian group and $\phi:T\longrightarrow G$ a semigroup homomorphism. Define $S:=\mathcal{P}_{fin}(T)$. Then $S$ is a semigroup with respect to the operation \[AB:=\{ab~|~a\in A~\text{and}~b\in B\}.\] Note that $S$ is ordered by set-inclusion and the empty subset $\emptyset$ is simultaneously the zero of this semigroup and the least element with respect to the order. For each $x\in G$, define \[S_x:=\{A\in S~|~A\subseteq \phi^{-1}(x)\}.\] Using these components, it is easy to see that $(S,\subseteq)$ is graded-ordered by $G$. 

$(iv)$ Suppose $S=\displaystyle{\bigcup_{x\in G}}~S_x$ is a $G$-graded inverse semigroup with $0$ (see Subsection \ref{ssec semigroups}). Note that for any $x\in G$ and $0\neq s\le t$, $s\in S_x$ if and only if $t\in S_x$. This is because $\deg(s)=\deg(ss^{-1}t)=\deg(ss^{-1})\deg(t)=\epsilon\deg(t)=\deg(t)$, where $\epsilon$ is the identity of $G$. In view of this fact, it can be observed that each $S_x$ is join-closed with respect to $\le$ and all the conditions of Definition \ref{def graded-ordered semigroup} are satisfied by considering these components. In fact, since each nonzero element $s\in S$ lies in a unique $S_x$, the join decomposition is trivial. Hence $(S,\le)$ is graded-ordered by $G$. 
\end{examples}
The following lemma gives the most important graded-ordered inverse semigroup in line with our discussion.
\begin{lem}\label{lem S_Lambda is graded-ordered}
For any finitely aligned $k$-graph $\Lambda$, the ordered inverse semigroup $(S_\Lambda,\subseteq)$ is graded-ordered by $\mathbb{Z}^k$ with \[(S_\Lambda)_\N:=\{F\in S_\Lambda~|~d(\lambda)-d(\mu)=\N~\text{for all}~(\lambda,\mu)\in F\}\] for all $\N\in \mathbb{Z}^k$.   
\end{lem}
\begin{proof}
Note that $\emptyset\in (S_\Lambda)_\N$ for all $\N\in \mathbb{Z}^k$. Also if $F,G\in (S_\Lambda)_\N$ and $F\cup G\in S_\Lambda$, then clearly, $F\cup G\in (S_\Lambda)_\N$. Thus, by Proposition \ref{pro natural partial order and set-inclusion}, each $(S_\Lambda)_\N$ is join-closed with respect to $\subseteq$. Let $\N,\M\in \mathbb{Z}^k$, $F\in (S_\Lambda)_\N$ and $G\in (S_\Lambda)_\M$. If $FG=\emptyset$, then $FG\in (S_\Lambda)_{\N+\M}$. So assume that $FG\neq \emptyset$ and let $(\lambda,\mu)\in FG$. Then there are $(\gamma,\delta)\in F$ and $(\xi,\eta)\in G$ such that $(\lambda,\mu)\in (\gamma,\delta)\cdot (\xi,\eta)$. This implies that $\lambda=\gamma\alpha$ and $\mu=\eta\beta$ for some $(\alpha,\beta)\in \Lambda^{\min}(\delta,\xi)$. Since $\delta\alpha=\xi\beta$, $d(\delta)+d(\alpha)=d(\xi)+d(\beta)$. Consequently, we have \[d(\lambda)-d(\mu)=d(\gamma)+d(\alpha)-d(\eta)-d(\beta)=d(\gamma)+d(\xi)-d(\delta)-d(\eta)=\left(d(\gamma)-d(\delta)\right)+\left(d(\xi)-d(\eta)\right)=\N+\M.\] It follows that $FG\in (S_\Lambda)_{\N+\M}$. This shows that $(S_\Lambda)_\N (S_\Lambda)_\M\subseteq (S_\Lambda)_{\N+\M}$. Finally, for any $\emptyset\neq F\in S_\Lambda$ and $\N\in \mathbb{Z}^k$, letting \[F_\N=\{(\lambda,\mu)\in F~|~d(\lambda)-d(\mu)=\N\},\] we have $F_n\in (S_\Lambda)_\N$ and $F=\displaystyle{\bigvee_{\N\in \mathbb{Z}^k}}~F_\N$ by Proposition \ref{pro natural partial order and set-inclusion}. Clearly, this is the unique expression of $F$ as a join of elements from $(S_\Lambda)_\N$, $\N\in \mathbb{Z}^k$. Finally, if $F\in (S_\Lambda)_\N$ for some $\N\in \mathbb{Z}^k$, then it is clear that $d(\lambda)-d(\mu)=-\N$ for all $(\lambda,\mu)\in F^{-1}$. So $F^{-1}\in (S_\Lambda)_{-\N}$, which shows that $(S_\Lambda)_\N^{-1}\subseteq (S_\Lambda)_{-\N}$. Hence the lemma follows.
\end{proof}
Suppose $(S,\leqslant)$ is graded-ordered by an abelian group $G$. We define \[S^h:=\displaystyle{\bigcup_{x\in G}}~S_x.\] From Condition $(ii)$ of Definition \ref{def graded-ordered semigroup}, it follows that $S^h$ is a subsemigroup of $S$. Moreover, if we set $S^h_x:=S_x$ for all $x\in G$, then $S^h$ becomes a $G$-graded semigroup. We call $S^h$, the \emph{graded core} of $(S,\leqslant)$. If $S=S^h$, then clearly $S$ becomes $G$-graded. For any directed graph $\Gamma$, if we consider the graph inverse semigroup $S(\Gamma)$ to be graded-ordered by $\mathbb{Z}$ as in Examples \ref{ex graded-ordered semigroup}$(iv)$, then $S(\Gamma)=S(\Gamma)^h$.

\begin{rmk}\label{rem the graded core functor}
Suppose $G$ is an abelian group and $S,T$ are semigroups graded-ordered by $G$. Let us call a $0$-preserving homomorphism $f:S\longrightarrow T$ \emph{exact}, if it satisfies the following conditions:

$(i)$ $f(S_x)\subseteq T_x$ for all $x\in G$;

$(ii)$ $f(s_1\vee s_2)=f(s_1)\vee f(s_2)$ for all $s_1,s_2\in S$ which admit a join in $S$.

If $G$ is fixed, then one can check that the collection of all semigroups with $0$ graded-ordered by $G$ together with exact homomorphisms forms a category. We denote this category by $G$-\textbf{GrOS}. If $f:S\longrightarrow T$ is an exact homomorphism, then by condition $(i)$, $f|_{S^h}:S^h\longrightarrow T^h$ is a $0$-preserving graded homomorphism. Indeed, the graded core construction can be viewed as a functor $\mathcal{H}$ from $G$-\textbf{GrOS} to $G$-\textbf{GrS}, the category of $G$-graded semigroups with $0$ together with $0$-preserving graded homomorphisms. If $S$ is an inverse semigroup graded-ordered by $G$, then the condition $S_x^{-1}\subseteq S_{x^{-1}}$ for all $x\in G$, forces $S^h$ to be an inverse subsemigroup of $S$. Hence, if we denote by $G$-\textbf{GrOIS}, the full subcategory of $G$-\textbf{GrOS} consisting of all inverse semigroups with $0$ graded-ordered by $G$, then $\mathcal{H}$ restricts to a functor from $G$-\textbf{GrOIS} to $G$-\textbf{GrIS}, the full subcategory of $G$-\textbf{GrS} consisting of all graded inverse semigroups with $0$.\qed
\end{rmk}

We now proceed to relate the graded core $S_\Lambda^h$ of the inverse semigroup $S_\Lambda$ with the Kumjian--Pask algebra of $\Lambda$. In what follows, $\mathsf{k}$ denotes an arbitrary commutative unital ring.

If $S$ is an inverse semigroup with $0$, then recall that the \emph{tight ideal} of $\mathsf{k}_0S$ \cite{Steinberg} is defined by 
\begin{align*}
\mathcal{T}_\mathsf{k}(S)&=\left\langle \overline{e}-\left(\displaystyle{\bigvee_{f\in \mathcal{Z}}}~\overline{f}\right)~|~e\in \mathcal{E}(S)~\text{and}~\mathcal{Z}~\text{is a finite cover of}~e\right\rangle\\
&=\left\langle \displaystyle{\prod_{f\in \mathcal{Z}}}(\overline{e}-\overline{f})~|~e\in \mathcal{E}(S)~\text{and}~\mathcal{Z}~\text{is a finite cover of}~e\right\rangle.    
\end{align*}
The next proposition is in order.
\begin{prop}\label{pro the tight ideal of the graded core algebra}
Let $\Lambda$ be a finitely aligned $k$-graph. Then

$(i)$ $S_\Lambda^h$ is a full inverse subsemigroup of $S_\Lambda$,

$(ii)$ $\mathcal{T}_\mathsf{k}(S_\Lambda^h)$ is the ideal of $\mathsf{k}_0 S_\Lambda^h$ generated by the elements
\begin{itemize}
    \item[\textup{(I)}] $\overline{F\cup G}-(\overline{F}+\overline{G})$ for all $F,G\in \mathcal{E}(S_\Lambda^h)$ with $F\cup G\in S_\Lambda$ and $F\cap G=\emptyset$;

    \item[\textup{(II)}] $\displaystyle{\prod_{\tau\in X}}\left(\overline{\{(v,v)\}}-\overline{\{(\tau,\tau)\}}\right)$ for all $v\in \Lambda^0$ and finite exhaustive subset $X$ of $v\Lambda$.
\end{itemize}

$(iii)$ $\overline{F}-\displaystyle{\sum_{i=1}^{n}}~\overline{F_i}\in \mathcal{T}_\mathsf{k}(S_\Lambda^h)$ for all $F,F_i\in S_\Lambda^h$, $i=1,2,\ldots,n$ with $F=\displaystyle{\bigsqcup_{i=1}^{n}}~F_i$. 
\end{prop}
\begin{proof}
$(i)$ Clearly, $S_\Lambda^h$ is an inverse subsemigroup of $S_\Lambda$ (see Remark \ref{rem the graded core functor}). From Lemma \ref{lem idempotents of the semigroup}, it follows that $\mathcal{E}(S_\Lambda)\subseteq (S_\Lambda)_0\subseteq S_\Lambda^h$. Therefore, $S_\Lambda^h$ is full. 

$(ii)$ By part $(i)$ of this proposition, we have $\mathcal{E}(S_\Lambda^h)=\mathcal{E}(S_\Lambda)$. Also note that the set $\{\overline{F}~|~F\in \mathcal{E}(S_\Lambda^h)\}$ generates a generalized Boolean algebra inside $\mathsf{k}_0S_\Lambda^h$, where the join is given by $\overline{E}\vee \overline{F}=\overline{E}+\overline{F}-\overline{EF}$. In view of these facts, we now run a proof similar to that of Theorem \ref{th description of Cuntz-Krieger ideal}. Let $\mathcal{I}$ be the ideal of $\mathsf{k}_0S_\Lambda^h$ generated by elements in (I) and (II). From the definition of the tight ideal, it follows that the elements in (I) and (II) are inside $\mathcal{T}_\mathsf{k}(S_\Lambda^h)$. Hence, $\mathcal{I}\subseteq \mathcal{T}_\mathsf{k}(S_\Lambda^h)$. Now let $F\in \mathcal{E}(S_\Lambda^h)$ and $\mathcal{Z}$ any finite cover of $F$ in $S_\Lambda^h$. We have 
\begin{equation}\label{eq TI1}
\overline{F}-\displaystyle{\sum_{\lambda\in \mathsf{D}(F)}}\overline{\{(\lambda,\lambda)\}}\in \mathcal{I}.   
\end{equation}
In view of $(i)$, we can proceed as in the proof of Theorem \ref{th description of Cuntz-Krieger ideal} to obtain the finite exhaustive subset \[X_\lambda=\{\tau\in \Lambda~|~(\lambda\tau,\lambda\tau)\in G_{(\lambda,\lambda)}~\text{for some}~G\in \mathcal{Z}\}\] of $s(\lambda)\Lambda$ for any $\lambda\in \mathsf{D}(F)$. Then for any $\lambda\in \mathsf{D}(F)$, 
\begin{align*}
&\overline{\{(\lambda,\lambda)\}}-\displaystyle{\bigvee_{\tau\in X_\lambda}}\overline{\{(\lambda\tau,\lambda\tau)\}}\\
=&\overline{\{(\lambda,s(\lambda))\}\{(s(\lambda),s(\lambda))\}\{(s(\lambda),\lambda)\}}-\displaystyle{\bigvee_{\tau\in X_\lambda}}\left(\overline{\{(\lambda,s(\lambda))\}\{(\tau,\tau)\}\{(s(\lambda),\lambda)\}}\right)\\
=&\overline{\{(\lambda,s(\lambda))\}\{(s(\lambda),s(\lambda))\}\{(s(\lambda),\lambda)\}}-\overline{\{(\lambda,s(\lambda))\}}\left(\displaystyle{\bigvee_{\tau\in X_\lambda}}\overline{\{(\tau,\tau)\}}\right)\overline{\{(s(\lambda),\lambda)\}}\\
=&\overline{\{(\lambda,s(\lambda))\}}\left(\displaystyle{\prod_{\tau\in X_\lambda}}\left(\overline{\{(s(\lambda),s(\lambda)\}}-\overline{\{(\tau,\tau)\}}\right)\right)\overline{\{(s(\lambda),\lambda)\}}\in \mathcal{I}.
\end{align*}
The second equality above holds since for any two distinct $\tau,\tau'\in X_\lambda$, we have 
\[\{(\tau,\tau)\}\{(s(\lambda),\lambda)\}\{(\lambda,s(\lambda))\}\{(\tau',\tau')\}=\{(\tau,\tau)\}\{(\tau',\tau')\}.\]
Now it follows that 
\begin{equation}\label{eq:TI2}
\displaystyle{\sum_{\lambda\in \mathsf{D}(F)}}\left(\overline{\{(\lambda,\lambda)\}}-\displaystyle{\bigvee_{\tau\in X_\lambda}}\overline{\{(\lambda\tau,\lambda\tau)\}}\right)\in \mathcal{I}.
\end{equation}
Note that for distinct $\lambda,\mu\in \mathsf{D}(F)$, since $\Lambda^{\min}(\lambda,\mu)=\emptyset$, $\Lambda^{\min}(\lambda\tau,\mu\tau')=\emptyset$ for any $\tau\in X_\lambda$ and $\tau'\in X_\mu$. Thus $\displaystyle{\sum_{\lambda\in \mathsf{D}(F)}}\left(\displaystyle{\bigvee_{\tau\in X_\lambda}}\overline{\{(\lambda\tau,\lambda\tau)\}}\right)=\displaystyle{\bigvee_{\lambda\in \mathsf{D}(F)}}\left(\displaystyle{\bigvee_{\tau\in X_\lambda}}\overline{\{(\lambda\tau,\lambda\tau)\}}\right)$, whence 
\begin{align}
\displaystyle{\sum_{\lambda\in \mathsf{D}(F)}}\left(\displaystyle{\bigvee_{\tau\in X_\lambda}}\overline{\{(\lambda\tau,\lambda\tau)\}}\right)-\displaystyle{\bigvee_{G\in \mathcal{Z}}}~\overline{G}&=\displaystyle{\bigvee_{\lambda\in \mathsf{D}(F)}}\left(\displaystyle{\bigvee_{G\in \mathcal{Z}}}\overline{G_{(\lambda,\lambda)}}\right)-\displaystyle{\bigvee_{G\in \mathcal{Z}}}~\overline{G} \nonumber \\
&=\displaystyle{\bigvee_{G\in \mathcal{Z}}}\left(\displaystyle{\bigvee_{\lambda\in \mathsf{D}(F)}}\overline{G_{(\lambda,\lambda)}}\right)-\displaystyle{\bigvee_{G\in \mathcal{Z}}}~\overline{G} \nonumber\\
&=\displaystyle{\bigvee_{G\in \mathcal{Z}}}\left(\left(\displaystyle{\sum_{\lambda\in \mathsf{D}(F)}}\overline{G_{(\lambda,\lambda)}}-\overline{G}\right)+\overline{G}\right)-\displaystyle{\bigvee_{G\in \mathcal{Z}}}~\overline{G}\in \mathcal{I}
\label{eq TI3}
\end{align}
since $\displaystyle{\sum_{\lambda\in \mathsf{D}(F)}}\overline{G_{(\lambda,\lambda)}}-\overline{G}\in \mathcal{I}$. Now using (\ref{eq TI1})--(\ref{eq TI3}), we have $\overline{F}-\displaystyle{\bigvee_{G\in \mathcal{Z}}}~\overline{G}\in \mathcal{I}$. Hence $\mathcal{T}_\mathsf{k}(S_\Lambda^h)\subseteq \mathcal{I}$

$(iii)$ Suppose $F,F_i\in S_\Lambda^h$ for $i=1,2,\ldots,n$ and $F=\displaystyle{\bigsqcup_{i=1}^{n}}~F_i$. Then we have $F^{-1}F=\displaystyle{\bigsqcup_{i=1}^{n}}~(F_i^{-1}F_i)$. Now applying part $(ii)$ of this proposition, it follows that $\overline{F^{-1}F}-\displaystyle{\sum_{i=1}^{n}}~\overline{(F_i^{-1}F_i)}\in \mathcal{T}_\mathsf{k}(S_\Lambda^h)$. Using this together with the fact that $F_i\subseteq F$ and hence $F_i\le F$, we have \[\overline{F}-\displaystyle{\sum_{i=1}^{n}}~\overline{F_i}=\overline{F}\left(\overline{F^{-1}F}-\displaystyle{\sum_{i=1}^{n}}~\overline{(F_i^{-1}F_i)}\right)\in \mathcal{T}_\mathsf{k}(S_\Lambda^h).\qedhere\]
\end{proof}
Since the graded core $S_\Lambda^h$ is graded by $\mathbb{Z}^k$, we can view the contracted semigroup algebra $\mathsf{k}_0S_\Lambda^h$ as a graded $\mathsf{k}$-algebra with respect to the grading induced by the grading of $S_\Lambda^h$. Moreover, the tight ideal being generated by homogeneous elements ($\mathcal{E}(S_\Lambda^h)\subseteq (S_\Lambda)_0$ and so each element in (I) and (II) of Proposition \ref{pro the tight ideal of the graded core algebra} is homogeneous of degree $0$ in $\mathsf{k}_0S_\Lambda^h$), is graded. It follows that the grading of $S_\Lambda^h$ induces a $\mathbb{Z}^k$-grading on the quotient $\mathsf{k}_0S_\Lambda^h/{\mathcal{T}_\mathsf{k}(S_\Lambda^h)}$.  

The following theorem can be viewed as an algebraic analogue of \cite[Theorem 4.7]{Larki} in the setting of higher-rank graphs. This also extends the well-known realization of the Leavitt path algebra of a graph as the contracted semigroup algebra of the tight Booleanization of the corresponding graph inverse semigroup.
\begin{thm}\label{th the relationship between S_Lambda and KPA}
Let $\Lambda$ be a finitely aligned $k$-graph. Then there is a graded isomorphism of $\mathbb{Z}^k$-graded algebras \[\Theta:\mathsf{k}_0S_\Lambda^h/{\mathcal{T}_\mathsf{k}(S_\Lambda^h)}\longrightarrow \KP_\mathsf{k}(\Lambda)\] such that $(\Theta\circ \sigma)(F)=\displaystyle{\sum_{(\lambda,\mu)\in F}}s_\lambda s_{\mu^*}$ for all $F\in S_\Lambda^h$, where $\sigma:S_\Lambda^h\longrightarrow \mathsf{k}_0S_\Lambda^h/{\mathcal{T}_\mathsf{k}(S_\Lambda^h)}$ is the canonical universal homomorphism sending covers of idempotents to joins.
\end{thm}
\begin{proof}
Define a map $\theta:S_\Lambda^h\longrightarrow \KP_\mathsf{k}(\Lambda)$ by $\theta(F):=\displaystyle{\sum_{(\lambda,\mu)\in F}}s_\lambda s_{\mu^*}$ for all $F\in S_\Lambda^h$. Viewing $\KP_\mathsf{k}(\Lambda)$ as a semigroup under multiplication, we verify that $\theta$ is a semigroup homomorphism. If $\N\in \mathbb{Z}^k$ and $F_1,F_2\in (S_\Lambda)_\N$ are such that $F_1\cap F_2=\emptyset$ and $F_1\cup F_2\in S_\Lambda$, then note that
\begin{equation}\label{eq disjoint preservation}
\theta(F_1\cup F_2)=\theta(F_1)+\theta(F_2).
\end{equation}
in $\KP_\mathsf{k}(\Lambda)$. Let $F,G\in S_\Lambda^h$. If $(\lambda,\mu),(\lambda',\mu')\in F$ and $(\xi,\eta),(\xi',\eta')\in G$ are such that $\left((\lambda,\mu),(\xi,\eta)\right)\neq \left((\lambda',\mu'),(\xi',\eta')\right)$, then it follows that $\left((\lambda,\mu)\cdot (\xi,\eta)\right)\cap \left((\lambda',\mu')\cdot (\xi',\eta')\right)=\emptyset$. Consequently, the union in (\ref{eq the multiplication}) is a disjoint union. Now using (\ref{eq disjoint preservation}) and (KP3), we get
\begin{align*}
\theta(FG)=\theta\left(\displaystyle{\bigsqcup_{\substack{(\lambda,\mu)\in F\\(\xi,\eta)\in G}}}(\lambda,\mu)\cdot (\xi,\eta)\right)&=\displaystyle{\sum_{\substack{(\lambda,\mu)\in F\\(\xi,\eta)\in G}}}\theta\left((\lambda,\mu)\cdot (\xi,\eta)\right)\\
&=\displaystyle{\sum_{\substack{(\lambda,\mu)\in F\\(\xi,\eta)\in G}}}\left(\displaystyle{\sum_{(\alpha,\beta)\in \Lambda^{\min}(\mu,\xi)}}s_{\lambda\alpha}s_{(\eta\beta)^*}\right)\\
&=\displaystyle{\sum_{\substack{(\lambda,\mu)\in F\\(\xi,\eta)\in G}}} s_\lambda\left(\displaystyle{\sum_{(\alpha,\beta)\in \Lambda^{\min}(\mu,\xi)}}s_\alpha s_{\beta^*}\right)s_{\eta^*}\\
&=\displaystyle{\sum_{\substack{(\lambda,\mu)\in F\\(\xi,\eta)\in G}}}s_\lambda (s_{\mu^*}s_\xi)s_{\eta^*}=\left(\displaystyle{\sum_{(\lambda,\mu)\in F}}s_\lambda s_{\mu^*}\right)\left(\displaystyle{\sum_{(\xi,\eta)\in G}}s_\xi s_{\eta^*}\right)=\theta(F)\theta(G).
\end{align*}
Hence $\theta$ is a homomorphism as desired. By the universal property of the contracted semigroup algebra, there is a unique $\mathsf{k}$-algebra homomorphism $\tilde{\theta}:\mathsf{k}_0S_\Lambda^h\longrightarrow \KP_\mathsf{k}(\Lambda)$ that extends $\theta$. For any $\N\in \mathbb{Z}^k$ and $F\in (S_\Lambda)_\N$, it is clear that $\theta(F)\in \KP_\mathsf{k}(\Lambda)_\N$. Since $(\mathsf{k}_0S_\Lambda^h)_\N$ is the $\mathsf{k}$-subspace of $\mathsf{k}_0S_\Lambda^h$ generated by $(S_\Lambda)_\N$, it follows that $\tilde{\theta}$ is a graded homomorphism.

Next we claim that $\mathcal{T}_\mathsf{k}(S_\Lambda^h)\subseteq \Ker(\tilde{\theta})$. By Proposition \ref{pro the tight ideal of the graded core algebra}$(ii)$, it suffices to check that the elements in (I) and (II) of the said proposition are inside the kernel. That the elements in (I) are in the kernel readily follows from (\ref{eq disjoint preservation}). Let $v\in \Lambda^0$ and $X$ a finite exhaustive subset of $v\Lambda$. Then by (KP4), \[\tilde{\theta}\left(\displaystyle{\prod_{\tau\in X}}(\overline{\{(v,v)\}}-\overline{\{(\tau,\tau)\}}\right)=\displaystyle{\prod_{\tau\in X}}(s_v-s_\tau s_{\tau^*})=0.\] Hence the claim is established. As a consequence, $\tilde{\theta}$ induces a well-defined graded $\mathsf{k}$-algebra homomorphism \[\Theta:\mathsf{k}_0S_\Lambda^h/{\mathcal{T}_\mathsf{k}(S_\Lambda^h)}\longrightarrow \KP_\mathsf{k}(\Lambda).\] If $\pi:\mathsf{k}_0S_\Lambda^h\longrightarrow \mathsf{k}_0S_\Lambda^h/{\mathcal{T}_\mathsf{k}(S_\Lambda^h)}$ is the canonical quotient map, then it follows that $(\Theta\circ \sigma) (F)=\Theta(\pi(\overline{F}))=\tilde{\theta}(\overline{F})=\theta(F)$ for all $F\in S_\Lambda^h$.

To show that $\Theta$ is an isomorphism, we construct its inverse. For any $\lambda,\mu\in \Lambda$, we define $T_\lambda:=\sigma\left(\{(\lambda,s(\lambda))\}\right)$ and $T_{\mu^*}:=\sigma\left(\{(s(\mu),\mu)\}\right)$. We observe that $\{T_\lambda,T_\mu^*~|~\lambda,\mu\in \Lambda$ is a Kumjian--Pask $\Lambda$-family in the $\mathsf{k}$-algebra $\mathsf{k}_0S_\Lambda^h/{\mathcal{T}_\mathsf{k}(S_\Lambda^h)}$. It is easy to check that this family satisfies all the relations (KP1)--(KP4). However, for the reader's convenience, we work out the cases of (KP3) and (KP4). 

Let $\lambda,\mu\in \Lambda$. Then in $\mathsf{k}_0S_\Lambda^h/{\mathcal{T}_\mathsf{k}(S_\Lambda^h)}$, we have 
\begin{align*}
T_{\lambda^*}T_\mu&=\sigma\left(\{(s(\lambda),\lambda)\}\right)\sigma\left(\{(\mu,s(\mu))\}\right)\\
&=\sigma\left((s(\lambda),\lambda)\cdot (\mu,s(\mu))\right)\\
&=\sigma\left(\Lambda^{\min}(\lambda,\mu)\right)\\
&=\displaystyle{\sum_{(\alpha,\beta)\in \Lambda^{\min}(\lambda,\mu)}}\sigma\left(\{(\alpha,\beta)\}\right)~(\text{using Proposition \ref{pro the tight ideal of the graded core algebra}}(iii))\\
&=\displaystyle{\sum_{(\alpha,\beta)\in \Lambda^{\min}(\lambda,\mu)}}\left(\sigma\left(\{(\alpha,s(\alpha))\}\right)\sigma\left(\{(s(\beta),\beta)\}\right)\right)=\displaystyle{\sum_{(\alpha,\beta)\in \Lambda^{\min}(\lambda,\mu)}}T_\alpha T_{\beta^*}.
\end{align*}
Hence, the family $\{T_\lambda,T_{\mu^*}\}$ satisfies (KP3).

For (KP4), let $v\in \Lambda^0$ and $X$ a finite exhaustive subset of $v\Lambda$. Then 
\begin{align*}
\displaystyle{\prod_{\tau\in X}}(T_v-T_\tau T_{\tau^*})&=\displaystyle{\prod_{\tau\in X}}\left(\sigma\left(\{(v,v)\}\right)-\sigma\left(\{(\tau,s(\tau))\}\right)\sigma\left(\{(s(\tau),\tau)\}\right)\right)\\
&=\displaystyle{\prod_{\tau\in X}}\left(\sigma\left(\{(v,v)\}\right)-\sigma\left(\{(\tau,\tau)\}\right)\right)\\
&=\displaystyle{\prod_{\tau\in X}}\left(\pi\left(\overline{\{(v,v)\}}\right)-\pi\left(\overline{\{(\tau,\tau)\}}\right)\right)\\
&=\pi\left(\displaystyle{\prod_{\tau\in X}} \left(\overline{\{(v,v)\}}-\overline{\{(\tau,\tau)\}}\right)\right)=0~(\text{using Proposition \ref{pro the tight ideal of the graded core algebra}}(ii)).
\end{align*}
Thus (KP4) is also verified. By the universal property of $\KP_\mathsf{k}(\Lambda)$, there is a unique $\mathsf{k}$-algebra homomorphism \[\Phi:\KP_\mathsf{k}(\Lambda)\longrightarrow \mathsf{k}_0S_\Lambda^h/{\mathcal{T}_\mathsf{k}(S_\Lambda^h)}\] such that $\Phi(s_\lambda)=T_\lambda$ and $\Phi(s_{\mu^*})=T_{\mu^*}$ for all $\lambda,\mu\in \Lambda$. Since $\Phi$ preserves the degrees of the generators of $\KP_\mathsf{k}(\Lambda)$, it is a graded homomorphism. For any $\lambda\in \Lambda$, \[\Theta(\Phi(s_\lambda))=\Theta(T_\lambda)=(\Theta\circ \sigma)\left(\{(\lambda,s(\lambda))\}\right)=\theta\left(\{(\lambda,s(\lambda))\}\right)=s_\lambda.\] Similarly, $\Theta(\Phi(s_{\mu^*}))=s_{\mu^*}$. It follows that $\Theta\circ \Phi$ is the identity map on $\KP_\mathsf{k}(\Lambda)$. On the other hand, for any $F\in S_\Lambda^h$, we have $\overline{F}-\displaystyle{\sum_{(\lambda,\mu)\in F}}\overline{\{(\lambda,\mu)\}}\in \mathcal{T}_\mathsf{k}(S_\Lambda^h)$ by Proposition \ref{pro the tight ideal of the graded core algebra}$(iii)$. Hence \[\Phi(\Theta(\sigma(F)))=\Phi(\theta(F))=\Phi\left(\displaystyle{\sum_{(\lambda,\mu)\in F}}s_\lambda s_{\mu^*}\right)=\displaystyle{\sum_{(\lambda,\mu)\in F}}T_\lambda T_{\mu^*}
=\displaystyle{\sum_{(\lambda,\mu)\in F}}\sigma\left(\{(\lambda,\mu)\}\right)=\sigma(F).\]
Now since any element of $\mathsf{k}_0S_\Lambda^h/{\mathcal{T}_\mathsf{k}(S_\Lambda^h)}$ can be written as a $\mathsf{k}$-linear combination of elements of the form $\sigma(F)$, it follows that $\Phi\circ \Theta$ is the identity on $\mathsf{k}_0S_\Lambda^h/{\mathcal{T}_\mathsf{k}(S_\Lambda^h)}$. This completes the proof.
\end{proof}
In the next section, after defining suitable categories and functors, we will show that the isomorphism $\Theta$ of Theorem \ref{th the relationship between S_Lambda and KPA} is natural.

For now, as an immediate corollary to Theorem \ref{th the relationship between S_Lambda and KPA}, we record the following result.
\begin{cor}\label{cor relation of KP with SA of tight groupoid of FMY}
Let $\Lambda$ be a finitely aligned $k$-graph and $S_\Lambda$ the Farthing--Muhly--Yeend inverse semigroup of $\Lambda$. Then $\KP_\mathsf{k}(\Lambda)\cong^{\gr} A_\mathsf{k}(\mathcal{G}_{tight}(S_\Lambda))$.
\end{cor}
\begin{proof}
Let $S$ be any $G$-graded inverse semigroup with $0$. Then $\mathsf{k}_0S/{\mathcal{T}_\mathsf{k}(S)}\cong A_\mathsf{k}(\mathcal{G}_{tight}(S))$ by \cite[Corollary 2.14]{Steinberg}. Moreover, this isomorphism maps the canonical image of $s\in S$ in $\mathsf{k}_0S/{\mathcal{T}_\mathsf{k}(S)}$ to the characteristic function $1_{\left(s,D(s^{-1}s)~\cap~\widehat{\mathcal{E}(S)}_{tight}\right)}$, where $D(s^{-1}s)=\{\chi\in \widehat{\mathcal{E}(S)}~|~\chi(s^{-1}s)=1\}$ and $\left(s,D(s^{-1}s)\cap \widehat{\mathcal{E}(S)}_{tight}\right)$ is the bisection consisting of the germs of $s$ at the characters in $D(s^{-1}s)\cap \widehat{\mathcal{E}(S)}_{tight}$. In view of the natural grading of $\mathcal{G}_{tight}(S)$ induced by the grading of $S$ (see \cite[\S 6.2]{HazMes}), we have $\left(s,D(s^{-1}s)\cap \widehat{\mathcal{E}(S)}_{tight}\right)\subseteq \mathcal{G}_{tight}(S)_{\deg(s)}$. It follows that the isomorphism between $\mathsf{k}_0S/{\mathcal{T}_\mathsf{k}(S)}$ and $A_\mathsf{k}(\mathcal{G}_{tight}(S))$ is graded. This, together with Theorem \ref{th the relationship between S_Lambda and KPA}, completes the proof.
\end{proof}

\section{Extending \texorpdfstring{$k$}{k}-graph morphisms to semigroup homomorphisms}\label{sec extending morphisms}
Suppose $\Lambda_1$ and $\Lambda_2$ are finitely aligned $k$-graphs and $p:\Lambda_1\longrightarrow \Lambda_2$ is a $k$-graph morphism. For any $F\in S_{\Lambda_1}$, we write \[p(F):=\{(p(\lambda),p(\mu))~|~(\lambda,\mu)\in F\}.\] It is natural to ask when $p(F)\in S_{\Lambda_2}$ for any $F\in S_{\Lambda_1}$. The following lemma provides the answer.

\begin{lem}\label{lem induced elements}
Let $(\Lambda_1,d_1),(\Lambda_2,d_2)$ be finitely aligned $k$-graphs and $p:\Lambda_1\longrightarrow \Lambda_2$ a $k$-graph morphism. Then the following are equivalent.

$(i)$ $p$ is injective and $p(F)\in S_{\Lambda_2}$ for all $F\in S_{\Lambda_1}$.

$(ii)$ $p$ preserves orthogonality in the sense that $\Lambda_1^{\min}(\lambda,\mu)=\emptyset$ implies $\Lambda_2^{\min}(p(\lambda),p(\mu))=\emptyset$.
\end{lem}
\begin{proof}
$(i)\Rightarrow (ii)$ Let $\lambda,\mu\in \Lambda_1$ be such that $\Lambda_1^{\min}(\lambda,\mu)=\emptyset$. Then $F:=\{(\lambda,\lambda),(\mu,\mu)\}\in S_{\Lambda_1}$ and hence $p(F)\in S_{\Lambda_2}$. Again since $p$ is injective, $p(\lambda)\neq p(\mu)$ and so $(p(\lambda),p(\lambda))$, $(p(\mu),p(\mu))$ are distinct pairs in $p(F)$. It now follows that $\Lambda_2^{\min}(p(\lambda),p(\mu))=\emptyset$. 

$(ii)\Rightarrow (i)$ We first show that $p$ is injective. Let $\lambda,\mu\in \Lambda_1$ be such that $\lambda\neq \mu$. If $d_1(\lambda)\neq d_1(\mu)$, then $d_2(p(\lambda))\neq d_2(p(\mu))$ since $p$ is a $k$-graph morphism. It follows that $p(\lambda)\neq p(\mu)$. If $d_1(\lambda)=d_1(\mu)$, then $\Lambda_1^{\min}(\lambda,\mu)=\emptyset$. Since $p$ preserves orthogonality, $\Lambda_2^{\min}(p(\lambda),p(\mu))=\emptyset$, whence $p(\lambda)\neq p(\mu)$. Hence $p$ is injective. Now suppose $(\lambda,\mu),(\xi,\eta)\in F$ are such that $(p(\lambda),p(\mu))\neq (p(\xi),p(\eta))$. Then $(\lambda,\mu)$ and $(\xi,\eta)$ are distinct pairs in $F$. So $(\lambda,\mu)\perp (\xi,\eta)$, which implies $(p(\lambda),p(\mu))\perp (p(\xi),p(\eta))$ since $p$ preserves orthogonality. Hence $p(F)\in S_{\Lambda_2}$. 
\end{proof} 
Assume that $p:\Lambda_1\longrightarrow \Lambda_2$ preserves orthogonality. In view of Lemma \ref{lem induced elements}, $p$ induces a well-defined map 
\begin{align*}
\tilde{p}:S_{\Lambda_1} &\longrightarrow S_{\Lambda_2}\\
F &\longmapsto p(F).
\end{align*}
In this section, we seek an answer to the following question.
\begin{ques}\label{query k-graph morphism inducing semigroup morphism}
Given finitely aligned $k$-graphs $\Lambda_1$, $\Lambda_2$ and a $k$-graph morphism $p:\Lambda_1\longrightarrow \Lambda_2$ preserving orthogonality, when is the induced map $\tilde{p}:S_{\Lambda_1} \longrightarrow S_{\Lambda_2}$ a $0$-preserving semigroup homomorphism?     
\end{ques}

When $\Lambda_1,\Lambda_2$ are $1$-graphs, then a $1$-graph morphism $p:\Lambda_1\longrightarrow \Lambda_2$ is precisely a graph homomorphism between the underlying directed graphs. In this case, $p$ preserves orthogonality if and only if $p$ is injective. However, for $k\ge 2$, an injective $k$-graph morphism may fail to preserve orthogonality. For example, consider the obvious $2$-graph morphism from the $2$-graph after Proposition \ref{pro u-closed ideals are join closed} to the $2$-graph of Remark \ref{rem recovering D of graph inverse semigroups}$(ii)$.

For directed graphs and graph inverse semigroups, the analogous version of Question \ref{query k-graph morphism inducing semigroup morphism} has been answered completely in \cite[Theorem 20]{Mesyan}: any injective graph homomorphism (equivalently, orthogonality preserving $1$-graph morphism) naturally induces a homomorphism between the corresponding graph inverse semigroups. However, as we find in the following example, the map induced by an orthogonality preserving $k$-graph morphism between two general $k$-graphs may not be a homomorphism between their inverse semigroups.

\begin{example}\label{ex orthogonal map does not produce semigroup hom}
Let $\Lambda_1$ and $\Lambda_2$ be $2$-graphs with $1$-skeletons shown below (the left one is of $\Lambda_1$ and the right one is of $\Lambda_2$):
\[
\begin{tikzpicture}[
    >=stealth,
    line width=1.2pt,
    every node/.style={font=\normalsize}
]


\node[circle, fill=black, inner sep=1.5pt, label=below:$a$] (a) at (0,0) {};
\node[circle, fill=black, inner sep=1.5pt, label=below:$b$] (b) at (2,0) {};

\draw[blue, ->] (b) -- node[above] {\textcolor{black}{$f$}} (a);

\path[] (a) edge[red,->, dashed, out=135, in=225, looseness=8, distance=1cm] node[left, black] {$e_1$} (a);

\path[] (b) edge[red,->, dashed, out=315, in=45, looseness=8, distance=1cm] node[right, black] {$e_2$} (b);


\node[circle, fill=black, inner sep=1.5pt, label=below:$u$] (u) at (5,0) {};
\node[circle, fill=black, inner sep=1.5pt, label=above:$v$] (v) at (7,2) {};
\node[circle, fill=black, inner sep=1.5pt, label=below:$z$] (z) at (7,0) {};
\node[circle, fill=black, inner sep=1.5pt, label=above:$w$] (w) at (5,2) {};

\draw[blue, ->]
    (w) -- node[left] {\textcolor{black}{$\alpha_1$}} (u);

\draw[blue, ->]
    (z) -- node[below] {\textcolor{black}{$\alpha_2$}} (u);

\draw[blue, ->]
    (v) -- node[above left] {\textcolor{black}{$\mu$}} (u);

\draw[red, dashed, ->]
    (w) -- node[above] {\textcolor{black}{$\beta_1$}} (v);

\draw[red, dashed, ->]
    (z) -- node[right] {\textcolor{black}{$\beta_2$}} (v);

\path[] (u) edge[red,->, dashed, out=135, in=225, looseness=8, distance=1cm] node[left, black] {$\lambda$} (u);

\path[] (v) edge[red,->, dashed, out=315, in=45, looseness=8, distance=1cm] node[right, black] {$\gamma$} (v);

\end{tikzpicture}
\]

In each of these $2$-graphs, the factorization rules are uniquely determined as evident from the skeletons. Define a graph homomorphism $p$ from the $1$-skeleton of $\Lambda_1$ to that of $\Lambda_2$ by $p(a):=u$, $p(b):=v$, $p(e_1):=\lambda$, $p(f):=\mu$ and $p(e_2):=\gamma$. Since $p$ is color preserving, it defines a unique $2$-graph homomorphism $p$ from $\Lambda_1$ to $\Lambda_2$. 

Suppose $x,y\in \Lambda_1$ are such that $r(x)=r(y)=a$. If $x=e_1^i$, $y=e_1^{i'}$ for some $i,i'\in \mathbb{N}$, then clearly $\Lambda_1^{\min}(x,y)\neq \emptyset$. If $x=e_1^i$, $y=e_1^{i'}fe_2^j$ for some $i,i',j\in \mathbb{N}$, then \[x e_1^{(i\vee i')-i}fe_2^j=ye_2^{(i\vee i')-i'}\] and so $\Lambda_1^{\min}(x,y)\neq \emptyset$. Lastly, if $x=e_1^ife_2^j$ and $y=e_1^{i'}fe_2^{j'}$ for some $i,i',j,j'\in \mathbb{N}$, then \[xe_2^{(i\vee i')-i+j'}=ye_2^{(i\vee i')-i'+j}\] and hence $\Lambda_1^{\min}(x,y)\neq \emptyset$. If $r(x)=r(y)=b$, then they are obviously of the form $x=e_2^j$, $y=e_2^{j'}$ for some $j,j'\in \mathbb{N}$, which implies $\Lambda_1^{\min}(x,y)\neq \emptyset$. Therefore if $\Lambda_1^{\min}(x,y)=\emptyset$, $x,y$ must have different ranges. But then $r(p(x))=p(r(x))\neq p(r(y))=r(p(y))$, whence $\Lambda_2^{\min}(p(x),p(y))=\emptyset$. Thus it follows that $p$ preserves orthogonality. 

Let $F=\{(e_1,e_1)\}$, $G=\{(f,e_2)\}$. Then $p(F)=\{(\lambda,\lambda)\}$ and $p(G)=\{(\mu,\gamma)\}$. Now $FG=\{(e_1f,e_2^2)\}$ and so $\tilde{p}(FG)=\{(\lambda\mu,\gamma^2)\}$, whereas \[p(F)p(G)=(\lambda,\lambda)\cdot (\mu,\gamma)=\{(\lambda\alpha_1,\gamma\beta_1),(\lambda\alpha_2,\gamma\beta_2),(\lambda\mu,\gamma^2)\}.\] Hence $\tilde{p}(FG)\neq \tilde{p}(F)\tilde{p}(G)$, which means that $\tilde{p}$ is not a semigroup homomorphism. 
\end{example}

The above example suggests that we need something more apart from the preservation of orthogonality in order to have a naturally induced semigroup homomorphism. The extra bit which is required is the following special property of a $k$-graph morphism.
\begin{dfn}\label{def the range path lifting}
Let $\Lambda_1,\Lambda_2$ be finitely aligned $k$-graphs and $p:\Lambda_1\longrightarrow \Lambda_2$ a $k$-graph morphism. We say that $p$ has \emph{unique $r$-path lifting} (resp., \emph{unique $s$-path lifting}) property if for any $v\in \Lambda_1^0$, $p|_{v\Lambda_1}:v\Lambda_1\longrightarrow p(v)\Lambda_2$ (resp., $p|_{\Lambda_1v}:\Lambda_1 v\longrightarrow \Lambda_2p(v)$) is bijective.
\end{dfn}
\begin{rmk}\label{rem path lifting and covering}
If a $k$-graph morphism $p:\Omega\longrightarrow \Lambda$ has both the unique $r$- and $s$-path lifting properties, then the pair $(\Omega,p)$ is called a covering of $\Lambda$ (see \cite[Definition 2.1]{Pask}). For directed graphs $\Gamma_1,\Gamma_2$, a graph homomorphism $p:\Gamma_1\longrightarrow \Gamma_2$ extends to a $1$-graph morphism with unique $r$-path lifting between the respective path categories if and only if $p$ is a directed cover in the sense of Meakin and Wang \cite{Meakin-Wang}.\qed    
\end{rmk}
We have seen in Lemma \ref{lem induced elements} that any $k$-graph morphism $p$ which preserves orthogonality is injective. Under the assumption that $p$ has unique $r$-path lifting property, the converse is also true.
\begin{lem}\label{lem injective+path lifting imply orthogonality preservation}
Suppose $p:\Lambda_1\longrightarrow \Lambda_2$ is a $k$-graph morphism which has the unique $r$-path lifting property. Then $p$ is injective if and only if $p$ preserves orthogonality.     
\end{lem}
\begin{proof}
Suppose $p$ is injective. Let $\lambda,\mu\in \Lambda_1$ be such that $\Lambda_2^{\min}(p(\lambda),p(\mu))\neq \emptyset$. Choose any $(\alpha,\beta)\in \Lambda_2^{\min}(p(\lambda),p(\mu))$. Since $p(s(\lambda))=s(p(\lambda))=r(\alpha)$, by unique $r$-path lifting, there exists a unique $\eta\in s(\lambda)\Lambda_1$ such that $p(\eta)=\alpha$. By similar argument, there exists a unique $\delta\in s(\mu)\Lambda_1$ with $p(\delta)=\beta$. Now \[p(\lambda\eta)=p(\lambda)\alpha=p(\mu)\beta=p(\mu\delta),\] which implies $\lambda\eta=\mu\delta$ by the injectivity of $p$. This implies $\Lambda_1^{\min}(\lambda,\mu)\neq \emptyset$ and so $p$ preserves orthogonality. The converse is already proved in Lemma \ref{lem induced elements}.
\end{proof}

We now state and prove the main result of this section, which gives an answer to Question \ref{query k-graph morphism inducing semigroup morphism}.
\begin{thm}\label{th answer to the query}
Let $(\Lambda_1,d_1)$, $(\Lambda_2,d_2)$ be finitely aligned $k$-graphs and $p:\Lambda_1\longrightarrow \Lambda_2$ a $k$-graph morphism such that 

$(i)$ $p$ is injective and 

$(ii)$ $p$ has the unique $r$-path lifting property.

Then the induced map $\tilde{p}:S_{\Lambda_1}\longrightarrow S_{\Lambda_2}$ is an injective $0$-preserving semigroup homomorphism.
\end{thm}
\begin{proof}
By Lemma \ref{lem induced elements} and Lemma \ref{lem injective+path lifting imply orthogonality preservation}, $p(F)\in S_{\Lambda_2}$ for all $F\in S_{\Lambda_1}$. Hence the map $\tilde{p}$ is well defined. Since $p(\emptyset)=\emptyset$, $\tilde{p}$ is $0$-preserving. We now show that $\tilde{p}(FG)=\tilde{p}(F)\tilde{p}(G)$ for all $F,G\in S_{\Lambda_1}$ by dealing with two cases.

\emph{Case 1}: Suppose $FG=\emptyset$. If $F=\emptyset$ or $G=\emptyset$, then obviously $p(F)p(G)=\emptyset$ and hence $\tilde{p}(FG)=\tilde{p}(F)\tilde{p}(G)$. So assume that $F,G\neq \emptyset$. By (\ref{eq the multiplication}), it follows that $\Lambda_1^{\min}(\mu,\xi)=\emptyset$ for all $\mu\in \mathsf{R}(F)$ and $\xi\in \mathsf{D}(G)$. Now $p$ preserves orthogonality by Lemma \ref{lem injective+path lifting imply orthogonality preservation}, and so $\Lambda_2^{\min}(p(\mu),p(\xi))=\emptyset$. Since $\mathsf{R}(p(F))=p(\mathsf{R}(F))$ and $\mathsf{D}(p(G))=p(\mathsf{D}(G))$, this implies $p(F)p(G)=\emptyset$. Therefore we are done in this case.

\emph{Case 2}: Suppose $FG\neq \emptyset$. Let $(\lambda,\mu)\in F$ and $(\xi,\eta)\in G$. We claim that $p\left((\lambda,\mu)\cdot (\xi,\eta)\right)\subseteq p(F)p(G)$. For this, we take any $(\alpha,\beta)\in \Lambda_1^{\min}(\mu,\xi)$. Note that \[p(\mu)p(\alpha)=p(\mu\alpha)=p(\xi\beta)=p(\xi)p(\beta)\] and \[d_2(p(\mu)p(\alpha))=d_2(p(\mu\alpha))=d_1(\mu\alpha)=d_1(\mu)\vee d_1(\xi)=d_2(p(\mu))\vee d_2(p(\xi)),\] which imply $(p(\alpha),p(\beta))\in \Lambda_2^{\min}(p(\mu),p(\xi))$. Hence \[(p(\lambda\alpha),p(\eta\beta))=(p(\lambda)p(\alpha),p(\eta)p(\beta))\in (p(\lambda),p(\mu))\cdot (p(\xi),p(\eta))\subseteq p(F)p(G)\] and the claim is established. Now \[\tilde{p}(FG)=\displaystyle{\bigcup_{\substack{(\lambda,\mu)\in F\\(\xi,\eta)\in G}}}p\left((\lambda,\mu)\cdot (\xi,\eta)\right)\subseteq p(F)p(G)=\tilde{p}(F)\tilde{p}(G).\] We now prove the reverse inclusion. Let $(x,y)\in \tilde{p}(F)\tilde{p}(G)$. Then there are $(\gamma,\nu)\in F$, $(\tau,\delta)\in G$ and $(\alpha',\beta')\in \Lambda_2^{\min}(p(\nu),p(\tau))$ such that $x=p(\gamma)\alpha'$ and $y=p(\delta)\beta'$. Since $\alpha'\in p(s(\nu))\Lambda_2$, by unique $r$-path lifting of $p$, there exists unique $\alpha\in s(\nu)\Lambda_1$ such that $p(\alpha)=\alpha'$. Similarly, there exists unique $\beta\in s(\tau)\Lambda_1$ such that $p(\beta)=\beta'$. Now since $p(\nu\alpha)=p(\nu)\alpha'=p(\tau)\beta'=p(\tau\beta)$ and $p$ is injective, we have $\nu\alpha=\tau\beta$. Again since \[d_1(\nu\alpha)=d_2(p(\nu\alpha))=d_2(p(\nu))\vee d_2(p(\tau))=d_1(\nu)\vee d_1(\tau),\] $(\alpha,\beta)\in \Lambda_1^{\min}(\nu,\tau)$. This implies that \[(\gamma\alpha,\delta\beta)\in (\gamma,\nu)\cdot (\tau,\delta)\subseteq FG\] and so \[(x,y)=(p(\gamma)\alpha',p(\delta)\beta')=(p(\gamma\alpha),p(\delta\beta))\in p\left((\gamma,\nu)\cdot (\tau,\delta)\right)\subseteq p(FG)=\tilde{p}(FG).\] Hence $\tilde{p}(F)\tilde{p}(G)\subseteq \tilde{p}(FG)$. 

We have proved that $\tilde{p}$ is a semigroup homomorphism. To show that it is injective, let $F,G\in S_{\Lambda_1}$ be such that $\tilde{p}(F)=\tilde{p}(G)$. If $F=\emptyset$, then $p(G)=p(F)=\emptyset$ and so $G=\emptyset$. Similarly if $G=\emptyset$, then $F=\emptyset$. Assume that $F,G\neq \emptyset$. Let $(\lambda,\mu)\in F$. Then $(p(\lambda),p(\mu))\in p(F)=p(G)$. So there exists $(\xi,\eta)\in G$ such that $p(\lambda)=p(\xi)$ and $p(\mu)=p(\eta)$. Since $p$ is injective, $\lambda=\xi$ and $\mu=\eta$. Hence $(\lambda,\mu)\in G$. This shows that $F\subseteq G$. Reversing the role of $F$ and $G$, we have $G\subseteq F$. Therefore $F=G$ and consequently, $\tilde{p}$ is injective. 
\end{proof}

The unique $r$-path lifting criterion is not necessary for the induced map $\tilde{p}$ to be a semigroup homomorphism. To see this, let $\M,\N\in \mathbb{N}^k$ be $k$-tuples such that $\M <\N$. Then we have a natural inclusion map $\iota:\Omega_{k,\M}\hookrightarrow \Omega_{k,\N}$ defined by $\iota((\bP,\Q)):=(\bP,\Q)$. For $\lambda,\mu\in \Omega_{k,\M}$, $\Omega_{k,\M}^{\min}(\lambda,\mu)=\emptyset$ if and only if $r(\lambda)\neq r(\mu)$, from which it follows that $\iota$ preserves orthogonality and hence $\tilde{\iota}$ is well defined. Note that these $k$-graphs are singly aligned \cite[Definition 20.1]{Exel}, and for $\lambda=(\bP,\Q_1)$, $\mu=(\bP,\Q_2)$ in $\Omega_{k,\M}$, we have \[\Omega_{k,\N}^{\min}(\iota(\lambda),\iota(\mu))=\left\{((\Q_1,\Q_1\vee \Q_2),(\Q_2,\Q_1\vee \Q_2))\right\}=\iota\left(\Omega_{k,\M}^{\min}(\lambda,\mu)\right).\] Using this, it can be observed that $\tilde{\iota}$ is a semigroup homomorphism. However, $\iota$ does not have the unique $r$-path lifting since $r((\M,\N))=(\M,\M)=\iota((\M,\M))$ but there is no path $\lambda\in (\M,\M)\Omega_{k,\M}$ with $\iota(\lambda)=(\M,\N)$. 

We now exhibit a specific feature of the semigroup homomorphisms induced by $k$-graph morphisms with unique $r$-path lifting property, namely they preserve the local submonoids of the inverse semigroups. In fact, this feature characterizes the unique path lifting property of a $k$-graph morphism inducing a semigroup homomorphism. 

For a finitely aligned $k$-graph $\Lambda$ and $v\in \Lambda^0$, we write \[vS_\Lambda v:=\{(v,v)\}S_\Lambda \{(v,v)\}=\{F\in S_\Lambda~|~r(\mathsf{D}(F))=r(\mathsf{R}(F))=\{v\}\}.\] Clearly, $vS_\Lambda v$ is an inverse submonoid of $S_\Lambda$ with identity $\{(v,v)\}$. 

\begin{thm}\label{th local submonoid preservation}
Let $\Lambda_1$, $\Lambda_2$ be finitely aligned $k$-graphs and $p:\Lambda_1\longrightarrow \Lambda_2$ a $k$-graph morphism preserving orthogonality. Then the following are equivalent.

$(i)$ $p$ has the unique $r$-path lifting property.

$(ii)$ $\tilde{p}$ is an injective $0$-preserving semigroup homomorphism such that for each $v\in \Lambda_1^0$, $\tilde{p}$ induces a monoid isomorphism from $vS_{\Lambda_1}v$ onto $p(v)S_{\Lambda_2}p(v)$. 
\end{thm}
\begin{proof}
$(i)\Rightarrow (ii)$ That $\tilde{p}$ is an injective $0$-preserving homomorphism is exactly proved in Theorem \ref{th answer to the query}. Let $v\in \Lambda_1^0$ and $F\in vS_{\Lambda_1}v$. Since we already have that $\tilde{p}$ is injective and $\tilde{p}(\{(v,v)\})=\{(p(v),p(v))\}$, it only remains to show that $\tilde{p}(vS_{\Lambda_1}v)=p(v)S_{\Lambda_2}p(v)$. Since \[r(\mathsf{D}(p(F)))=r(p(\mathsf{D}(F)))=p(r(\mathsf{D}(F)))=\{p(v)\}\] and similarly $r(\mathsf{R}(p(F)))=\{p(v)\}$, it follows that $\tilde{p}(vS_{\Lambda_1}v)\subseteq p(v)S_{\Lambda_2}p(v)$. Now let $G\in p(v)S_{\Lambda_2}p(v)$ and $(\mu,\nu)\in G$. Then $r(\mu)=r(\nu)=p(v)$. By the unique $r$-path lifting of $p$, there are unique $\lambda_\mu,\lambda_\nu\in v\Lambda_1$ such that $p(\lambda_\mu)=\mu$ and $p(\lambda_\nu)=\nu$. Define \[F:=\{(\lambda_\mu,\lambda_\nu)~|~(\mu,\nu)\in G\}.\] Suppose $(\lambda_\mu,\lambda_\nu)$ and $(\lambda_{\mu'},\lambda_{\nu'})$ are distinct pairs in $F$. Then $(\mu,\nu)$ and $(\mu',\nu')$ are distinct in $G$. Now if $\Lambda_1^{\min}(\lambda_\mu,\lambda_{\mu'})\neq \emptyset$, then $\lambda_\mu \alpha=\lambda_{\mu'}\beta$ for some $\alpha,\beta\in \Lambda_1$. Upon applying $p$, this gives $\mu p(\alpha)=\mu' p(\beta)$. But this is a contradiction since $\Lambda_2^{\min}(\mu,\mu')=\emptyset$. Therefore, $\Lambda_1^{\min}(\lambda_\mu,\lambda_{\mu'})=\emptyset$. Similarly, $\Lambda_1^{\min}(\lambda_\nu,\lambda_{\nu'})=\emptyset$. This implies that $(\lambda_\mu,\lambda_\nu)\perp (\lambda_{\mu'},\lambda_{\nu'})$ and so $F\in vS_{\Lambda_1}v$. It is clear that $p(F)=G$. Thus we have the desired equality and $(ii)$ follows.

$(ii)\Rightarrow (i)$ Let $v\in \Lambda^0$ and $\gamma\in p(v)\Lambda_2$. Now $\{(\gamma,\gamma)\}\in p(v)S_{\Lambda_2}p(v)$ and so by $(ii)$, there is a unique $F\in vS_{\Lambda_1}v$ such that $\tilde{p}(F)=p(F)=\{(\gamma,\gamma)\}$. Since $p$ is injective, $F$ must be of the form $\{(\lambda,\lambda)\}$ for some $\lambda\in \Lambda_1$. Then clearly $\lambda$ is the unique path in $v\Lambda_1$ such that $p(\lambda)=\gamma$. Hence $(i)$ follows. 
\end{proof}

Theorem \ref{th answer to the query} hints that the correspondence $\Lambda\longmapsto S_\Lambda$ is functorial. We now make this precise. For any positive integer $k$, let us denote by $k$-\textbf{Gra}, the category of all finitely aligned $k$-graphs with $k$-graph morphisms satisfying Conditions $(i)$ and $(ii)$ of Theorem \ref{th answer to the query}. On the other hand, we consider the category $\mathbb{Z}^k$-\textbf{GrOIS} of inverse semigroups with $0$ graded-ordered by $\mathbb{Z}^k$ together with exact homomorphisms (see Remark \ref{rem the graded core functor}).

\begin{prop}\label{pro Inverse semigroup construction as a functor}
Consider the correspondence $\mathcal{S}:k\textup{-\textbf{Gra}}\longrightarrow \mathbb{Z}^k\textup{-\textbf{GrOIS}}$ such that $\mathcal{S}(\Lambda)=S_\Lambda$ for all $\Lambda\in \Obj(k\textup{-\textbf{Gra}})$ and $\mathcal{S}(p)=\tilde{p}$ for all $p\in \Hom_{k\textup{-\textbf{Gra}}}(\Lambda_1,\Lambda_2)$. Then $\mathcal{S}$ is a faithful conservative functor.
\end{prop}
\begin{proof}
Let $F\in (S_{\Lambda_1})_\N$. For any $(x,y)\in p(F)$, we have $x=p(\lambda)$ and $y=p(\mu)$ for some $(\lambda,\mu)\in F$. Then $d_2(x)-d_2(y)=d_1(\lambda)-d_1(\mu)=\N$. It follows that $\tilde{p}(F)\in (S_{\Lambda_2})_\N$. Suppose $F,G\in S_{\Lambda_1}$ admit a join with respect to $\subseteq$. Then by Proposition \ref{pro natural partial order and set-inclusion}, $F\cup G$ is in $S_{\Lambda_1}$ and is the join of $F$ and $G$. Clearly, $p(F\cup G)=p(F)\cup p(G)$. These imply that $\tilde{p}$ is an exact homomorphism and hence the correspondence $\mathcal{S}$ is meaningful. Obviously $\mathcal{S}$ preserves identity morphisms and compositions of $k$-graph morphisms. 

To show that $\mathcal{S}$ is faithful, let $p,q\in \Hom_{k\textup{-\textbf{Gra}}}(\Lambda_1,\Lambda_2)$ be such that $\tilde{p}=\tilde{q}$. For any $\lambda\in \Lambda_1$, $\{(\lambda,\lambda)\}\in S_{\Lambda_1}$. By the assumption, $\{(p(\lambda),p(\lambda))\}=\{(q(\lambda),q(\lambda))\}$ in $S_{\Lambda_2}$. Hence, $p(\lambda)=q(\lambda)$. Since $\lambda$ is arbitrary, this shows that $p=q$. 

Let $p\in \Hom_{k\textup{-\textbf{Gra}}}(\Lambda_1,\Lambda_2)$ be such that $\tilde{p}$ is an isomorphism between $S_{\Lambda_1}$ and $S_{\Lambda_2}$. Then for any $\xi\in \Lambda_2$, there is an $F\in S_{\Lambda_1}$ such that $\tilde{p}(F)=\{(\xi,\xi)\}$. Consequently, $\xi=p(\lambda)$ for some $\lambda\in \Lambda_1$ showing that $p$ is surjective. Being a morphism in $k\textup{-\textbf{Gra}}$, it is already injective. Hence $p$ is a $k$-graph isomorphism. Therefore $\mathcal{S}$ reflects isomorphism and so is conservative. 
\end{proof}
\begin{rmk}\label{rem the functor S is not full}
Let $\M<\N$ and consider the injective semigroup homomorphism $\tilde{\iota}:S_{\Omega_{k,\M}}\longrightarrow S_{\Omega_{k,\N}}$ discussed after the proof of Theorem \ref{th answer to the query}. Using Theorem \ref{th local submonoid preservation}, it is not hard to see that there is no morphism $p\in \Hom_{k\textup{-\textbf{Gra}}}(\Omega_{k,\M},\Omega_{k,\N})$ such that $\mathcal{S}(p)=\tilde{\iota}$. Hence the functor $\mathcal{S}$ is not full. \qed
\end{rmk}
By considering the functor $\mathcal{H}$ of Remark \ref{rem the graded core functor}, we get the functor $\mathcal{H}\circ \mathcal{S}:k\textup{-\textbf{Gra}}\longrightarrow \mathbb{Z}^k\textup{-\textbf{GrIS}}$. To locate the image of this functor, we introduce the following notion. 
\begin{dfn}\label{def cover preserving}
Let $S_1,S_2$ be inverse semigroups with zeros and $\phi:S_1\longrightarrow S_2$ a homomorphism. Then $\phi$ is called \emph{cover-to-cover} if for any finite cover $\mathcal{Z}$ of $s\in S_1$, $\phi(\mathcal{Z})$ is a cover of $\phi(s)$.  
\end{dfn}
\begin{rmk}\label{rem cover-to-cover but not cover-to-join}
A cover-to-cover homomorphism need not be cover-to-join. For a simple example, consider the directed graphs $\Gamma_1$ and $\Gamma_2$ shown below. 
\[
\begin{tikzpicture}[
    >={Stealth[length=3mm]},
    vertex/.style={circle, fill=black, inner sep=1.5pt},
    edge/.style={->, thick},
    every node/.style={font=\large}
]

\node[vertex, label=left:\normalsize{$u$}] (u) at (0,0) {};
\node[vertex, label=above right:\normalsize{$v$}] (v) at (1.5,0.8) {};
\node[vertex, label=right:\normalsize{$w$}] (w) at (1.5,-0.8) {};

\draw[edge] (u) -- node[above left] {\normalsize{$e$}} (v);
\draw[edge] (u) -- node[below left] {\normalsize{$f$}} (w);

\node[vertex, label=left:\normalsize{$a$}] (a) at (5,0) {};
\node[vertex, label=below:\normalsize{$b$}] (b) at (6.7,0) {};
\node[vertex, label=above right:\normalsize{$c$}] (c) at (8.2,0.8) {};
\node[vertex, label=below right:\normalsize{$d$}] (d) at (8.2,-0.8) {};

\draw[edge] (a) -- node[above] {\normalsize{$x$}} (b);
\draw[edge] (b) -- node[above left] {\normalsize{$y$}} (c);
\draw[edge] (b) -- node[below left] {\normalsize{$z$}} (d);

\node[] at (-1.2,0) {\normalsize{$\Gamma_1\equiv$}};
\node[] at (3.8,0) {\normalsize{$\Gamma_2\equiv$}};

\end{tikzpicture}
\]
It can be checked that the path homomorphism sending $u\mapsto a$, $v\mapsto c$, $w\mapsto d$, $e\mapsto xy$, $f\mapsto xz$ induces a $0$-preserving cover-to-cover homomorphism $\phi:S(\Gamma_1)\longrightarrow S(\Gamma_2)$. However, $\phi$ is not cover-to-join since $\{ee^{-1},ff^{-1}\}$ is a finite cover of $u$ in $S(\Gamma_1)$ but $\phi(u)=a\neq xy(xy)^{-1} \vee xz(xz)^{-1}=\phi(ee^{-1})\vee \phi(ff^{-1})$ in $S(\Gamma_2)$.\qed  
\end{rmk}
For an abelian group $G$, suppose $G\textup{-\textbf{GrIS}}_{cc}$ denotes the subcategory of $G\textup{-\textbf{GrIS}}$, where the objects are the same and the morphisms are $0$-preserving graded homomorphisms which are also cover-to-cover. The following lemma shows that the image of the functor $\mathcal{H}\circ \mathcal{S}$ is inside this subcategory.
\begin{lem}\label{lem the induced homomorphism is cover-to-cover}
Let $p\in \Hom_{k\textup{-\textbf{Gra}}}(\Lambda_1,\Lambda_2)$. Then $\tilde{p}:S_{\Lambda_1}\longrightarrow S_{\Lambda_2}$ is cover-to-cover. Consequently, $\tilde{p}|_{S_{\Lambda_1}^h}:S_{\Lambda_1}^h\longrightarrow S_{\Lambda_2}^h$ is a graded cover-to-cover homomorphism.   
\end{lem}
\begin{proof}
Let $F\in S_{\Lambda_1}$ and $\mathcal{Z}$ a finite cover of $F$. Suppose $\emptyset\neq K\le p(F)$. Choose any $(\tau,\nu)\in K$. Then by Proposition \ref{pro natural partial order on S_Lambda}, there exist $(\lambda,\mu)\in F$ and $\beta\in \Lambda_2$ such that $\tau=p(\lambda)\beta$ and $\nu=p(\mu)\beta$. Note that $p(s(\lambda))=r(\beta)$ and so by the unique $r$-path lifting of $p$, there is a unique $\gamma\in s(\lambda)\Lambda_1$ such that $p(\gamma)=\beta$. Now $H:=\{(\lambda\gamma,\mu\gamma)\}\le \{(\lambda,\mu)\}\le F$. Since $\mathcal{Z}$ is a cover of $F$, there exists some $G\in \mathcal{Z}$ such that $H\wedge G\neq \emptyset$. Choose any $(\xi,\eta)\in H\wedge G$. Since $\tilde{p}$ is a homomorphism, we have $\{(p(\xi),p(\eta))\}\le \{(\tau,\nu)\}\le K$ and $\{(p(\xi),p(\eta))\}\le p(G)$. It follows that $K\wedge p(G)\neq \emptyset$. Hence $\tilde{p}(\mathcal{Z})$ is a cover of $p(F)$. 

Since $\tilde{p}$ is an exact homomorphism, $\tilde{p}|_{S_{\Lambda_1}^h}$ is a graded homomorphism (see Remark \ref{rem the graded core functor}). Now let $F\in S_{\Lambda_1}^h$ and $\mathcal{Z}$ is a finite cover of $F$ in $S_{\Lambda_1}^h$. By the first part, it suffices to show that $\mathcal{Z}$ is a cover of $F$ in $S_{\Lambda_1}$. Assume that $F\in (S_{\Lambda_1})_\N$. If $\emptyset \neq H\le F$, then for any $(x,y)\in H$, we have $(\lambda,\mu)\in F$ and $\tau\in \Lambda_1$ such that $x=\lambda\tau$ and $y=\mu\tau$. Then $d(x)-d(y)=d(\lambda)-d(\mu)=\N$ and so $H\in (S_{\Lambda_1})_\N\subseteq S_{\Lambda_1}^h$. What we want to show is now immediate.
\end{proof}

Suppose $G$ is an abelian group and $S_1,S_2$ are $G$-graded inverse semigroups with zeros. A graded homomorphism $\psi:S_1\longrightarrow S_2$ naturally gives a graded $\mathsf{k}$-algebra homomorphism $\overline{\psi}:\mathsf{k}_0S_1\longrightarrow \mathsf{k}_0S_2$ such that $\overline{\psi}(\overline{s})=\overline{\psi(s)}$ for all $s\in S_1$. If $\psi$ is cover-to-cover, then we observe that $\mathcal{T}_\mathsf{k}(S_1)\subseteq \Ker(\pi_2\circ \overline{\psi})$, where $\pi_2:\mathsf{k}_0S_2\longrightarrow \mathsf{k}_0S_2/{\mathcal{T}_\mathsf{k}(S_2)}$ is the canonical quotient map. Indeed, if $e\in \mathcal{E}(S_1)$ and $\mathcal{Z}$ is a finite cover of $e$, then \[(\pi_2\circ \overline{\psi})\left(\displaystyle{\prod_{f\in \mathcal{Z}}}(e-f)\right)=\pi_2\left(\displaystyle{\prod_{f\in \mathcal{Z}}}(\psi(e)-\psi(f))\right)=0\] since $\psi(\mathcal{Z})$ is a finite cover of $\psi(e)$. It follows that the map $\psi_*:\mathsf{k}_0S_1/{\mathcal{T}_\mathsf{k}(S_1)}\longrightarrow \mathsf{k}_0S_2/{\mathcal{T}_\mathsf{k}(S_2)}$ defined by $\psi_*(\pi_1(x))=\pi_2(\overline{\psi}(x))$ is well defined. 

If $G\textup{-\textbf{GrAlg}}_\mathsf{k}$ is the category of all $G$-graded $\mathsf{k}$-algebras with graded algebra homomorphisms, then the above discussion gives us a functor $\mathcal{A}:G\textup{-\textbf{GrIS}}_{cc}\longrightarrow G\textup{-\textbf{GrAlg}}_{\mathsf{k}}$ such that $\mathcal{A}(S):=\mathsf{k}_0S/{\mathcal{T}_\mathsf{k}(S)}$ for any $S\in \Obj(G\textup{-\textbf{GrIS}}_{cc})$ and $\mathcal{A}(\psi):=\psi_*$ for any $\psi\in \Hom_{G\textup{-\textbf{GrIS}}_{cc}}(S_1,S_2)$. Specializing this to $k$-graphs, now we have the chain of functors \[k\textup{-\textbf{Gra}}\overset{\mathcal{S}}{\xrightarrow{\hspace*{1.2cm}}}\mathbb{Z}^k\textup{-\textbf{GrOIS}}\overset{\mathcal{H}}{\xrightarrow{\hspace*{1.2cm}}}\mathbb{Z}^k\textup{-\textbf{GrIS}}\supset \mathbb{Z}^k\textup{-\textbf{GrIS}}_{cc}\overset{\mathcal{A}}{\xrightarrow{\hspace*{1.2cm}}} \mathbb{Z}^k\textup{-\textbf{GrAlg}}_{\mathsf{k}}\]
Now by Lemma \ref{lem the induced homomorphism is cover-to-cover}, this gives the composition 
\begin{equation}\label{eq func1}
\mathcal{A}\circ \mathcal{H}\circ \mathcal{S}:k\textup{-\textbf{Gra}}\longrightarrow \mathbb{Z}^k\textup{-\textbf{GrAlg}}_{\mathsf{k}}
\end{equation}
such that $(\mathcal{A}\circ \mathcal{H}\circ \mathcal{S})(\Lambda)=\mathsf{k}_0S_\Lambda^h/{\mathcal{T}_\mathsf{k}(S_\Lambda^h)}$ for any $\Lambda\in \Obj(k\textup{-\textbf{Gra}})$ and $(\mathcal{A}\circ \mathcal{H}\circ \mathcal{S})(p)=\tilde{p}_*$ for any $p\in \Hom_{k\textup{-\textbf{Gra}}}(\Lambda_1,\Lambda_2)$. Note that for simplicity, we are writing the restriction of $\tilde{p}:S_{\Lambda_1}\longrightarrow S_{\Lambda_2}$ to the graded core $S_{\Lambda_1}^h$ also by $\tilde{p}$.  

The following proposition gives another functor between the categories in (\ref{eq func1}).
\begin{thm}\label{th KP functor}
Let $\Lambda_1,\Lambda_2$ be finitely aligned $k$-graphs and $p\in \Hom_{k\textup{-\textbf{Gra}}}(\Lambda_1,\Lambda_2)$. Let $\{s_\lambda,s_{\mu^*}~|~\lambda,\mu\in \Lambda_1\}$ and $\{t_\alpha,t_{\beta^*}~|~\alpha,\beta\in \Lambda_2\}$ be the universal Kumjian--Pask families in $\KP_\mathsf{k}(\Lambda_1)$ and $\KP_\mathsf{k}(\Lambda_2)$ respectively. Then 

$(i)$ there is a graded injective $k$-algebra homomorphism $\widehat{p}:\KP_\mathsf{k}(\Lambda_1)\longrightarrow \KP_\mathsf{k}(\Lambda_2)$ such that $\widehat{p}(s_\lambda)=t_{p(\lambda)}$ and $\widehat{p}(s_{\mu^*})=t_{p(\mu)^*}$ for all $\lambda,\mu\in \Lambda_1$;

$(ii)$ the correspondence 
\begin{equation}\label{eq func2}
\mathcal{K}\mathcal{P}:k\textup{-\textbf{Gra}}\longrightarrow \mathbb{Z}^k\textup{-\textbf{GrAlg}}_{\mathsf{k}}
\end{equation}
with $\mathcal{K}\mathcal{P}(\Lambda):=\KP_\mathsf{k}(\Lambda)$ for all $\Lambda\in \Obj(k\textup{-\textbf{Gra}})$ and $\mathcal{K}\mathcal{P}(p)=\widehat{p}$ for all $p\in \Hom_{k\textup{-\textbf{Gra}}}(\Lambda_1,\Lambda_2)$ is a faithful functor. 
\end{thm}
\begin{proof}
$(i)$ Define $T_\lambda:=t_{p(\lambda)}$ and $T_{\mu^*}:=t_{p(\mu)^*}$ for all $\lambda,\mu\in \Lambda_1$. It suffices to show that $\{T_\lambda,T_{\mu^*}\}$ is a Kumjian--Pask $\Lambda_1$ family in $\KP_\mathsf{k}(\Lambda_2)$. The verifications for (KP1) and (KP2) are straightforward and follow from the fact that $p$ is an injective $k$-graph morphism. For (KP3), let $\lambda,\mu\in \Lambda_1$. Consider the map 
\begin{align*}
\rho:\Lambda_1^{\min}(\lambda,\mu)&\longrightarrow \Lambda_2^{\min}(p(\lambda),p(\mu))\\
(\alpha,\beta)&\longmapsto (p(\alpha),p(\beta)).
\end{align*}
Using the unique $r$-path lifting property of $p$, it can be shown that $\rho$ is a bijection (see the proof of Theorem \ref{th answer to the query}). Now \[T_{\lambda^*}T_\mu=t_{p(\lambda)^*}t_{p(\mu)}=\displaystyle{\sum_{(\tau,\eta)\in \Lambda_2^{\min}(p(\lambda),p(\mu))}}t_\tau t_{\eta^*}=\displaystyle{\sum_{(\alpha,\beta)\in \Lambda_1^{\min}(\lambda,\mu)}}t_{p(\alpha)}t_{p(\beta)^*}=\displaystyle{\sum_{(\alpha,\beta)\in \Lambda_1^{\min}(\lambda,\mu)}}T_\alpha T_{\beta^*}.\] Hence (KP3) is verified. Now for (KP4), let $v\in \Lambda_1$ and $X$ a finite exhaustive subset of $v\Lambda_1$. We claim that $p(X)$ is an exhaustive subset of $p(v)\Lambda_2$. Let $\gamma\in p(v)\Lambda_2$. Then by the unique $r$-path lifting, there is a unique $\lambda\in v\Lambda_1$ such that $p(\lambda)=\gamma$. Since $X$ is exhaustive subset of $v\Lambda_1$, there exists $\tau\in X$ such that $\Lambda_1^{\min}(\tau,\lambda)\neq \emptyset$. In view of the map $\rho$ defined above, $\Lambda_2^{\min}(p(\tau),\gamma)\neq \emptyset$. This proves our claim. Now using this, we have \[\displaystyle{\prod_{\tau\in X}}\left(T_v-T_\tau T_{\tau^*}\right)=\displaystyle{\prod_{\tau\in X}}\left(t_{p(v)}-t_{p(\tau)}t_{p(\tau)^*}\right)=\displaystyle{\prod_{\eta\in p(X)}}\left(t_{p(v)}-t_\eta t_{\eta^*}\right)=0.\] Hence we have (KP4). The result now follows from the universal property of $\KP_\mathsf{k}(\Lambda_1)$. Since $\widehat{p}$ preserves the degrees of generators, it is a graded homomorphism. Moreover, since $\widehat{p}(r s_v)=rt_{p(v)}\neq 0$ for all $r\in \mathsf{k}\setminus \{0\}$ and $v\in \Lambda_1^0$, it follows from \cite[Theorem 3.7(b) \& Theorem 4.1]{Clark} that $\widehat{p}$ is injective.

$(ii)$ From the definition of $\widehat{p}$, it follows that $\mathcal{K}\mathcal{P}$ takes the identity morphism on any $\Lambda$ to the identity morphism on $\KP_\mathsf{k}(\Lambda)$, and preserves composition of $k$-graph morphisms. Hence $\mathcal{K}\mathcal{P}$ is a functor. Now suppose $p,q\in \Hom_{k\textup{-\textbf{Gra}}}(\Lambda_1,\Lambda_2)$ are such that $\widehat{p}=\widehat{q}$. Then for any $\lambda\in \Lambda_1$, $t_{p(\lambda)}=\widehat{p}(s_\lambda)=\widehat{q}(s_\lambda)=t_{q(\lambda)}$. Multiplying by $t_{q(\lambda)^*}$, we get $t_{q(\lambda)^*}t_{p(\lambda)}=t_{s(q(\lambda))}\neq 0$. Since $d_2(p(\lambda))=d_2(q(\lambda))=d_1(\lambda)$, by (KP3), it must happen that $p(\lambda)=q(\lambda)$. Hence $p=q$. This shows that $\mathcal{K}\mathcal{P}$ is faithful.  
\end{proof}

Now we have the perfect categorical set-up to interpret the isomorphism established as Theorem \ref{th the relationship between S_Lambda and KPA}. For any finitely aligned $k$-graph $\Lambda$, let us denote the isomorphism of Theorem \ref{th the relationship between S_Lambda and KPA} by $\Theta_\Lambda$. 
\begin{rmk}\label{rem the natural isomorphism}
Suppose $p\in \Hom_{k\textup{-\textbf{Gra}}}(\Lambda_1,\Lambda_2)$. Then for any $F\in S_{\Lambda_1}^h$, using the notations set throughout this section, we have \[\Theta_{\Lambda_2}\left(\tilde{p}_*\left(\pi_1\left(\overline{F}\right)\right)\right)=\Theta_{\Lambda_2}\left(\pi_2\left(\overline{p(F)}\right)\right)=\displaystyle{\sum_{(\lambda,\mu)\in F}}t_{p(\lambda)}t_{p(\mu)^*}.\] On the other hand, \[\widehat{p}\left(\Theta_{\Lambda_1}\left(\pi_1\left(\overline{F}\right)\right)\right)=\widehat{p}\left(\displaystyle{\sum_{(\lambda,\mu)\in F}}s_\lambda s_{\mu^*}\right)=\displaystyle{\sum_{(\lambda,\mu)\in F}}t_{p(\lambda)}t_{p(\mu)^*}.\] Hence the square shown in the following diagram is commutative.
\[
\begin{tikzcd}[
  column sep=2cm,
  row sep=1.7cm
]
\Lambda_1
  \arrow[d, "p"']
&
\mathsf{k}_0S_{\Lambda_1}^h/{\mathcal{T}_\mathsf{k}(S_{\Lambda_1}^h)}
  \arrow[d, "\tilde{p}_*"']
  \arrow[r, "\Theta_{\Lambda_1}"]
&
\KP_\mathsf{k}(\Lambda_1)
  \arrow[d, "\widehat{p}"]
\\
\Lambda_2
&
\mathsf{k}_0S_{\Lambda_2}^h/{\mathcal{T}_\mathsf{k}(S_{\Lambda_2}^h)}
  \arrow[r, "\Theta_{\Lambda_2}"']
&
\KP_\mathsf{k}(\Lambda_2)
\end{tikzcd}
\]
By Theorem \ref{th the relationship between S_Lambda and KPA}, each $\Theta_\Lambda$ is an isomorphism in the category $\mathbb{Z}^k\textup{-\textbf{GrAlg}}_\mathsf{k}$. It follows that the functors $\mathcal{A}\circ \mathcal{H}\circ \mathcal{S}$ and $\mathcal{K}\mathcal{P}$ (see (\ref{eq func1}) and (\ref{eq func2})) are naturally isomorphic.\qed 
\end{rmk}

\section{A few topics for further study}\label{sec future scopes}
We end this paper by mentioning some probable future lines of investigation, which originate from this work.

\begin{itemize}
\item[$\bullet$] In Section \ref{sec Relation to KP-algebra}, the notion of a graded-ordered semigroup has been introduced. The Farthing--Muhly--Yeend inverse semigroup is a nice example of such a semigroup. As the concept perfectly blends order theory, group grading on semigroups and semirings (see Examples \ref{ex graded-ordered semigroup}), an independent study of graded-ordered semigroups might be helpful to unravel the interconnections among these theories.

\item[$\bullet$] The graded core $S_\Lambda^h$ is naturally $\mathbb{Z}^k$-graded as described in this paper. It is interesting to find when $S_\Lambda^h$ is strongly graded, locally strongly graded, saturated strongly graded (see \cite{HazMes} for the definitions) under this $\mathbb{Z}^k$-grading. In \cite{CHR}, Clark, Hazrat and Rigby characterized strongly graded Kumjian--Pask algebras of finitely aligned $k$-graphs. It was shown that $\KP_\mathsf{k}(\Lambda)$ is strongly $\mathbb{Z}^k$-graded with respect to the natural $\mathbb{Z}^k$-grading if and only if $\Lambda$ is row-finite without sources and satisfies a certain geometric condition (known as Condition (Y)). In the level of directed graph, Hazrat and Mesyan proved that Condition (Y) is equivalent to saying that the graph inverse semigroup is saturated strongly graded with respect to the natural $\mathbb{Z}$-grading (see \cite[Lemma 8.11]{HazMes}). In view of the isomorphism of Theorem \ref{th the relationship between S_Lambda and KPA}, it is natural to ask if an analogous result holds for the graded core of the Farthing--Muhly--Yeend inverse semigroup.

\item[$\bullet$] Given a directed graph $\Gamma$, Kudryavtseva in \cite{Ganna} established that the graph monoid of $M_\Gamma$ can be realized as the type monoid of $\mathsf{B}_{tight}(S(\Gamma))$. In \cite{HMPS}, the $k$-graph monoid $M_\Lambda$ of a row-finite $k$-graph $\Lambda$ has been defined. Given the description of the Cuntz--Krieger ideal of $\mathsf{B}(S_\Lambda)$ in Theorem \ref{th description of Cuntz-Krieger ideal}, can we obtain the isomorphism $\Typ(\mathsf{B}_{tight}(S_\Lambda))\cong M_\Lambda$? 
\end{itemize}

\section*{Acknowledgements} I sincerely thank my doctoral advisor Prof. Sujit Kumar Sardar for his continuous support throughout my work. I also thank Prof. Roozbeh Hazrat for being an active source of inspiration and for all the fruitful discussions we have.

\end{document}